\documentclass[a4paper,12pt,twoside]{article}
\usepackage[pdftex,colorlinks,bookmarksnumbered,bookmarksopen,citecolor=red,urlcolor=red]{hyperref}
\usepackage{graphicx}\usepackage{a4wide}
\usepackage{amssymb,amsmath,amsthm}
\usepackage[only,llbracket,rrbracket,interleave]{stmaryrd}
\usepackage[pdftex,colorlinks,bookmarksnumbered,bookmarksopen,citecolor=red,urlcolor=red]{hyperref}
\usepackage{subfigure,color,colordvi,graphicx,xcolor}
\usepackage[only,llbracket,rrbracket,interleave]{stmaryrd}
\usepackage[sort,compress]{cite}
\usepackage{enumitem}
\usepackage{sidecap, caption}
\usepackage{multirow}

\usepackage{setspace}

\makeatletter
\newcommand*\wt[2][0.2ex]{%
        \begingroup
        \mathchoice{\wt@helper{#1}{#2}{\displaystyle}{\textfont}}
                   {\wt@helper{#1}{#2}{\textstyle}{\textfont}}
                   {\wt@helper{#1}{#2}{\scriptstyle}{\scriptfont}}
                   {\wt@helper{#1}{#2}{\scriptscriptstyle}{\scriptscriptfont}}%
        \endgroup
        #2%
}
\newcommand*\wt@helper[4]{%
        \def\currentfont{\the#41}%
        \def\currentskewchar{\char\the\skewchar\currentfont}%
        \setbox\tw@\hbox{\currentfont#2\currentskewchar}%
        \dimen@ii\wd\tw@
        \setbox\tw@\hbox{\currentfont#2{}\currentskewchar}%
        \advance\dimen@ii-\wd\tw@
        \rlap{\raisebox{-#1}{$\m@th#3\kern\dimen@ii\widetilde{\phantom{#2}}$}}%
}
\makeatother

\usepackage{booktabs,array,dcolumn}
\newcommand{\PreserveBackslash}[1]{\let\temp=\\#1\let\\=\temp}
\newcolumntype{C}[1]{>{\PreserveBackslash\centering}p{#1}}
\newcolumntype{R}[1]{>{\PreserveBackslash\raggedleft}p{#1}}
\newcolumntype{L}[1]{>{\PreserveBackslash\raggedright}p{#1}}

\usepackage{fancyhdr}\usepackage{lastpage} 

\newenvironment{keyword}{\begin{quote}\emph{\textbf{Keywords:}}}{\end{quote}}
\newtheorem{remark}{Remark}

\newcommand{\hl}[1]{{\color{black} #1\color{black}}}

\newcommand{\bm}[1]{\text{\boldmath $#1$\unboldmath}}
\newcommand{\abs}[1]{\lvert#1\rvert}
\newcommand{\norm}[1]{\lVert#1\rVert}
\DeclareMathOperator{\adj}{adj}
\DeclareMathOperator{\Det}{det}

\newcommand{\Ga}[1]{\Gamma_{\!\!#1}}
\newcommand{\GaEta}[1]{\Gamma_{\!\!#1}^{\eta}}
\newcommand{\OmEta}{\Omega^{\eta}}
\newcommand{\Iset}{\bm{\mathcal{I}}}
\newcommand{\I}{\mathcal{I}}

\newcommand{\Aset}{\mathcal{A}}
\newcommand{\Eset}{\mathcal{E}}

\newcommand{\mat}[1]{\mathbf{#1}}

\newcommand{\grad}{\bm{\nabla}}

\newcommand{\RR}{\mathbb{R}}

\newcommand{\D}{\mathrm{D}}
\newcommand{\Vbar}{\bar{\mathrm{V}}}
\newcommand{\Gt}{\mat{G}_{\theta}}
\newcommand{\Vspace}{\mathcal{V}}
\newcommand{\Wspace}{\mathcal{W}}

\newcommand{\bgEta}{\bm{g}^{\eta}}

\newcommand{\bnEta}{\bm{n}^{\!\eta}}

\newcommand{\bu}{\bm{u}}
\newcommand{\buEta}{\bm{u}^{\!\eta}}

\newcommand{\buS}{\bm{u}^\star}
\newcommand{\bv}{\bm{v}}

\newcommand{\bx}{\bm{x}}

\newcommand{\beps}{\bm{\varepsilon}}
\newcommand{\bEta}{\bm{\eta}}

\newcommand{\bsigma}{\bm{\sigma}}
\newcommand{\bsigmaEta}{\bm{\sigma}^{\eta}}

\newcommand{\bsigmaS}{\bm{\sigma}^\star}

\newcommand{\bxi}{\bm{\xi}}

\newcommand{\eltwo}{\ensuremath{\mathcal{L}_2}}

\newcommand{\nsd}  {\ensuremath{\texttt{n}_{\texttt{sd}}}}
\newcommand{\npa}  {\ensuremath{\texttt{n}_{\texttt{pa}}}}
\newcommand{\nS}  {\ensuremath{\texttt{n}_{\texttt{S}}}}
\newcommand{\nT}  {\ensuremath{\texttt{n}_{\theta}}}
\newcommand{\nNNpar}[1][]{\ensuremath{\texttt{n}_{\omega#1}}}
\newcommand{\nA}  {\ensuremath{\texttt{n}_{\alpha}}}
\newcommand{\msd}  {\ensuremath{\texttt{m}_{\texttt{sd}}}}
\newcommand{\numel}{\ensuremath{\texttt{n}_{\texttt{el}}}}

\newcommand{\niter}[1][]{\ensuremath{\texttt{n}_{\texttt{iter}}^{#1}}}

\newcommand{\sInit}[1]{\ensuremath{\texttt{s}_{#1}}}

\newcommand{\Uad}{\mathcal{U}_{\text{ad}}}
\newcommand{\UadS}{\mathcal{U}_{\text{ad}}^{\star}}
\newcommand{\Js}{J^{\star}}

\newcommand{\Insd}{\mat{I}_{\nsd}}

\newcommand{\matC}{\mat{C}}
\newcommand{\matA}{\mat{A}}
\newcommand{\matAs}{\mat{A}^{\!\!\star}}
\newcommand{\matAc}[1]{\mat{A}^{\!\!c}_{#1}}
\newcommand{\matAhs}{\mat{A}_{\texttt{HS}}}

\newcommand{\lamSig}[1]{\lambda^{\!\sigma}_{#1}}
\newcommand{\vSig}[1]{\mathbf{v}^{\sigma}_{\! #1}}

\newcommand{\lambdaE}{\lambda_{\texttt{eff}}}
\newcommand{\muE}{\mu_{\texttt{eff}}}

\newcommand{\thetaT}{\tilde{\theta}}

\newcommand{\thetaBar}{\overline{\theta}}
\newcommand{\xiT}{\tilde{\xi}}

\newcommand{\NNpar}[1]{\omega_{#1}}

\newcommand{\Loss}[1]{\mathcal{L}_{#1}}
\newcommand{\LAE}{\Loss{\texttt{AE}}}
\newcommand{\LEta}{\Loss{\eta}}
\newcommand{\LCode}{\Loss{\alpha}}
\newcommand{\Reg}{\mathcal{R}}

\newcommand{\thetaEta}{\theta^{\eta}}
\newcommand{\btheta}{\bm{\theta}}
\newcommand{\bthetaIn}{\btheta^{\text{in}}}
\newcommand{\bthetaAE}{\btheta^{\texttt{AE}}}
\newcommand{\bthetaEta}{\btheta^{\eta}}
\newcommand{\bthetaRef}{\btheta^{\text{ref}}}
\newcommand{\balpha}{\bm{\alpha}}
\newcommand{\balphaAE}{\balpha^{\texttt{AE}}}
\newcommand{\balphaEta}{\balpha^{\eta}}

\newcommand{\Jref}[1][]{J^{\text{ref}}_{#1}}
\newcommand{\Jeta}[1][]{J^{\eta}_{#1}}
\newcommand{\Vref}[1][]{V^{\text{ref}}_{#1}}
\newcommand{\Veta}[1][]{V^{\eta}_{#1}}

\newcommand{\nTrain}  {\ensuremath{\texttt{n}_{\texttt{train}}}}
\newcommand{\nVal}  {\ensuremath{\texttt{n}_{\texttt{val}}}}
\newcommand{\nTest}  {\ensuremath{\texttt{n}_{\texttt{test}}}}
\newcommand{\TTrain}  {\ensuremath{\texttt{T}_{\texttt{train}}}}
\newcommand{\TTrainAv}  {\ensuremath{\overline{\texttt{T}}_{\texttt{train}}}}

\newcommand{\Err}[1]{\mathrm{E}^{\theta}_{#1}}

\newcommand{\elemE}{\Omega_e}

\newcommand{\stressV}{\bm{\sigma}_{\texttt{V}}}
\newcommand{\strainV}{\bm{\varepsilon}_{\texttt{V}}}

\RequirePackage{algorithm}[0.1] 
\usepackage{algorithmic} 

\graphicspath{{figsPDF/}{figsPNG/}{figs/}{figs_arXiv/}}

\usepackage{hyphenat}
\begin{document}
\fancypagestyle{plain}{%
  \renewcommand{\headrulewidth}{0pt}%
  \fancyhead[L]{}%
  \fancyfoot[C]{\footnotesize Page \thepage\ of \pageref{LastPage}}%
}
\title{A surrogate model for topology optimisation of elastic structures via parametric autoencoders}
\author{
\renewcommand{\thefootnote}{\arabic{footnote}}
			Matteo Giacomini\footnotemark[1]\textsuperscript{ \ ,}\footnotemark[2]
                         \hspace{1ex}  and
			Antonio Huerta\footnotemark[1]\textsuperscript{ \ ,}\footnotemark[2]\textsuperscript{ \ ,}*
}


\date{\today}
\maketitle

\renewcommand{\thefootnote}{\arabic{footnote}}

\footnotetext[1]{Laboratori de C\`alcul Num\`eric (LaC\`aN), ETS de Ingenier\'ia de Caminos, Canales y Puertos, Universitat Polit\`ecnica de Catalunya - BarcelonaTech (UPC), Barcelona, Spain.}
\footnotetext[2]{Centre Internacional de M\`etodes Num\`erics en Enginyeria (CIMNE), Barcelona, Spain.
\vspace{5pt}\\
* Corresponding author: Antonio Huerta. \textit{E-mail:} \texttt{antonio.huerta@upc.edu}
}

\begin{abstract}
A surrogate-based topology optimisation algorithm for linear elastic structures under parametric loads and boundary conditions is proposed.
Instead of learning the parametric solution of the state (and adjoint) problems or the optimisation trajectory as a function of the iterations,  the proposed approach devises a surrogate version of the entire optimisation pipeline. 
First, the method predicts a \emph{quasi-optimal} topology for a given problem configuration as a surrogate model of high-fidelity topologies optimised with the homogenisation method. 
This is achieved by means of a feed-forward net learning the mapping between the input parameters characterising the system setup and a latent space determined by encoder/decoder blocks reducing the dimensionality of the parametric topology optimisation problem and reconstructing a high-dimensional representation of the topology.
Then,  the predicted topology is used as an \emph{educated} initial guess for a computationally efficient algorithm penalising the intermediate values of the design variable, while enforcing the governing equations of the system.
This step allows the method to correct potential errors introduced by the surrogate model, eliminate artifacts,  and refine the design in order to produce topologies consistent with the underlying physics.
Different architectures are proposed and the approximation and generalisation capabilities of the resulting models are numerically evaluated.
The \emph{quasi-optimal} topologies allow to outperform the high-fidelity optimiser by reducing the average number of optimisation iterations by $53\%$ while achieving discrepancies below $4\%$ in the optimal value of the objective functional, even in the challenging scenario of testing the model to extrapolate beyond the training and validation domain.
\end{abstract}

\section{Introduction}
\label{sc:Intro}

Reduced order models (ROMs)~\cite{AH-CHRW:17} are commonly employed to solve many-queries problems arising in outer-loop applications such as optimisation pipelines.
For instance, shape optimisation problems have been extensively tackled by means of ROMs, devising surrogate models of geometrically parametrised problems~\cite{Rozza-LR-10,SZ-ZDMH:15,sevilla2020solution,RS-SBGH-20,Rozza-SZR-20} and accelerating the resulting optimisation procedures by replacing the computationally-demanding full-order solver with a parsimonious reduced model~\cite{Rozza-MQR-12,AH-AHCCL:14,Farhat-ZF-15}.
Nonetheless, the above mentioned strategies rely on the assumption of being able to characterise shape deformations by means of a parametric transformation preserving the domain connectivity, thus limiting their applicability to fixed topologies.

Research on ROM-accelerated topology optimisation (see, e.g., ~\cite{Sigmund-LSJA-21}) has received less attention in the literature given the intrinsic extreme high-dimensionality of the problem (potentially, one extra dimension for each cell employed in the spatial discretisation).
To circumvent this issue and avoid the need to precompute and process a large dataset of snapshots,  \emph{on-the-fly} solutions based on proper orthogonal decomposition (POD) have been presented in~\cite{Gogu-15,Breitkopf-XLBRDZ-20,Breitkopf-XLDRZD-20}. These methods construct a reduced basis during optimisation and adaptively enrich it if the solution of the reduced system is inaccurate.
Alternative strategies propose to accelerate the solid isotropic material with penalisation (SIMP) optimiser by providing an informed initial guess, either constructed by means of POD on a fixed coarse mesh~\cite{Perotto-FMP-19} or identified as the \emph{closest} configuration within a dictionary of previously generated snapshots~\cite{Rodenas-MMNACRN-25}.
All these approaches share a common shortcoming, that is, the limited generalisation capabilities of the described surrogate models when extrapolation to unseen configurations and outside the domain of training is performed.

Machine learning (ML) and deep learning (DL) techniques provide promising solutions to this challenge, see the recent reviews~\cite{Sigmund-WABS-22,Sin-SSK-23}.
A methodology mimicking the rationale of substituting the state (and possibly adjoint) full-order solver by means of a surrogate version has been discussed in~\cite{Jeong-JBBRZG-23,Jeong-JBBXRZG-23} using physics-informed neural networks to replace the finite element computation.
Neural networks (NNs) have alternatively been exploited to introduce novel parametrisations of the design field.  In~\cite{Zhang-ZLZCYZ-21,Sigmund-HSLVK-24},  the design variable is described by means of a neural network mapping the spatial coordinates to material densities and sensitivities are computed using automatic differentiation, whereas \cite{To-DT-21} proposes to implicitly define the level-set function using a NN and update its weights and biases instead of solving the Hamilton-Jacobi equation. Analogously,  \cite{Chandrasekhar-CS-21} employs a NN to achieve a density field independent of the underlying computational mesh but relies on a finite element solver to compute sensitivities.
Similarly to the framework in~\cite{Gogu-15,Breitkopf-XLBRDZ-20,Breitkopf-XLDRZD-20}, ML approaches exploiting a hierarchy of coarse and fine grids with convolutional neural nets and \emph{on-the-fly} training during optimisation were proposed in~\cite{Paulino-ACZMTP-22,Paulino-CZTMDSP-21,Zhang-ZJLXGWZ-24}.
An alternative approach relies on learning \emph{offline} the evolution of the optimisation trajectory in order to bypass the intermediate iterations during the \emph{online} phase. This has been achieved using deep belief networks~\cite{Lagaros-KKL-20},  convolutional NNs with U-Net architecture for the spatial features coupled with recurrent NNs with long short-term memory units for the temporal evolution~\cite{Qiu-QDY-21}, and proximal policy optimisation to train a reinforcement learning agent to sequentially discover novel designs~\cite{Jang-JYK-22}.
Finally, a hybrid approach combining tensor decomposition inspired by a posteriori proper generalised decomposition~\cite{Giacomini-GBSH-21} to reduce the dimensionality of the problem with NNs has been introduced in~\cite{Liu-LKLPGMLCAL-23}.

This work aims to devise a non-intrusive surrogate model of the optimised layout for a class of topology optimisation problems with parametrised loads and boundary conditions.
In contrast to the previously mentioned approaches that replace the state and adjoint solutions by means of reduced models or learn the progression of the optimisation iterations, the proposed method focuses on \emph{reducing} the computational burden of the entire optimisation process through a surrogate model to generate a \emph{quasi-optimal} topology.

\hl{
The motivation of this study stems from the large amount of data owned by industries concerning optimised topologies.
This information can be naturally employed to train surrogate models, without the need to first execute the corresponding expensive finite element analyses to populate datasets. Moreover, using optimised topologies as data (instead of the solutions of the state and adjoint equations) significantly reduces storage requirements since only the final, optimised configurations need to be saved.
This article proposes a demonstrator of the feasibility of \emph{learning} fundamental features of optimised topologies from existing data.
The objective is to devise a data-driven initial condition to accelerate the physics-based high-fidelity optimisation algorithm when new configurations of the system, associated with parameters not seen during training, are studied.
To this end, a synthetic benchmark problem consisting of a two-dimensional linear elastic cantilever beam with two parameters controlling the external load is considered.
}

\hl{
The proposed surrogate model trains a feed-forward NN to map the input parameters into a low-dimensional latent space, whereas a decoder-type module is in charge of reconstructing a high-dimensional representation of the topology.
}
The model is trained starting from a dataset of optimised topologies computed using the homogenisation method. The \emph{quasi-optimal} topology provided as outcome is then employed as an \emph{educated} initial guess for a novel algorithm penalising intermediate density values, while significantly reducing the overall computational cost of the optimisation procedure.
The proposed approach shares similarities with existing works in the literature learning the relation between problem setup and optimised topology while compressing input information: \cite{Yu-YHJJ-19} trains an autoencoder coupled with a generative adversarial network,  \cite{Wang-WXPCZZ-22} relies on a convolutional NN with U-Net architecture, whereas \cite{Corigliano-GTBC-24} employs support vector regression and k-nearest neighbours.
All these strategies perform predictions in a single step, with no mechanism to correct errors, eliminate artifacts or refine the design, contrary to the novel surrogate-based algorithm that employs the \emph{quasi-optimal} solution as initial guess to improve the topology description.
Moreover,  these methods are purely data-driven and do not enforce the underlying governing equations in the \emph{online} phase, possibly leading to invalid topologies or designs not physically consistent.
This is avoided in this work by means of the surrogate-based algorithm that solves the linear elastic equation while penalising intermediate material distributions.
\hl{
In this context, the proposed methodology does not aim to substitute the high-fidelity optimiser but to complement it, accelerating the involved iterative procedure.
}

The remainder of this article is structured as follows.
Section~\ref{sc:Problem} introduces the parametric topology optimisation problem under analysis and its homogenised version.
Section~\ref{sc:DataHomog} describes the dataset of optimised topologies and the high-fidelity optimisation algorithm employed to construct it.
\hl{
The design of a linear elastic structure in two spatial dimensions with two parameters is considered to construct a demonstrator of the surrogate model.
}
The techniques to generate the \emph{quasi-optimal} topologies via surrogate modelling are presented in Section~\ref{sc:Surrogate}, together with the novel surrogate-based optimisation algorithm.
Section~\ref{sc:Simulations} numerically validates the proposed surrogate models, in both interpolation and extrapolation scenarios, and showcases the accuracy and computational efficiency of the surrogate-based optimisation strategy, when applied to cases not seen during training.
Finally, Section~\ref{sc:Conclusion} summarises the contributions of the work and two appendices provide technical details on the computation of the homogenised elasticity tensor in the high-fidelity optimisation algorithm and on the architecture of the networks composing the surrogate models.

\section{Parametric topology optimisation of structures}
\label{sc:Problem}

Consider a fixed, open, bounded domain $\D \subset \RR^{\nsd}$ in $\nsd$ spatial dimensions.  The set $\OmEta \subset \D$ describes the reference configuration of an elastic body such that $\partial\OmEta = \GaEta{D} \cup \GaEta{N} \cup \GaEta{0}$,  with $\GaEta{D} \neq \emptyset$, and the three boundary portions being disjoint by pairs.  For the sake of simplicity, $\GaEta{D}$ and $\GaEta{N}$ are henceforth assumed to be also subsets of $\partial\D = \GaEta{D} \cup \GaEta{N} \cup \GaEta{}$, with $\GaEta{} \cap \GaEta{D} = \emptyset$, and $\GaEta{} \cap \GaEta{N} = \emptyset$.  Moreover, the tuple $\bEta=(\eta_1,\ldots,\eta_{\npa})^\top$ of $\npa$ parameters, defined in the set $\Iset = \I_1 \times \cdots \times \I_{\npa}$, is introduced to control the working conditions of the system. 


The system of equations describing the mechanical behaviour of the homogeneous, isotropic, linear elastic structure $\OmEta$ is given by
\begin{equation}\label{eq:elasticity}
\left\{
\begin{aligned}
- \grad \cdot \bsigmaEta &= \bm{0} && \text{in} \ \OmEta \times \Iset, \\
\bsigmaEta &= \matA\beps(\buEta) && \text{in} \ \OmEta \times \Iset, \\
\buEta &= \bm{0} && \text{on} \ \GaEta{D} \times \Iset, \\
\bsigmaEta\bnEta &= \bgEta && \text{on} \ \GaEta{N} \times \Iset, \\
\bsigmaEta\bnEta &= \bm{0} && \text{on} \ \GaEta{0} \times \Iset, 
\end{aligned}
\right.
\end{equation}
where $\buEta = \bu(\bx,\bEta)$ and $\bsigmaEta  = \bsigma(\bx,\bEta)$ are the displacement field and the Cauchy stress tensor, and $\beps(\buEta) := \left( \nabla \buEta + [\nabla \buEta]^\top \right)/2$ is the linearised strain tensor. The fourth-order elasticity tensor $\matA$ describes the Hooke's law
\begin{equation}\label{eq:hooke_lame} 
\matA \beps(\buEta) = 2\mu \beps(\buEta) + \lambda (\grad \cdot \buEta)\Insd ,
\end{equation}
with Lam\'e coefficients $\lambda := E\nu[(1+\nu)(1-2\nu)]^{-1}$ and $\mu := E[2(1+\nu)]^{-1}$,  $E$ being the Young's modulus and $\nu$ the Poisson's ratio.

The structure is clamped on $\GaEta{D}$ and subject to a load $\bgEta$ on $\GaEta{N}$,  whereas a free-boundary condition is applied on $\GaEta{0}$. 
It is worth noticing that the position of the boundary portions $\GaEta{D}$, $\GaEta{N}$, and $\GaEta{0}$, as well as the boundary load $\bgEta$ depend upon the parameters $\bEta$, that consequently affect the physical fields $\buEta$ and $\bsigmaEta$ fulfilling equation~\eqref{eq:elasticity}.

The parametric topology optimisation of elastic structures consists of minimising, for any value $\bEta \in \Iset$ of the parameters, an objective functional $J$, e.g., the compliance
\begin{equation}\label{eq:compliance} 
J(\OmEta) = \int_{\GaEta{N}} \bgEta \cdot \buEta \ d\Ga{} ,
\end{equation}
under a constraint on the maximum volume fraction $\Vbar$ occupied by the structure $\OmEta$ in $\D$, that is,
\begin{equation}\label{eq:topOpt} 
\min_{\OmEta \in \Uad} J(\OmEta) , \ \forall \bEta \in \Iset ,
\end{equation}
where the set of admissible topologies is defined as $\Uad := \{ \OmEta \subset \D \text{ s.t. } \abs{\OmEta} \leq \Vbar\abs{\D} \}$.

\subsection{Parametric homogenised problem}
\label{sc:HomogenPb}


For a fixed set of parameters $\bEta \in \Iset$, it is well known that problem~\eqref{eq:topOpt} usually does not admit any minimiser because the performance of composite structures with small microstructures is always superior to homogeneous structures consisting of a plain material~\cite{Allaire-book12}.
In this context, the homogenisation method~\cite{Kikuchi-BK-88} aims to extend the framework of problem~\eqref{eq:topOpt} to account for composite shapes described by the local material density (or local solid volume fraction) $\theta=\theta(\bx,\bEta)$ and by the homogenised elasticity tensor $\matAs=\matAs(\bx,\bEta)$, encapsulating the information on the microstructure at any point $\bx \in \D$,  and for all $\bEta \in \Iset$.
It follows that the volume fraction occupied by $\OmEta$ in $\D$ can be defined in terms of the material density as
\begin{equation}\label{eq:volFrac} 
V(\theta) = \frac{1}{\abs{\D}}\int_{\D} \theta \ d\bx .
\end{equation}

The macroscopic displacement $\buS(\bx,\bEta)$ and stress tensor $\bsigmaS(\bx,\bEta)$ of the structure are obtained by solving, in the entire domain $\D$, the system
\begin{equation}\label{eq:elasticityH}
\left\{
\begin{aligned}
- \grad \cdot \bsigmaS &= \bm{0} && \text{in} \ \D \times \Iset, \\
\bsigmaS &= \matAs\beps(\buS) && \text{in} \ \D \times \Iset, \\
\buS &= \bm{0} && \text{on} \ \GaEta{D} \times \Iset, \\
\bsigmaS\bnEta &= \bgEta && \text{on} \ \GaEta{N} \times \Iset, \\
\bsigmaS\bnEta &= \bm{0} && \text{on} \ \GaEta{} \times \Iset.
\end{aligned}
\right.
\end{equation}

The parametric topology optimisation problem is thus rewritten in terms of the material density $\theta :  \D \times \Iset \rightarrow [0,1]$ and the homogenised elasticity tensor $\matAs$ as
\begin{equation}\label{eq:topOptH} 
\min_{(\theta,\matAs) \in \UadS} \Js(\theta,\matAs) ,  \ \forall \bEta \in \Iset ,
\end{equation}
where $\Js$ is a relaxed objective functional accounting for the homogenised effect of the microscopic structure, namely,
\begin{equation}\label{eq:complianceHomog} 
\Js(\theta,\matAs) = \int_{\GaEta{N}} \bgEta \cdot \buS \ d\Ga{} = \int_{\D} [\matAs]^{-1} \bsigmaS : \bsigmaS \ d\bx ,
\end{equation}
for the case of the compliance, and the corresponding set of admissible topologies is given by $\UadS := \{ (\theta,\matAs) \in [0,1] \times \Gt \text{ s.t. } V(\theta)\leq \Vbar \}$, with $\Gt$ being the set of the homogenised Hooke's laws with microstructures of density $\theta$.

\begin{remark}
For general objective functionals $J$, the expression of the homogenised problem~\eqref{eq:topOptH} cannot be made fully explicit, see~\cite{Allaire-book12}.
\end{remark}

It is common, see~\cite{Pantz-PT-08,Pantz-AGP-19}, to limit  the set of admissible composites to periodic composites (e.g., a periodic unit cell square) in order to be able to numerically compute the corresponding homogenised Hooke's law.
For the case under analysis of compliance minimisation with $\nsd=2$, the set $\Gt$ consists of a special class of composites known as sequential laminates. 
For each set of parameters $\bEta \in \Iset$, a material distribution $\theta$ is thus determined in terms of optimal sequential laminates solving~\eqref{eq:topOptH} in order to compute one entry of a dataset of optimised topologies.

\section{Dataset construction via the homogenisation method}
\label{sc:DataHomog}

In this section,  the numerical method used to construct the dataset of optimised topologies is presented.
\hl{
This work considers as benchmark problem the topology optimisation of a two-dimensional cantilever beam composed of a linear elastic material under a parametrised surface load with varying position and angle, leading to a high-dimensional parametric problem in $\RR^4$.
}
For each set of parameters $\bEta \in \Iset$, the displacement field $\buS$ and the stress tensor $\bsigmaS$ are obtained from the solution of the parametric elastic problem~\eqref{eq:elasticityH}.
This information is then employed to determine the material density $\theta$ and the homogenised stress tensor $\matAs$ of the sequential laminates as described in the following subsections.

\subsection{Characterisation of optimal sequential laminates}
\label{sc:Laminates}

Assume $\theta \neq 0$ and let $(\lamSig{i},\vSig{i})$ denote the pair of the $i$-th eigenvalue and eigenvector of $\bsigmaS$. 
Sequential laminates are manufactured by successive layers of void and material along $\nsd$ orthogonal directions given by $\vSig{i}, \, i=1,\ldots,\nsd$.
In 2D, the optimal lamination proportions are defined as
\begin{equation}\label{eq:laminProp}
m_1 := \frac{\abs{\lamSig{2}}}{\abs{\lamSig{1}} + \abs{\lamSig{2}}} , \qquad
m_2 := \frac{\abs{\lamSig{1}}}{\abs{\lamSig{1}} + \abs{\lamSig{2}}} ,
\end{equation}
and represent the volume fractions of the two materials composing the laminate.
The resulting homogenised elasticity tensor is such that
\begin{equation}\label{eq:optA}
[\matAs]^{-1} := [\matA]^{-1} + \frac{1-\theta}{\theta} [m_1 \matAc{1} + m_2 \matAc{2}]^{-1} ,
\end{equation}
where $[\matA]^{-1}$ represents the compliance tensor of the base material corresponding to the full density $\theta=1$ and $\matAc{i}$ denotes the fourth-order elasticity tensor associated with the $i$-th material phase of the composite, weighted by the corresponding volume fraction $m_i$.
\begin{remark}
The coefficient $(1-\theta)/\theta$ is a relative scaling factor accounting for the contribution of the laminates when $\theta \neq 1$. It is worth noticing that, as $\theta$ decreases, the laminated composites contribute more significantly to the homogenised compliance $[\matAs]^{-1}$.  In the limit of $\theta \rightarrow 0$, the homogenised elasticity tensor becomes singular, $\matAs=\bm{0}$, and $[\matAs]^{-1}$ is not well-defined, justifying the assumption of $\theta \neq 0$ introduced at the beginning of the section.
\end{remark}
Finally, for any symmetric matrix $\bxi$ it holds that $\matAc{i}$ is such that
\begin{equation}\label{eq:Ac}
\matAc{i} \bxi : \bxi = \matA \bxi : \bxi - \frac{1}{\mu} \norm{ (\matA \bxi) \vSig{i} }_2^2 + \frac{\mu+\lambda}{\mu(2\mu+\lambda)} ([\vSig{i}]^\top \! (\matA\bxi) \vSig{i})^2 ,
\end{equation}
where $\norm{ \odot }_2 := \sqrt{\odot \cdot \odot}$ denotes the Euclidean norm of the vector $\odot$.
Technical details on the structure of the above mentioned tensors are presented in~\ref{sc:appHooke}.

It is worth noticing that the value of the density $\theta$ depends upon the maximum volume fraction $\Vbar$ for the admissible topologies.
To enforce the maximum volume constraint in $\UadS$, the penalty approach described in~\cite{Pantz-AP-06} is employed. Let $\gamma > 0$ denote a penalty parameter. 
The corresponding optimal value of the density, see~\cite{Allaire-book12}, is thus given by
\begin{equation}\label{eq:optTheta}
\theta := \min \Bigg\{ 1, \left[\frac{2\mu+\lambda}{4\mu(\mu+\lambda) \gamma}\right]^{1/2} (\abs{\lamSig{1}}+\abs{\lamSig{2}}) \Bigg\} .
\end{equation}

\subsection{Topology optimisation algorithm based on homogenisation}
\label{sc:HomogAlgo}

For any $\bEta \in \Iset$,  the topology optimisation problem~\eqref{eq:topOptH} is solved via the finite element-based homogenisation strategy described in~\cite{Pantz-AP-06}.
The method relies on an iterative algorithm that performs minimisation by alternatively computing the stress tensor $\bsigmaS$ and the design variables $(\theta,\matAs)$.
The initial condition is set starting from the elasticity tensor $\matAhs$ fulfilling the Hashin-Shtrikman bounds~\cite{Allaire-book12}.

Consider a discretisation of the computational domain $\D$ using a set of $\numel$ non-overlapping elements $\elemE, \, e=1, \ldots, \numel$. 
The following finite element spaces are introduced:
\begin{subequations}
\begin{align}
\Vspace &:= \{ v \in \mathcal{C}(\overline{\D}) \, : \, v{\mid}_{\elemE} \in \mathbb{P}^2(\elemE), \, e=1,\ldots,\numel \text{ and } v{\mid}_{\partial\elemE \cap \GaEta{D}} = 0 \} ,
\label{eq:spaceFEM} \\
\Wspace &:= \{ w \in \eltwo(\D) \, : \, w{\mid}_{\elemE} \in \mathbb{P}^0(\elemE), \, e=1,\ldots,\numel \} ,
\label{eq:spaceP0}
\end{align}
\end{subequations}
where $\Vspace$ denotes the space of continuous finite element functions, being polynomials of degree at most $2$ in each element $\Omega_e$, and vanishing on $\GaEta{D}$, whereas $\Wspace$ is the space of square-integrable functions in $\D$ with piecewise constant approximation in each mesh element.

At iteration $k+1$ of the optimisation loop (Algorithm~\ref{alg:topOptHiFi}, lines 2-15),  the homogenised elastic problem~\eqref{eq:elasticityH} is solved setting the value of the design parameters to the last computed ones $(\theta_{k},\matAs_{k}) \in \Wspace \times [\Wspace]^{\nsd \times \nsd}$. 
In particular, a continuous Galerkin finite element formulation is employed to determine $\buS_{k+1} \in [\Vspace]^{\nsd}$ such that 
\begin{equation}\label{eq:elastWeak}
\int_{\D}{\matAs_{k} \beps(\buS_{k+1}): \beps(\bv) \, d\bx} = \int_{\GaEta{N}}{\bgEta \cdot \bv \, d\Gamma} \quad \forall \bv \in [\Vspace]^{\nsd} .
\end{equation}

\begin{remark}
From a practical viewpoint,  the discrete form of equation~\eqref{eq:elastWeak} is obtained leveraging Voigt notation to exploit the symmetry of the fourth-order elasticity tensor and the second-order strain tensor, see~\ref{sc:appHooke}.
Note that the employed numerical method is known to experience stability issues for nearly incompressible materials. To avoid unstable scenarios, this work considers the value $\nu=0.3$ for the Poisson's ratio.
Alternative formulations based on mixed hybrid methods and robust in the incompressible limit, see, e.g.,~\cite{RS-SGKH:18,RS-SGH:19} could be considered to circumvent this issue.
\end{remark}

Then, the stress tensor $\bsigmaS_{\!k+1}$ is obtained by post-processing the displacement field $\buS_{k+1}$ computed in~\eqref{eq:elastWeak} as
\begin{equation}\label{eq:stressH}
\bsigmaS_{\!k+1} =\matAs_{k} \beps(\buS_{k+1}) .
\end{equation}
The eigenvalues and eigenvectors $(\lamSig{i},\vSig{i}), \, i=1,2$ of $\bsigmaS_{\!k+1}$ are employed to determine the new density $\theta_{k+1}$ using~\eqref{eq:optTheta} and the new homogenised elasticity tensor $\matAs_{k+1}$ according to~\eqref{eq:optA}.

Finally,  once the relative increments of the compliance and the volume fraction between two consecutive iterations are below user-defined tolerances $\xi_J$ and $\xi_V$, respectively, the resulting design is projected onto a set of \emph{classical} topologies by progressively forcing the density to assume only the values $0$ or $1$ (Algorithm~\ref{alg:topOptHiFi}, lines 16-27).
To achieve this goal,  the previously computed density $\theta_{k+1}$ is replaced by a modified one $\thetaT_{k+1}$ obtained by means of the transformation
\begin{equation}\label{eq:penal}
\thetaT_{k+1} = \frac{1}{2}\bigl(1-\cos(\pi \, \theta_{k+1})\bigr) ,
\end{equation}
until a more restrictive stopping criterion on the increments of the compliance and the volume fraction is fulfilled with tolerances $\xiT_J$ and $\xiT_V$, respectively.
The resulting high-fidelity topology optimisation algorithm based on the homogenisation method is sketched in algorithm~\ref{alg:topOptHiFi}.
\begin{algorithm}[!p]
\caption{High-fidelity topology optimisation algorithm}\label{alg:topOptHiFi}
\resizebox{!}{0.4\textheight}{%
\begin{minipage}{\linewidth}
\begin{algorithmic}[1]
\REQUIRE{Tolerances $\xi_J$, $\xi_V$, $\xiT_J$, and $\xiT_V$ for the stopping criteria on the compliance and the volume fraction, without and with penalisation; initial elasticity tensor $\matAhs$; initial material density distribution $\thetaBar$.}
\STATE{Set $\matAs_0 = \matAhs$ and $\theta_0 = \thetaBar$.} 
\\ \texttt{Optimisation loop allowing intermediate density values}
\FOR{$k=0,1,\ldots$}
\STATE{Compute the displacement $\buS_{k+1}$ solving~\eqref{eq:elastWeak} and the stress tensor $\bsigmaS_{\!k+1}$ as~\eqref{eq:stressH}.}
\STATE{Compute the compliance $\Js(\theta_{k},\matAs_{k})$ as~\eqref{eq:complianceHomog} and the volume fraction $V(\theta_{k})$ as~\eqref{eq:volFrac}.}
\IF{$\abs{ \Js(\theta_{k},\matAs_{k}) {-} \Js(\theta_{k-1},\matAs_{k-1}) } {>} \xi_J  \Js(\theta_{k},\matAs_{k})$
\OR
$\abs{ V(\theta_{k}) {-} V(\theta_{k-1}) } {>} \xi_V V(\theta_{k})$}
\STATE{Compute the updated design density distribution $\theta_{k+1}$ as~\eqref{eq:optTheta}.}
\STATE{Compute the updated homogenised elasticity tensor $\matAs_{k+1}$ as~\eqref{eq:optA} using $\theta_{k+1}$.}
\ELSE
\STATE{Compute the updated design density distribution $\theta_{k+1}$ as~\eqref{eq:optTheta}.}
\STATE{Compute the penalised design density distribution $\thetaT_{0}$ as~\eqref{eq:penal} using $\theta_{k+1}$.}
\STATE{Compute the new initial elasticity tensor $\matAs_{0}$ as~\eqref{eq:optA} using $\thetaT_{0}$.}
\STATE{\textbf{go to} line 16.}
\ENDIF
\STATE{$k \gets k+1$.}
\ENDFOR
\\ \texttt{Optimisation loop penalising intermediate density values}
\FOR{$\ell=0,1,\ldots$}
\STATE{Compute $\buS_{\ell+1}$ solving~\eqref{eq:elastWeak} and $\bsigmaS_{\!\ell+1}$ as~\eqref{eq:stressH}.}
\STATE{Compute the compliance $\Js(\thetaT_{\ell},\matAs_{\ell})$ as~\eqref{eq:complianceHomog} and the volume fraction $V(\thetaT_{\ell})$ as~\eqref{eq:volFrac}.}
\IF{$\abs{ \Js(\thetaT_{\ell},\matAs_{\ell}) {-} \Js(\thetaT_{\ell-1},\matAs_{\ell-1}) } {>} \xiT_J  \Js(\thetaT_{\ell},\matAs_{\ell})$
\OR
$\abs{ V(\thetaT_{\ell}) {-} V(\thetaT_{\ell-1}) } {>} \xiT_V V(\thetaT_{\ell})$}
\STATE{Compute the updated design density distribution $\theta_{\ell+1}$ as~\eqref{eq:optTheta}.}
\STATE{Compute the penalised design density distribution $\thetaT_{\ell+1}$ as~\eqref{eq:penal} using $\theta_{\ell+1}$.}
\STATE{Compute the updated homogenised elasticity tensor $\matAs_{\ell+1}$ as~\eqref{eq:optA} using $\thetaT_{\ell+1}$.}
\ELSE
\STATE{\textbf{return} $\thetaT_{\ell}$.}
\ENDIF
\STATE{$\ell \gets \ell+1$.}
\ENDFOR
\ENSURE{Optimal material density distribution $\thetaT_{\ell}$.}
\end{algorithmic}
\end{minipage}
}
\end{algorithm}

\subsection{Dataset of optimised topologies}
\label{sc:Dataset}

Algorithm~\ref{alg:topOptHiFi} provides the optimal design field for one configuration of the elastic problem~\eqref{eq:elasticityH}, that is, for one set of parameters $\bEta \in \Iset$.  Starting from this, a dataset of optimised topologies is constructed computing the solution of the optimisation problem~\eqref{eq:topOptH}, for different values of the parameters.
%
\hl{
In this work, a demonstrator is presented for the parametric design of a minimum compliance two-dimensional cantilever beam under maximum volume constraint ($\nsd = 2$), with each dataset entry being computed for a fixed pair of parameters $\bEta = (\eta_1,\eta_2)^\top$, respectively controlling the position where the external force $\bgEta$ is applied and its orientation ($\npa=2$).
Extensions to three-dimensional design with higher number of parameters, alternative objective functionals and geometric constraints, and nonlinear models will be the focus of future studies.
}

Let the computational domain be $\D := [-1,1] \times [0,1]$ and define $\Ga{b} := \{ (x,y) \in \partial\D : y=0 \}$, $\Ga{r} := \{ (x,y) \in \partial\D : x=1 \}$, and $\Ga{t} := \{ (x,y) \in \partial\D : y=1 \}$.
The structure, with material parameters $E=1$ and $\nu=0.3$, is clamped on $\GaEta{D} = \Ga{D} := \{ (x,y) \in \partial\D : x=-1 \}$, whereas the load is applied on $\GaEta{N}$, a subset of $\Ga{b} \cup \Ga{r} \cup \Ga{t}$ of length $\ell_N = 10^{-1}$. 
%
%
%
The boundary $\GaEta{N}$ is defined using the parameter $\eta_1$ as curvilinear abscissa of the contour, namely,
\begin{subequations}
\begin{equation}\label{eq:GammaN}
\GaEta{N} :=
\left\{
\begin{array}{lll}
x \in \left[\sInit{b} + \eta_1 - \dfrac{\ell_N}{2}, \sInit{b}+\eta_1+\dfrac{\ell_N}{2}\right], &y=0,  \\
& \text{for $(x,y) \in \Ga{b}$ and $\eta_1 \in \I_{1,b}$}, \\[1.2em]
x =1 ,& y \in \left[\sInit{r}+\eta_1 - \dfrac{\ell_N}{2}, \sInit{r}+\eta_1+\dfrac{\ell_N}{2}\right],  \\
& \text{for $(x,y) \in \Ga{r}$ and $\eta_1 \in \I_{1,r}$}, \\[1.2em]
x \in \left[\sInit{t}-\eta_1 - \dfrac{\ell_N}{2}, \sInit{t}-\eta_1+\dfrac{\ell_N}{2}\right], &y=1,  \\
& \text{for $(x,y) \in \Ga{t}$ and $\eta_1 \in \I_{1,t}$}, 
\end{array}
\right.
\end{equation}
with 
\begin{align}\label{eq:paramEta1}
\I_{1,b} &:= [0,6 \times 10^{-1}] , && \sInit{b}=3.5 \times 10^{-1} , \\
\I_{1,r} &:= [7 \times 10^{-1},16 \times 10^{-1}] , && \sInit{r}=-6.5 \times 10^{-1}, \\
\I_{1,t} &:= [17 \times 10^{-1},23 \times 10^{-1}] , && \sInit{t}=26.5 \times 10^{-1} ,
\end{align}
\end{subequations}
%
%
where each interval $\I_{1,\diamond}$ is discretised using 
\begin{equation}\label{eq:nPos}
\texttt{n}_{1,\diamond} = 1 + \frac{2}{\ell_N} \left( \max_{\eta \in \I_{1,\diamond}}{\eta} - \min_{\eta \in \I_{1,\diamond}}{\eta} \right)
\end{equation}
equally-spaced points (including the extrema).

The parameter $\eta_2$ controls the angle at which the force (with module $1$) is applied, namely
\begin{equation}\label{eq:ForceAngle}
\bgEta := \left[\cos\left(\dfrac{\pi}{2}\left(1 - \dfrac{\eta_2}{30}\right)\right), \sin\left(\dfrac{\pi}{2}\left(1 - \dfrac{\eta_2}{30}\right)\right)\right]^\top, \quad \eta_2 \in \I_2,
\end{equation}
where $\I_2 :=[0,59]$ is discretised with $\texttt{n}_{2}=60$ equally-spaced points (including the extrema), yielding a resolution for the angle of the force of $3^{\circ}$ within the interval $[-87^{\circ},90^{\circ}]$.
Summarising, the dataset contains $45$ positions where the load is applied with $60$ different angles, for a total of $\nS = (\texttt{n}_{1,b}+\texttt{n}_{1,r}+\texttt{n}_{1,t})\texttt{n}_{2} = 2,700$ optimised topologies obtained using the homogenisation method in algorithm~\ref{alg:topOptHiFi}, with $\xi_J=\xi_V=0.5 \times 10^{-2}$,  $\xiT_J=\xiT_V=0.5 \times 10^{-4}$, $\thetaBar=0.4$, and $\matAhs$ the elasticity tensor fulfilling the Hashin-Shtrikman bounds (see Appendix~\ref{sc:appHooke}).

The computation of each instance of the dataset requires the numerical solution of multiple partial differential equations. This is performed using the continuous Galerkin finite element method on a uniform mesh of $\numel=25,600$ triangular elements obtained by subdividing each squared element in a structured grid of $160 \times 80$ quadrilaterals with a diagonal.
The optimal material density distribution provided as output by algorithm~\ref{alg:topOptHiFi} is discretised as a piecewise constant finite element function in each triangular element of the mesh, that is, $\thetaT \in \Wspace$.

A preprocessing step projects the information of the $\numel$ mesh elements onto the underlying grid of quadrilaterals by averaging the values of $\thetaT$ corresponding to the two triangles within each squared element. This yields a vector $\btheta$ of dimension $\nT = 12,800$, with piecewise constant material density, which will be employed as \emph{ground truth} solution in the training process described in the following section.

\hl{
The $2,700$ optimised topologies are subdivided into training, validation, and test sets with a proportion of $80\%$, $10\%$, and $10\%$ of cases, respectively.
In order to employ the same simulation data for both the interpolation and the extrapolation studies presented in Section~\ref{sc:Simulations}, a Cartesian sampling is constructed in the parametric spaces $\I_{1,b} \times \I_2$, $\I_{1,r} \times \I_2$, and $\I_{1,t} \times \I_2$.
It is well-known that this sampling strategy is not optimal and its cost becomes computationally prohibitive in high-dimensional spaces~\cite{Liu-LPC-18}.
Reducing the number of samples required during training is outside the scope of the present work and can be achieved using standard approaches in the ROM literature, including statistical sampling to maximise the explored parametric space~\cite{Conover-79-Latin,Gunzburger-DFG-99,Gunzburger-SGB-07} and adaptive sampling to balance local exploitation and global exploration of the parametric domain~\cite{Sobester-SLK-05,Willcox-BWG-08,Breitkopf-PBBVZ-20}.
The rest of this work focuses on demonstrating the feasibility of \emph{learning} fundamental features of optimised topologies from an existing dataset of solutions, in the context of parametric topology optimisation problems.
}

\section{Surrogate models of optimised topologies}
\label{sc:Surrogate}

Whilst one instance of problem~\eqref{eq:topOptH}, that is, the topology optimisation of an elastic structure for a fixed set of parameters $\bEta \in \Iset$,  can be solved effectively using well-established topology optimisation strategies such as algorithm~\ref{alg:topOptHiFi}, the multiple queries introduced by the dependence of the problem on $\bEta$ yield the need for efficient parametric surrogate models of the topology.

Starting from the collection of optimised topologies presented in Section~\ref{sc:Dataset},  this section introduces a surrogate model for the parametric topology optimisation problem~\eqref{eq:topOptH}.
It is worth noticing that the proposed methodology goes beyond the classical notion of ROM-accelerated optimisation, commonly obtained by replacing the full-order solver in each iteration of the optimiser by means of a surrogate model of the physical problem~\cite{Breitkopf-FBK-08,Breitkopf-RHXBV-13}.
Indeed, this work proposes: (i) to construct a surrogate model of the optimised solution, that is,  a \emph{quasi-optimal} topology; (ii) to devise a surrogate-based optimisation algorithm employing it as an \emph{educated} initial condition in order to significantly reduce the number of iterations required by algorithm~\ref{alg:topOptHiFi} to attain convergence.

To achieve this goal,  different neural network architectures are proposed, leveraging dimensionality reduction capabilities of autoencoders to extract relevant features from $\btheta$ and functional approximation properties of feed-forward nets to create a map between the space of parameters $\Iset$ and the reduced version of the space of admissible topologies $\UadS$.

\subsection{A parametric encoder/feed-forward/decoder net}
\label{sc:Combined}

\begin{figure}[!htb]
	\centering
	\includegraphics[width=0.7\textwidth]{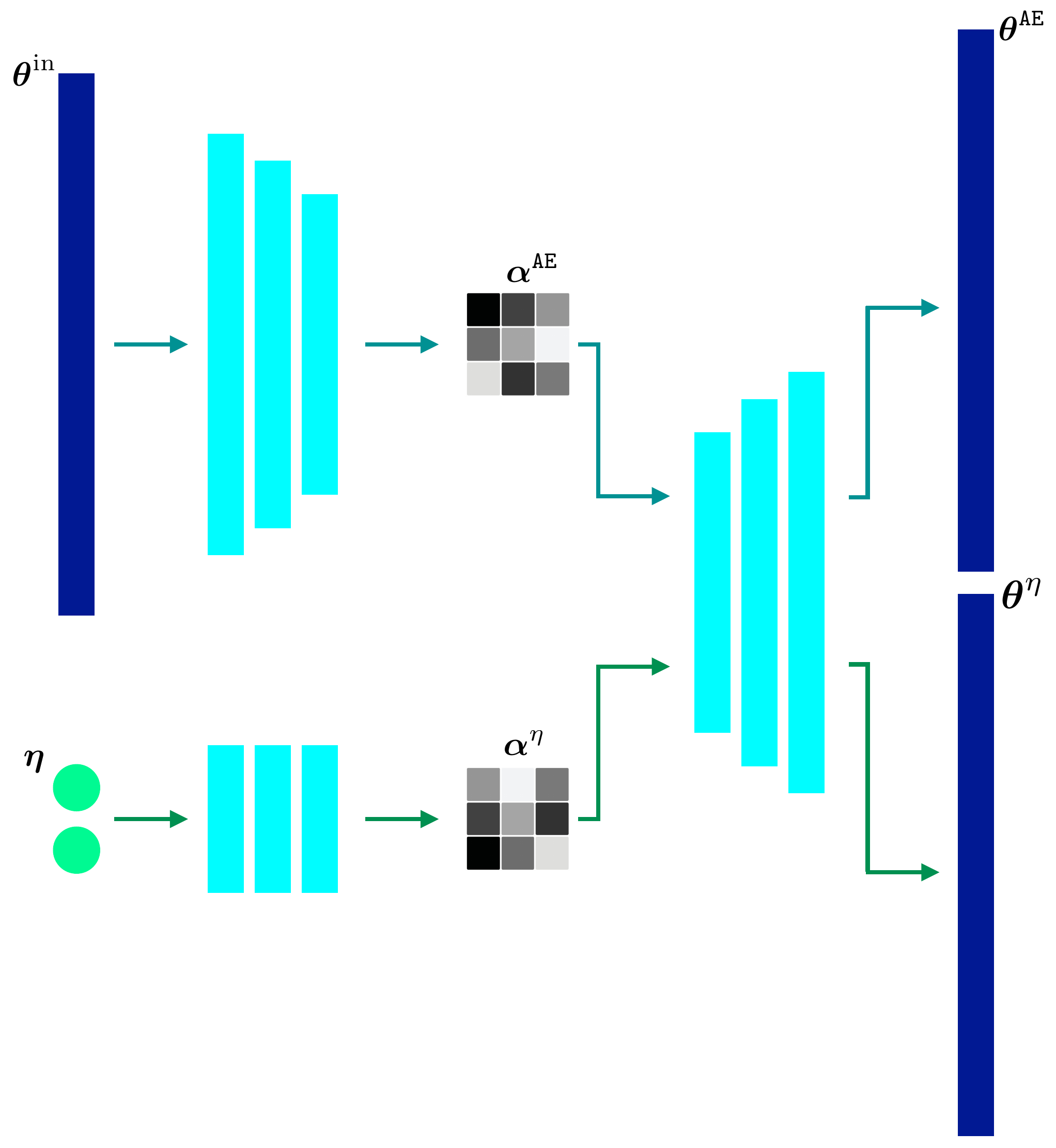}
	
	\caption{$\texttt{E-FF}_{\!\eta}\texttt{-D}$: Architecture of the network combining an encoder block $\texttt{E}$, a feed-forward net $\texttt{FF}_{\!\eta}$, and a decoder block $\texttt{D}$.}
	\label{fig:Combined}
\end{figure}

The first proposed architecture combines, in a monolithic approach, an encoder-decoder for dimensionality reduction and a feed-forward net to learn the map between the input parameters and the latent space (Figure~\ref{fig:Combined}). The net, denoted as $\texttt{E-FF}_{\!\eta}\texttt{-D}$, trains:
\begin{itemize}
\item an encoder $\texttt{E}$ to reduce the dimension $\nT$ of the input data $\bthetaIn$ to the dimension $\nA$ of the latent space $\balphaAE$;
\item a feed-forward net $\texttt{FF}_{\!\eta}$ learning the functional mapping from the input parameters $\bEta$ of dimension $\npa$ to another low-dimensional space $\balphaEta$ of dimension $\nA$;
\item a decoder-type block $\texttt{D}$ to reconstruct the high-dimensional representation of the input data of dimension $\nT$ from the low-dimensional space of dimension $\nA$, concurently upscaling the latent space descriptions $\balphaAE$ to $\bthetaAE$ and $\balphaEta$ to $\bthetaEta$.
\end{itemize}

The training of the net $\texttt{E-FF}_{\!\eta}\texttt{-D}$ is performed via the minimisation of the loss function
\begin{equation}\label{eq:lossCombined}
\Loss{\texttt{E-FF}_{\!\eta}\texttt{-D}} = \LAE + \LEta + \NNpar{\alpha} \frac{\nT}{\nA} \LCode + \NNpar{r} \frac{\nT}{\nNNpar} \Reg ,
\end{equation}
obtained as the linear combination of: (i) the standard loss function $\LAE$ of the autoencoder; (ii) the loss function $\LEta$ assessing the mismatch between the reconstruction $\bthetaEta$ and the \emph{ground truth} $\bthetaIn$; (iii) the loss function $\LCode$ measuring the discrepancy between the latent spaces $\balphaAE$ and $\balphaEta$, respectively constructed by the autoencoder and the feed-forward net; (iv) an $\ell^1$ or least absolute shrinkage and selection operator (LASSO) regularisation term $\Reg$,  providing a convex relaxation of the combinatorial $\ell^0$ regularisation to promote sparsity among the trainable parameters $\NNpar{i}, \ i=1, \ldots,\nNNpar$ of the network.
The expressions of the terms appearing in~\eqref{eq:lossCombined} are detailed below:
\begin{subequations}\label{eq:lossesDef}
\begin{equation}\label{eq:LAE}
\LAE := \frac{1}{\nS \nT} \sum_{i=1}^{\nS} \norm{ \bthetaAE_i - \bthetaIn_i }_2^2 ,
\end{equation}
\begin{equation}\label{eq:LEta}
\LEta := \frac{1}{\nS \nT} \sum_{i=1}^{\nS} \norm{ \bthetaEta_i - \bthetaIn_i }_2^2 ,
\end{equation}
\begin{equation}\label{eq:LCode}
\LCode := \frac{1}{\nS \nA} \sum_{i=1}^{\nS} \norm{ \balphaEta_i - \balphaAE_i }_2^2 ,
\end{equation}
\begin{equation}\label{eq:Reg}
\Reg := \sum_{j=1}^{\nNNpar} \abs{ \NNpar{j} } ,
\end{equation}
\end{subequations}
where the subindex $_i$ denotes the $i$-th entry of the dataset comprising $\nS$ optimised topologies described by the vectors $\bthetaIn_i \in \RR^{\nT}, \ i=1,\ldots,\nS$,  $\norm{ \cdot }_2$ is the Euclidean norm of the corresponding vectors,  $\nT/\nA$ denotes the compression factor executed during dimensionality reduction,  $\nT/\nNNpar$ represents the input-to-parameter ratio, and $\NNpar{\alpha}$ and $\NNpar{r}$ are two non-trainable user-defined coefficients.

The above network is designed with a \emph{high compression} rate to construct a latent representation much smaller than the input data, a low input-to-parameter ratio to retain \emph{high capacity} in order to accurately represent complex functional mappings, and a \emph{sparse} network exploiting the $\ell^1$ regularisation to reduce the number of active parameters while avoiding overfitting. The details of the architecture are provided in~\ref{sc:appNets}.

\subsubsection{Staggered approaches to encoder/feed-forward/decoder nets}
\label{sc:Staggered}

The training of $\texttt{E-FF}_{\!\eta}\texttt{-D}$ is performed monolithically by minimising the loss function~\eqref{eq:lossCombined} using a suitable strategy, for instance the Adam optimiser for the numerical results presented in Section~\ref{sc:Simulations}.
A variant of this approach can be devised following the idea of a staggered training algorithm proposed in~\cite{Duraisamy-XD-20,Corigliano-GTBC-24}.
\begin{figure}[!htb]
	\centering
	\includegraphics[width=0.7\textwidth]{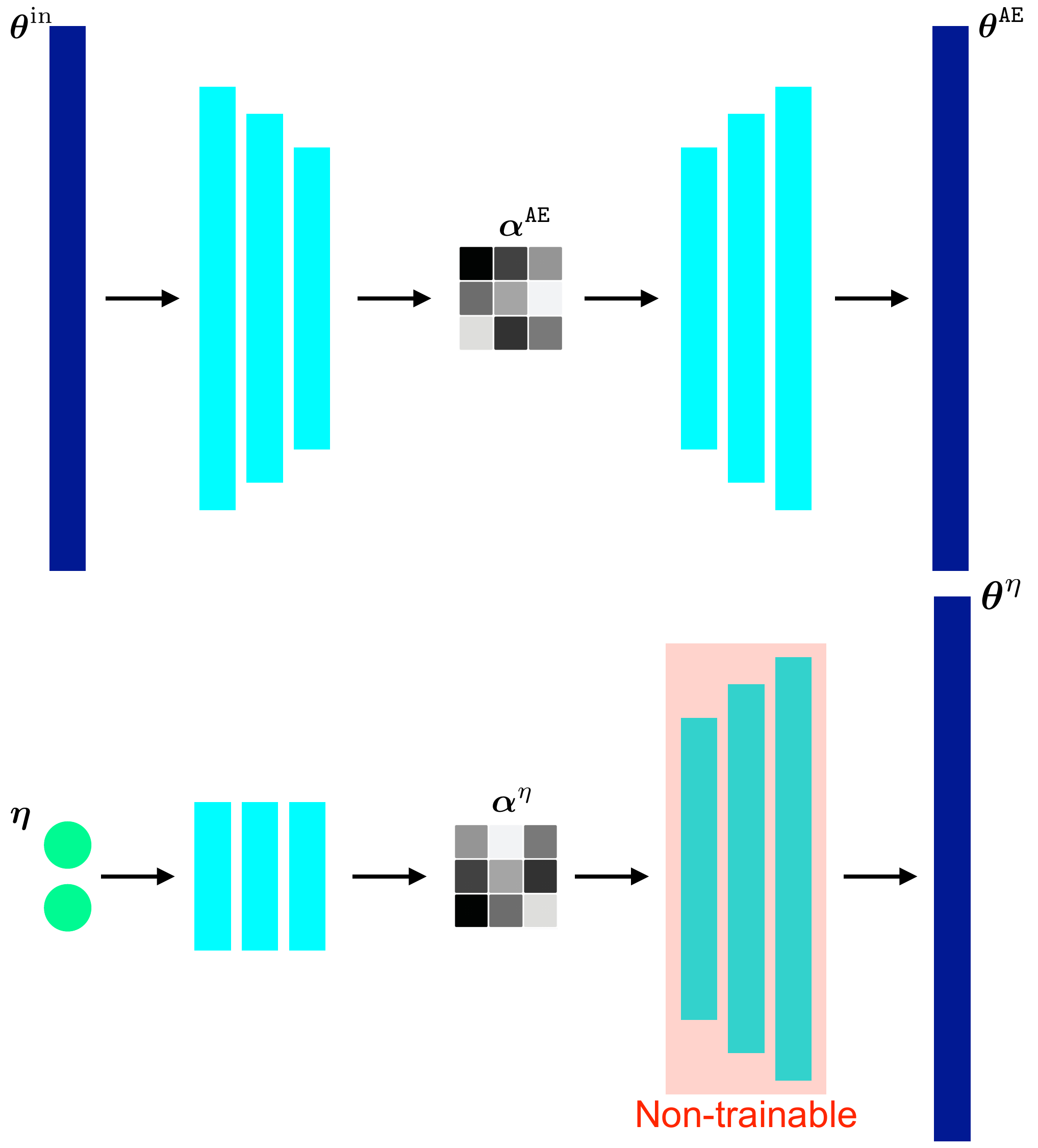}
	
	\caption{$\texttt{S-AE-FF}_{\!\eta}\texttt{-D}$: Architecture of the network comprising the decoupled autoencoder $\texttt{AE}$ and the feed-forward net $\texttt{FF}_{\!\eta}$, trained in a staggered ($\texttt{S}$) approach,  with a non-trainable decoder block $\texttt{D}$.}
	\label{fig:StaggeredWithDecoder}
\end{figure}
In particular, the network $\texttt{S-AE-FF}_{\!\eta}\texttt{-D}$ in Figure~\ref{fig:StaggeredWithDecoder} features:
\begin{itemize}
\item a classic autoencoder $\texttt{AE}$ reducing the dimension of input data $\bthetaIn$ from $\nT$ to the latent dimension $\nA$ of $\balphaAE$ and subsequently reconstructing the high-dimensional representation $\bthetaAE$ of dimension $\nT$;
\item a feed-forward net $\texttt{FF}_{\!\eta}$ learning the functional mapping between the input parameters $\bEta$ of dimension $\npa$ and another low-dimensional space $\balphaEta$ of dimension $\nA$;
\item a non-trainable decoder $\texttt{D}$ obtained from the autoencoder $\texttt{AE}$ to reconstruct the high-dimensional representation $\bthetaEta$ from the latent representation $\balphaEta$ of dimension $\nA$.
\end{itemize}
Note that the $\texttt{S}$ in the name $\texttt{S-AE-FF}_{\!\eta}\texttt{-D}$ denotes the staggered nature of this architecture,  since the net is trained following a two-step strategy. 
The first stage of the training of the $\texttt{S-AE-FF}_{\!\eta}\texttt{-D}$ net minimises the discrepancy between $\bthetaAE$ and $\bthetaIn$ with LASSO regularisation, namely,
\begin{subequations}\label{eq:lossStagWithD}
\begin{equation}\label{eq:lossStagWithD1}
\Loss{\texttt{S-AE-FF}_{\!\eta}\texttt{-D}}^{\mathrm{I}} = \LAE + \NNpar{r,\mathrm{I}} \frac{\nT}{\nNNpar[,\mathrm{I}]} \Reg_{\mathrm{I}} ,
\end{equation}
whereas the second step determines the learned mapping from the parameters $\bEta$ to the high-dimensional reconstruction $\bthetaEta$ through the minimisation of the loss function
\begin{equation}\label{eq:lossStagWithD2}
\Loss{\texttt{S-AE-FF}_{\!\eta}\texttt{-D}}^{\mathrm{II}} = \LEta + \NNpar{r,\mathrm{II}} \frac{\nA}{\nNNpar[,\mathrm{II}]} \Reg_{\mathrm{II}} .
\end{equation}
\end{subequations}
In the above equations,  $\Reg_{\mathrm{I}}$ (respectively, $\Reg_{\mathrm{II}}$) denotes the LASSO regularisation~\eqref{eq:Reg} particularised for the $\nNNpar[,\mathrm{I}]$ (respectively, $\nNNpar[,\mathrm{II}]$) trainable parameters of the autoencoder $\texttt{AE}$ (respectively,  feed-forward net $\texttt{FF}_{\!\eta}$). Moreover,  $\NNpar{r,\mathrm{I}}$ and $\NNpar{r,\mathrm{II}}$ are two user-defined regularisation coefficients.

An alternative staggered approach with same number of trainable parameters can be devised modifying the second step of the training procedure in equation~\eqref{eq:lossStagWithD} by neglecting the decoder block and exclusively training the map between the parameters $\bEta \in \Iset$ and the latent representation $\balphaEta$ of dimension $\nA$.
\begin{figure}[!htb]
	\centering
	\includegraphics[width=0.7\textwidth]{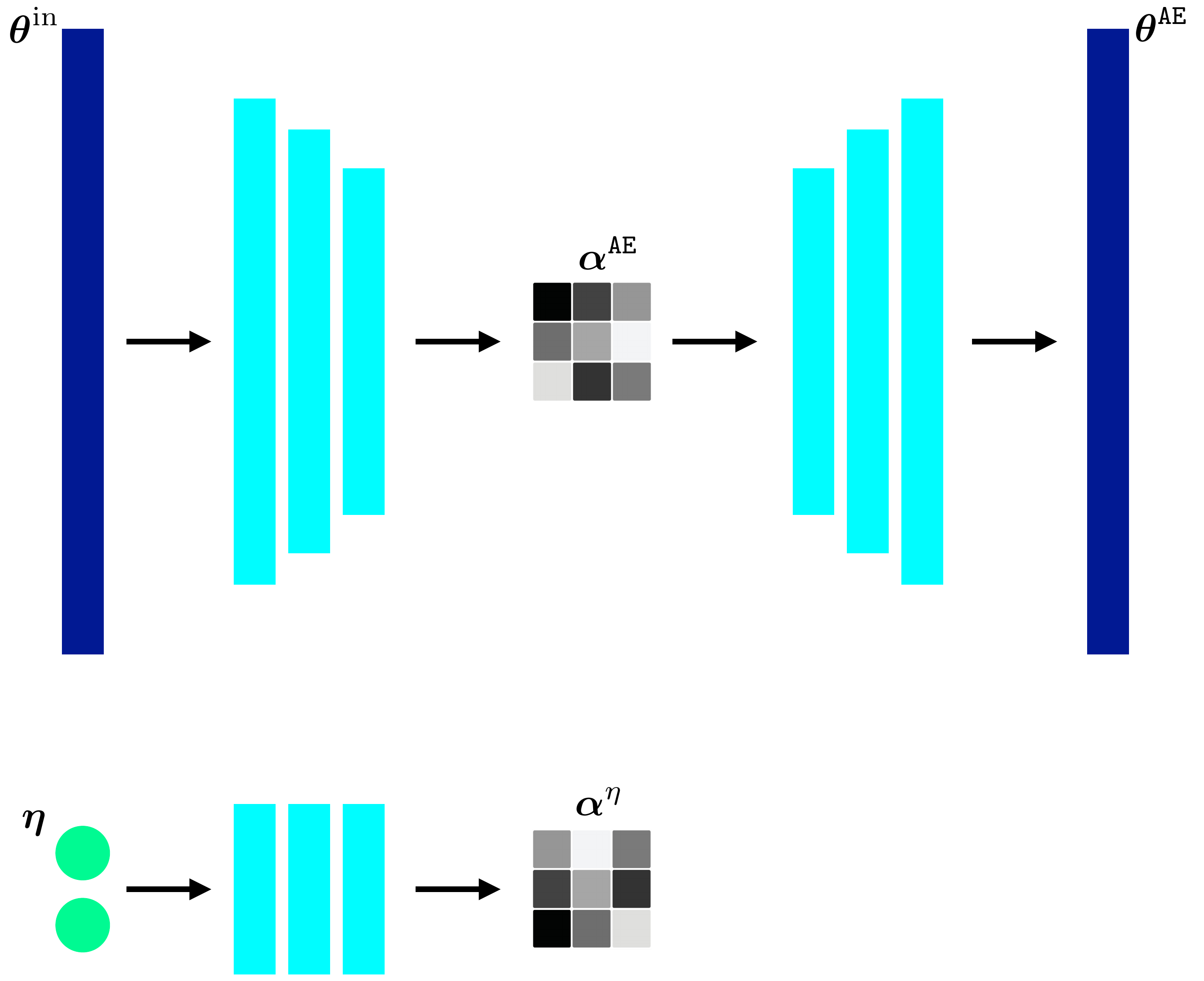}
	
	\caption{$\texttt{S-AE-FF}_{\!\eta}$: Architecture of the network comprising the decoupled autoencoder $\texttt{AE}$ and feed-forward net $\texttt{FF}_{\!\eta}$, trained in a staggered ($\texttt{S}$) approach.}
	\label{fig:StaggeredOnlyCode}
\end{figure}
The architecture of the resulting network $\texttt{S-AE-FF}_{\!\eta}$ is depicted in Figure~\ref{fig:StaggeredOnlyCode} and consists of:
\begin{itemize}
\item a classic autoencoder $\texttt{AE}$ reducing the dimension of input data $\bthetaIn$ from $\nT$ to the latent dimension $\nA$ of $\balphaAE$ and subsequently reconstructing the high-dimensional representation $\bthetaAE$ of dimension $\nT$;
\item a feed-forward net $\texttt{FF}_{\!\eta}$ learning the functional mapping between the input parameters $\bEta$ of dimension $\npa$ and the previously determined latent space of dimension $\nA$.
\end{itemize}
Also in this case,  $\texttt{S}$ denotes a staggered training strategy. First, a classical autoencoder is trained by minimising
\begin{subequations}\label{eq:lossStagComb}
\begin{equation}\label{eq:lossStagComb1}
\Loss{\texttt{S-AE-FF}_{\!\eta}}^{\mathrm{I}} = \LAE + \NNpar{r,\mathrm{I}} \frac{\nT}{\nNNpar[,\mathrm{I}]} \Reg_{\mathrm{I}} .
\end{equation}
Then,  for each entry of the dataset $\bthetaIn_i, \ i=1,\ldots,\nS$, the corresponding latent representation $\balphaAE_i$ is constructed using the previously trained encoder block. This is employed as \emph{ground truth} data for the second training minimising the loss function
\begin{equation}\label{eq:lossStagComb2}
\Loss{\texttt{S-AE-FF}_{\!\eta}}^{\mathrm{II}} = \LCode + \NNpar{r,\mathrm{II}} \frac{\nA}{\nNNpar[,\mathrm{II}]} \Reg_{\mathrm{II}} .
\end{equation}
\end{subequations}
It is worth noticing that, whilst the size of the output vector evaluated in the loss function~\eqref{eq:lossStagComb2} is significantly smaller than the one in equation~\eqref{eq:lossStagWithD2} (namely $\nA$ instead of $\nT$), the number of trainable parameters is the same in both cases, thus leading to a comparable algorithm complexity as verified also by the similar training times reported in Section~\ref{sc:Validation}.

\hl{
\subsubsection{Promoting sparsity in the surrogate models}
\label{sc:Sparsity}

The choice of $\ell^2$ norms in the loss functions~\eqref{eq:lossesDef} stems from the need of the network to construct a description of the topology that is accurate in the entire domain, not just in a few key points. This promotes completeness of the learned representation, penalising global large errors while fostering the identification of relevant features from raw data without a priori information.
Moreover, the LASSO regularisation applied to the network parameters guarantees the construction of parsimonious models with few coefficients that can be efficiently evaluated during the online parameter exploration procedures.

The use of $\ell^1$ norms for the loss functions has also been explored during the preliminary phase of the design of the network architecture. The results, not reported here for brevity, showed that this approach is suitable to preserve sharper details in the topology reconstruction. Nonetheless, the optimisation procedure was observed to be slower due to the non-smoothness of the associated gradient and sparsity was not achieved in the final representation of the \emph{quasi-optimal} topologies.

To preserve the optimisation advantages provided by the $\ell^2$ norm of the loss functions while promoting sparsity, alternative strategies include replacing $\Reg$ in~\eqref{eq:Reg} by means of the $\ell^1$ regularisation of the activation functions of the hidden layers~\cite{Jiang-JZZX-13} or with the Kullback-Leibler divergence penalty on the average activations~\cite{Qi-QLWL-17}.
The resulting sparse autoencoder architectures have shown promising results in the construction of reduced order models, see, e.g.,~\cite{Hernandez-HBGCC-21}.
These techniques, outside the scope of the present work, disentangle the effect of individual neurons encouraging feature learning, with each neuron in the net specialising in the detection of specific characteristics of the input data.
This represents a feature of interest to be explored in future works in order to enhance the generalisation capabilities of the proposed architectures and achieve better predictions of \emph{quasi-optimal} topologies for configurations unseen during training.
}

\subsection{A parametric feed-forward/decoder net}
\label{sc:Parametric}

In this section, a parsimonious architecture for a surrogate model of \emph{quasi-optimal} topologies is devised to enhance the performance of the solutions introduced in Section~\ref{sc:Combined}.

The surrogate models presented therein are characterised by three common blocks ($\texttt{E}$, $\texttt{FF}_{\!\eta}$, and $\texttt{D}$), suitably organised according to different architectures. 
In particular,  the feed-forward net $\texttt{FF}_{\!\eta}$ upscales the input parameters $\bEta$ from a space of dimension $\npa$ to the latent space of dimension $\nA$,  generally with $\nA \gg \npa$.
Moreover,  the space of dimension $\nA$ for the latent representation is constructed by training the encoder block $\texttt{E}$ to compress the input data of dimension $\nT$, with $\nT \gg \nA$.

It is worth noticing that the global number of trainable parameters is the same in the three surrogate models of Section~\ref{sc:Combined}, with most of the unknowns being associated with the encoder and decoder blocks. 
Indeed, fully-connected convolutional layers are employed to successfully compress (respectively, decompress) the information in the encoding (respectively, decoding) part.
This leads to a rapid growth of the number of parameters in the nets, starting from the input layer of dimension $\nT=12,800$. Detailed information on the specifics of the nets is provided in~\ref{sc:appNets}.

\begin{figure}[!htb]
	\centering
	\includegraphics[width=0.7\textwidth]{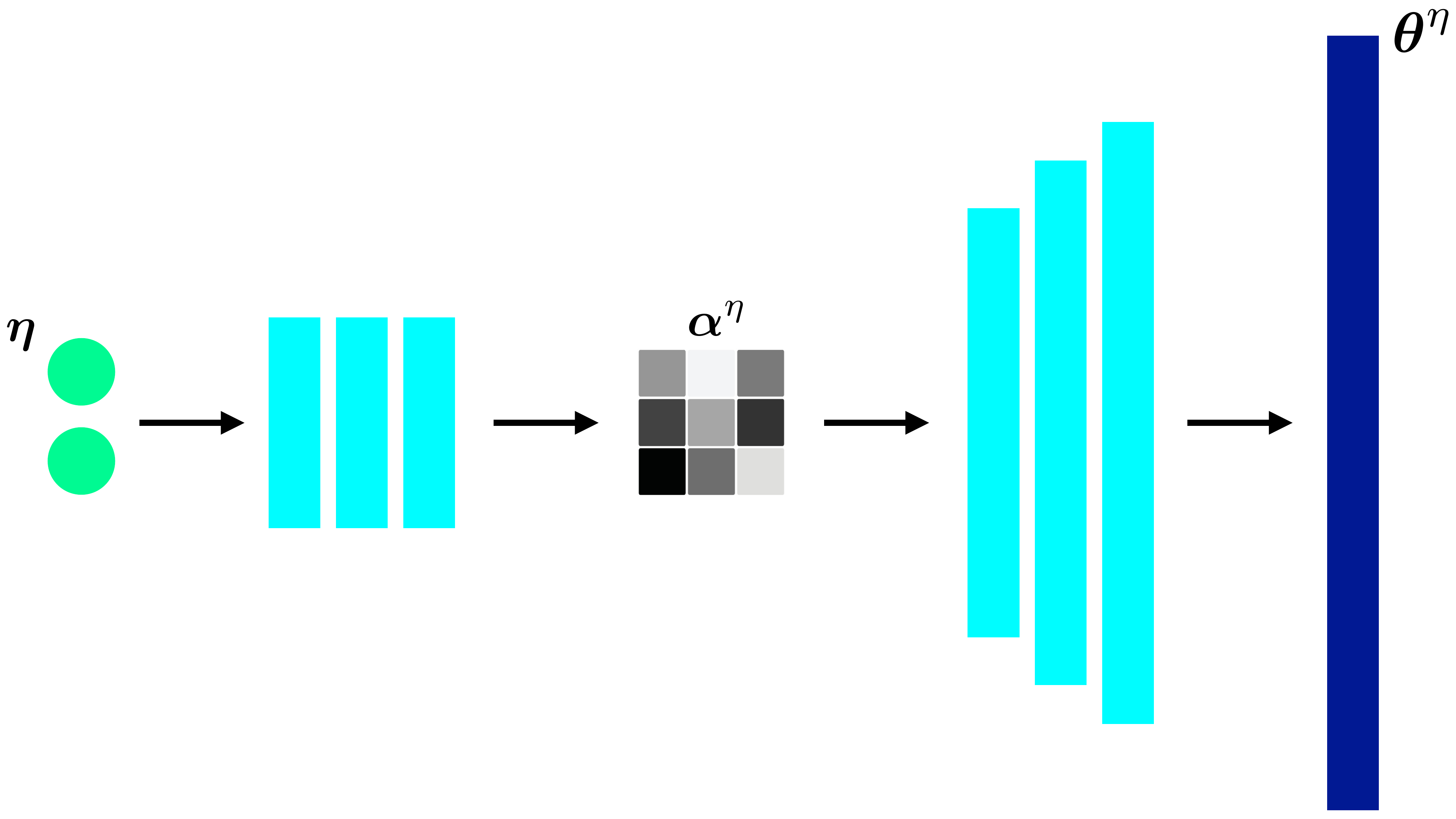}
	
	\caption{$\texttt{FF}_{\!\eta}\texttt{-D}$: Architecture of the network combining a feed-forward net $\texttt{FF}_{\!\eta}$ and a decoder block $\texttt{D}$.}
	\label{fig:Parametric}
\end{figure}

In order to circumvent this issue, a parsimonius architecture is devised by assuming that the latent representation $\balphaEta$ provides a parametric description of the input data $\bthetaIn$ in a latent space of dimension $\nA$, but this is not obtained as a classic autoencoder compression of the \emph{ground truth} data.
In this context, the encoder block $\texttt{E}$ can be neglected and the architecture of the resulting $\texttt{FF}_{\!\eta}\texttt{-D}$ net, displayed in Figure~\ref{fig:Parametric},  contains:
\begin{itemize}
\item a feed-forward net $\texttt{FF}_{\!\eta}$ learning the functional mapping from the input parameters $\bEta$ of dimension $\npa$ to a low-dimensional representation $\balphaEta$ of dimension $\nA$;
\item a block $\texttt{D}$ mimicking a decoder with the goal of reconstructing the high-dimensional representation $\bthetaEta$ of dimension $\nT$ from the low-dimensional space of dimension $\nA$.
\end{itemize}

Similarly to the case of $\texttt{E-FF}_{\!\eta}\texttt{-D}$, the training of $\texttt{FF}_{\!\eta}\texttt{-D}$ is performed monolithically, by minimising the loss function
\begin{equation}\label{eq:lossParametric}
\Loss{\texttt{FF}_{\!\eta}\texttt{-D}} = \LEta + \NNpar{r} \frac{\nT}{\nNNpar} \Reg ,
\end{equation}
measuring the discrepancy between the input data $\bthetaIn$ and the reconstruction $\bthetaEta$, under a suitably scaled $\ell^1$ regularisation.
The details of the architecture are provided in~\ref{sc:appNets}.

\subsection{Surrogate-based topology optimisation algorithm}
\label{sc:OptiAlg}

The surrogate models described in Section~\ref{sc:Combined} and~\ref{sc:Parametric} provide a prediction $\bthetaEta$ of a \emph{quasi-optimal} topology,  for a given value $\bEta \in \Iset$ of the parameters (Algorithm~\ref{alg:topOptSurrogate}, line 1).
It is worth noticing that, once the models have been trained, only the architectures $\texttt{FF}_{\!\eta}\texttt{-D}$, $\texttt{S-AE-FF}_{\!\eta}\texttt{-D}$, and $\texttt{S-AE-FF}_{\!\eta}$ can perform such a prediction using exclusively the information of the parameters $\bEta$ as input data.
On the contrary,  the $\texttt{E-FF}_{\!\eta}\texttt{-D}$ net requires both the parameters $\bEta$ and the \emph{ground truth} topology $\bthetaIn$ as input data to reconstruct the corresponding value of $\bthetaEta$, thus making this architecture unsuitable for the surrogate-based optimisation algorithm~\ref{alg:topOptSurrogate}.
Similarly,  the autoencoder $\texttt{AE}$ compresses the input topology $\bthetaIn$ into a low-dimensional latent space before reconstructing $\bthetaEta$ and cannot predict the \emph{quasi-optimal} topology starting from the parameters $\bEta$, which is the first step of this surrogate-based optimisation strategy (Algorithm~\ref{alg:topOptSurrogate}, line 1).

Denote by $\thetaEta$ the piecewise constant functional approximation of the material density in the finite element space $\Wspace$.  Note that the vector of elemental values associated with $\thetaEta$ has dimension $\numel$, whereas $\bthetaEta \in \RR^{\nT}$. As explained in Section~\ref{sc:Dataset},  the finite element mesh is generated subdividing each of the $\nT$ squared element in the computational grid into two triangular elements by means of its diagonal. Consequently,  the value of $\thetaEta$ in each mesh element $\Omega_e, \ e =1,\ldots, \numel$ is inherited by its corresponding parent cell in the structured grid.

The density distribution $\thetaEta$ predicted using the surrogate models is thus employed as an \emph{educated} initialisation for a surrogate-based topology optimisation algorithm designed starting from the high-fidelity strategy presented in algorithm~\ref{alg:topOptHiFi}.
In particular,  the method first constructs a mechanically-admissible initialisation by determining an homogenised elasticity tensor compatible with the initial material distribution (Algorithm~\ref{alg:topOptSurrogate}, lines 3-4). Then it performs an optimisation loop penalising the intermediate values of the density in order to attain a distribution featuring only the values $0$ or $1$ (Algorithm~\ref{alg:topOptSurrogate}, lines 5-16).
\begin{algorithm}[!htb]
\caption{Surrogate-based topology optimisation algorithm}\label{alg:topOptSurrogate}
\begin{algorithmic}[1]
\REQUIRE{Tolerances $\xiT_J$ and $\xiT_V$ for the stopping criteria on the compliance and the volume with penalisation; initial elasticity tensor $\matAhs$; values of the problem parameters $\bEta$.}
\STATE{Obtain an educated initialisation $\thetaEta$ evaluating the surrogate model for the input parameter $\bEta$.} 
\STATE{Set $\thetaT_{0}=\thetaEta$.}
\\ \texttt{Preprocess of NN output for a mechanically-admissible initialisation}
\STATE{Compute the displacement $\buS$ solving~\eqref{eq:elastWeak} and the stress tensor $\bsigmaS$ as~\eqref{eq:stressH} using the elasticity tensor $\matAhs$.}
\STATE{Compute the homogenised elasticity tensor $\matAs_0$ as~\eqref{eq:optA} using $\thetaT_{0}$.}
\\ \texttt{Optimisation loop penalising intermediate density values}
\FOR{$\ell=0,1,\ldots$}
\STATE{Compute $\buS_{\ell+1}$ solving~\eqref{eq:elastWeak} and $\bsigmaS_{\!\ell+1}$ as~\eqref{eq:stressH}.}
\STATE{Compute the compliance $\Js(\thetaT_{\ell},\matAs_{\ell})$ as~\eqref{eq:complianceHomog} and the volume fraction $V(\thetaT_{\ell})$ as~\eqref{eq:volFrac}.}
\IF{$\abs{ \Js(\thetaT_{\ell},\matAs_{\ell}) {-} \Js(\thetaT_{\ell-1},\matAs_{\ell-1}) } {>} \xiT_J  \Js(\thetaT_{\ell},\matAs_{\ell})$
\OR
$\abs{ V(\thetaT_{\ell}) {-} V(\thetaT_{\ell-1}) } {>} \xiT_V V(\thetaT_{\ell})$}
\STATE{Compute the updated design density distribution $\theta_{\ell+1}$ as~\eqref{eq:optTheta}.}
\STATE{Compute the penalised design density distribution $\thetaT_{\ell+1}$ as~\eqref{eq:penal} using $\theta_{\ell+1}$.}
\STATE{Compute the updated homogenised elasticity tensor $\matAs_{\ell+1}$ as~\eqref{eq:optA} using $\thetaT_{\ell+1}$.}
\ELSE
\STATE{\textbf{return} $\thetaT_{\ell}$.}
\ENDIF
\STATE{$\ell \gets \ell+1$.}
\ENDFOR
\ENSURE{Optimal material density distribution $\thetaT_{\ell}$.}
\end{algorithmic}
\end{algorithm}

\section{Numerical experiments}
\label{sc:Simulations}

In this section, the approximation properties of the surrogate models described in Section~\ref{sc:Surrogate} are numerically assessed using the dataset introduced in Section~\ref{sc:Dataset}. 
To this end, Section~\ref{sc:Validation} and~\ref{sc:Generalisation} respectively evaluate the interpolation and extrapolation (i.e., generalisation) properties of the surrogate models, that is,  their ability to efficiently extract relevant features from high-dimensional topology representations and accurately reproduce unseen topologies during the online phase. All presented architectures can perform this task, with the autoencoder $\texttt{AE}$ providing the target optimal performance in terms of compression and reconstruction.

Moreover, the surrogate-based optimisation algorithm~\ref{alg:topOptSurrogate} is evaluated when the \emph{quasi-optimal} topologies are used as \emph{educated} initial conditions, comparing accuracy and efficiency with the optimal solutions provided by the high-fidelity strategy in algorithm~\ref{alg:topOptHiFi}. 
This task, discussed in Section~\ref{sc:InterpolationOpti} and~\ref{sc:ExtrapolationOpti},  requires that the online evaluation of the surrogate model can be performed using only the parameters $\bEta$ as input data. Hence,  only the surrogate models $\texttt{FF}_{\!\eta}\texttt{-D}$, $\texttt{S-AE-FF}_{\!\eta}\texttt{-D}$, and $\texttt{S-AE-FF}_{\!\eta}$ can be employed in the surrogate-based optimisation algorithm~\ref{alg:topOptSurrogate}.

\subsection{Validation of the surrogate models}
\label{sc:Validation}

\begin{figure}[!htb]
	\centering

	\subfigure[Dataset $\Aset_1$]{\includegraphics[width=0.3\textwidth]{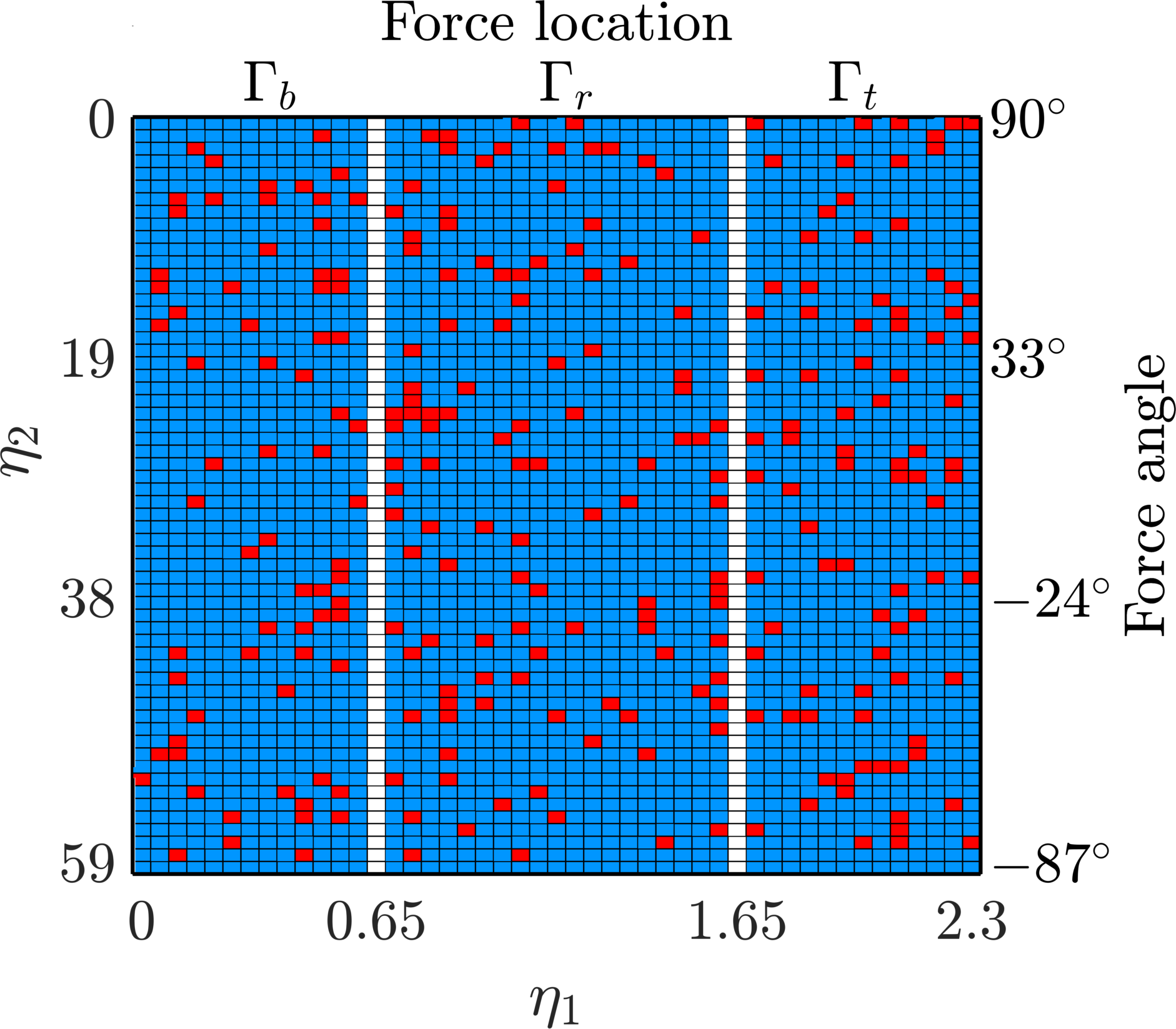}}
	\hspace{3pt}
	\subfigure[Dataset $\Aset_2$]{\includegraphics[width=0.3\textwidth]{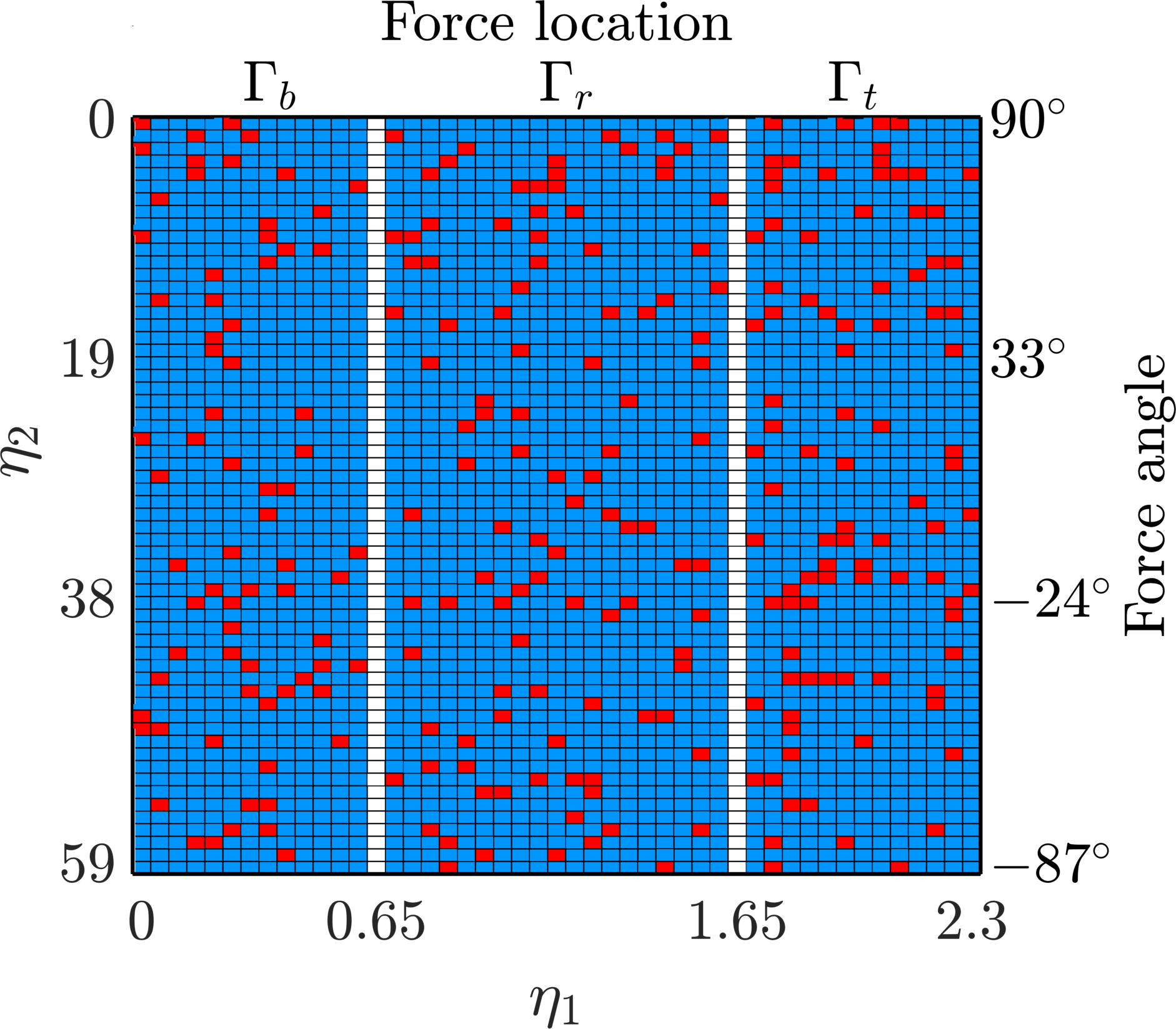}}
	\hspace{3pt}
	\subfigure[Dataset $\Aset_3$]{\includegraphics[width=0.3\textwidth]{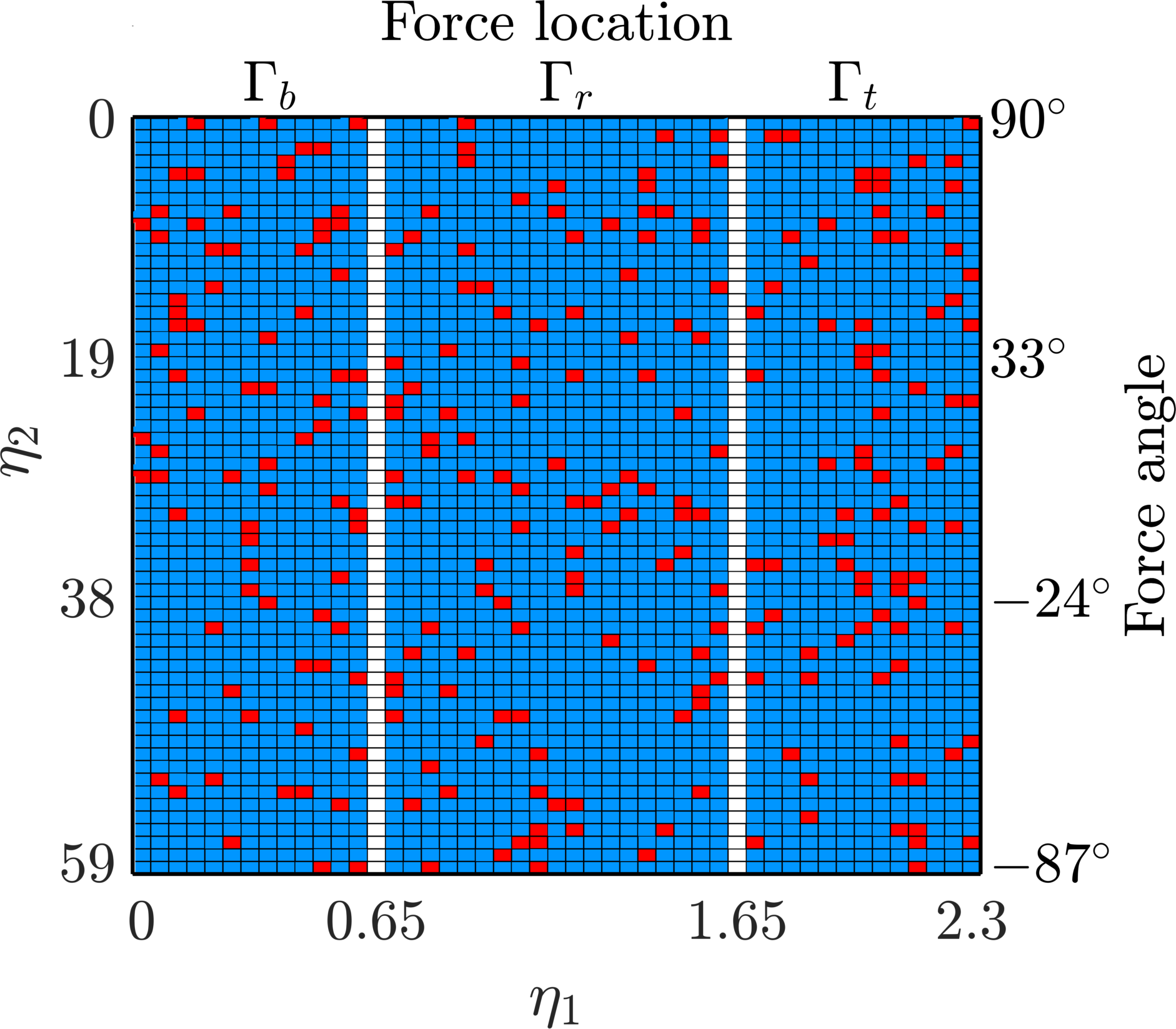}}
	
	\subfigure[Dataset $\Aset_4$]{\includegraphics[width=0.3\textwidth]{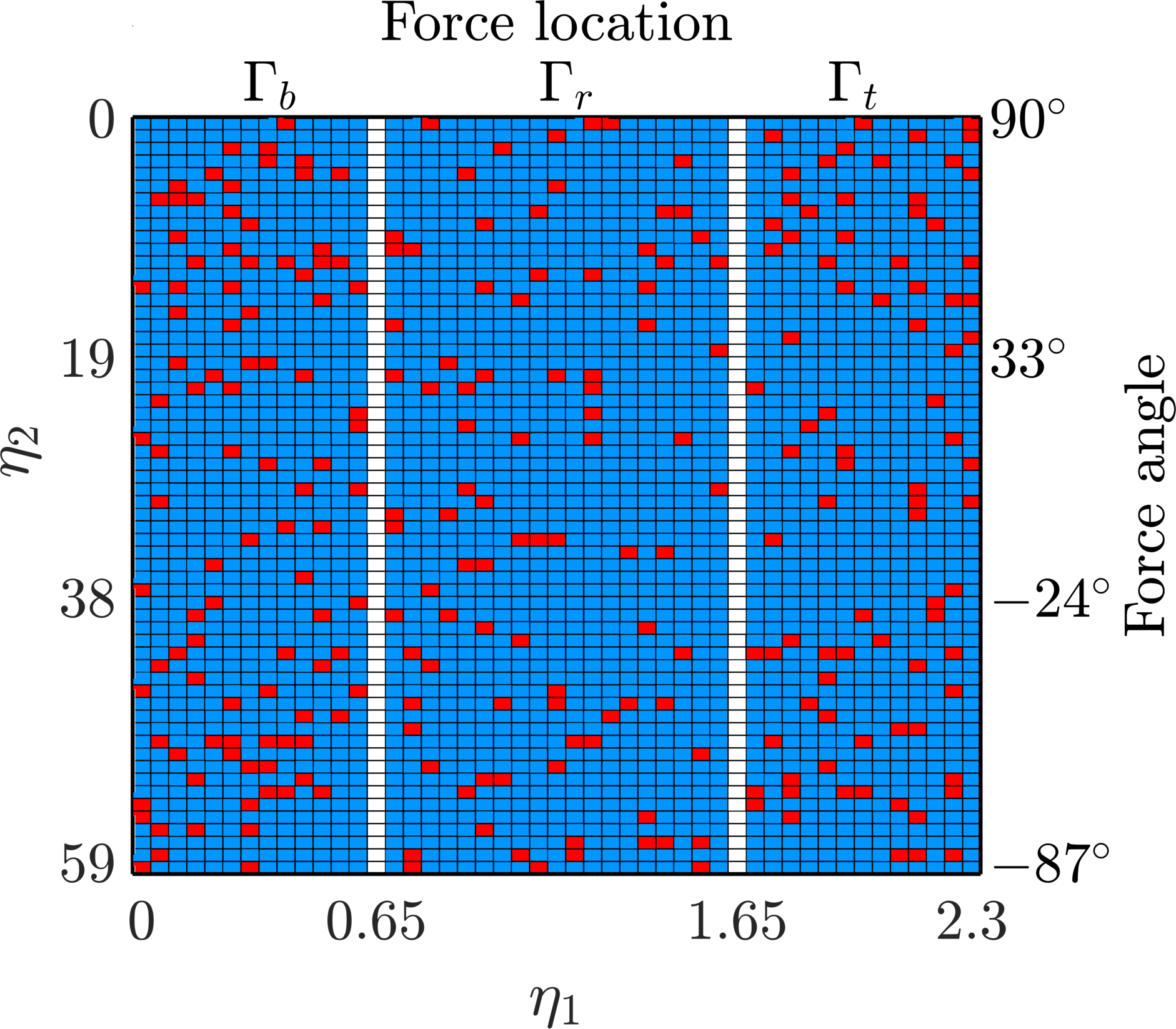}}
	\hspace{3pt}
	\subfigure[Dataset $\Aset_5$]{\includegraphics[width=0.3\textwidth]{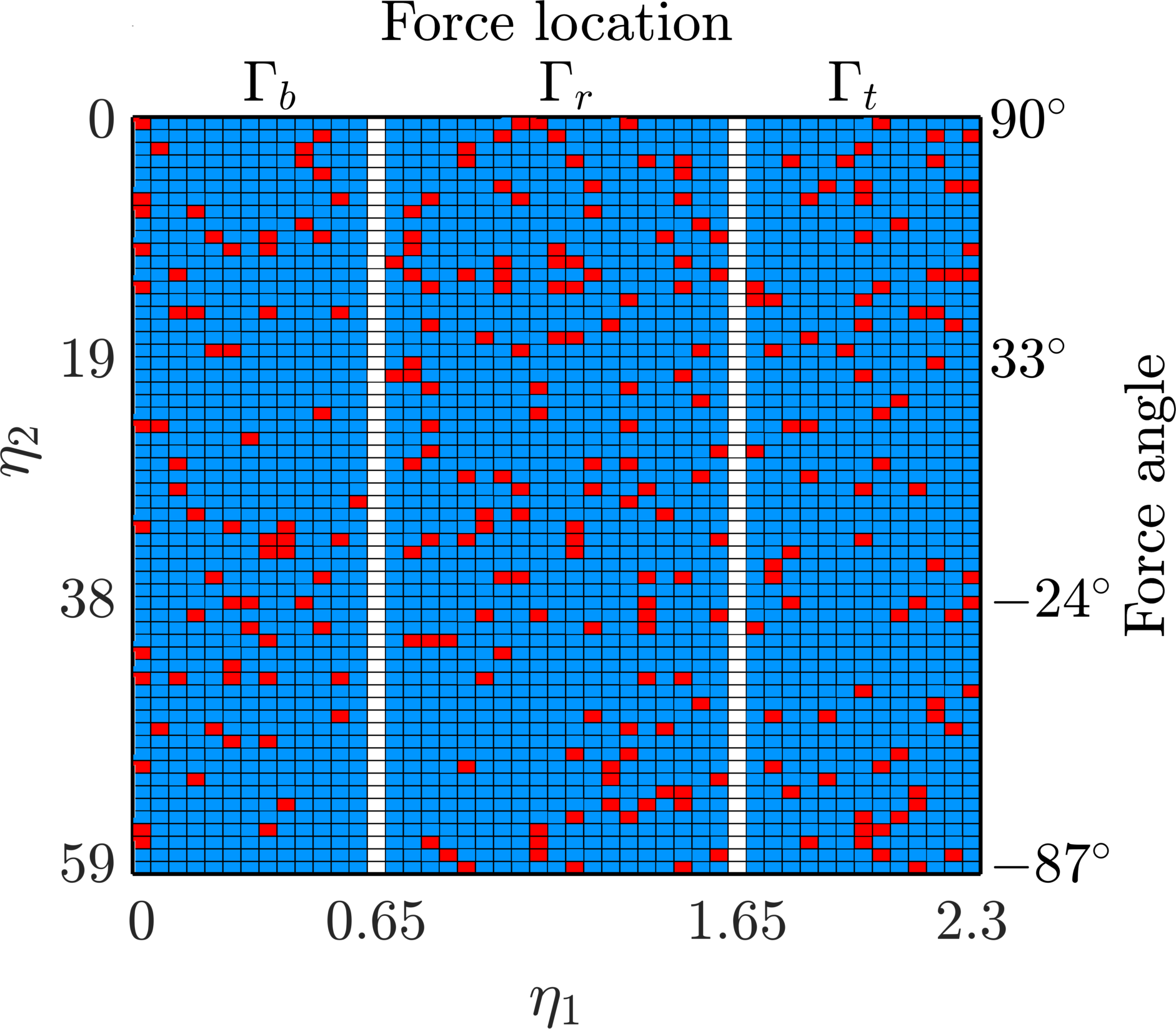}}
	
	\caption{Datasets for cross-validation. Five random splittings of the dataset of $2,700$ optimal topologies into training and validation cases (blue) and test cases (red).}
	\label{fig:CrossValData}
\end{figure}

To evaluate the accuracy of the \emph{quasi-optimal} topologies predicted by the surrogate models, a cross-validation study is performed. Starting from the dataset of $2,700$ optimal topologies, five datasets are generated by randomly splitting the entries into $\nTrain=2,160$ training cases,  $\nVal=270$ validation and $\nTest = 270$ test cases (Figure~\ref{fig:CrossValData}).

The accuracy of the surrogate models is measured as a function of the relative error between the prediction $\bthetaEta$ and the \emph{ground truth} $\bthetaIn$, averaged on the testing dataset that contains the $\nTest$ cases not seen by the NNs during training and validation.  In particular, the relative mean squared error (rMSE, $\Err{2}$) and the relative mean absolute error (rMAE, $\Err{1}$) are defined as
\begin{subequations}
\begin{align}
\Err{2} &:= \frac{1}{\nTest} \sum_{i=1}^{\nTest} \frac{\norm{ \bthetaEta_i - \bthetaIn_i }_2^2}{\norm{ \bthetaIn_i }_2^2} ,
\label{eq:errTheta2} \\
\Err{1} &:= \frac{1}{\nTest} \sum_{i=1}^{\nTest} \frac{\norm{ \bthetaEta_i - \bthetaIn_i }_1}{\norm{ \bthetaIn_i }_1} ,
\label{eq:errTheta1} 
\end{align}
\end{subequations}
where $\norm{ \cdot }_2$ and $\norm{ \cdot }_1$ respectively denote the Euclidean and the Manhattan norms.

\begin{figure}[!htb]
	\centering
	\subfigure[Relative mean squared error]{\includegraphics[height=0.45\textwidth]{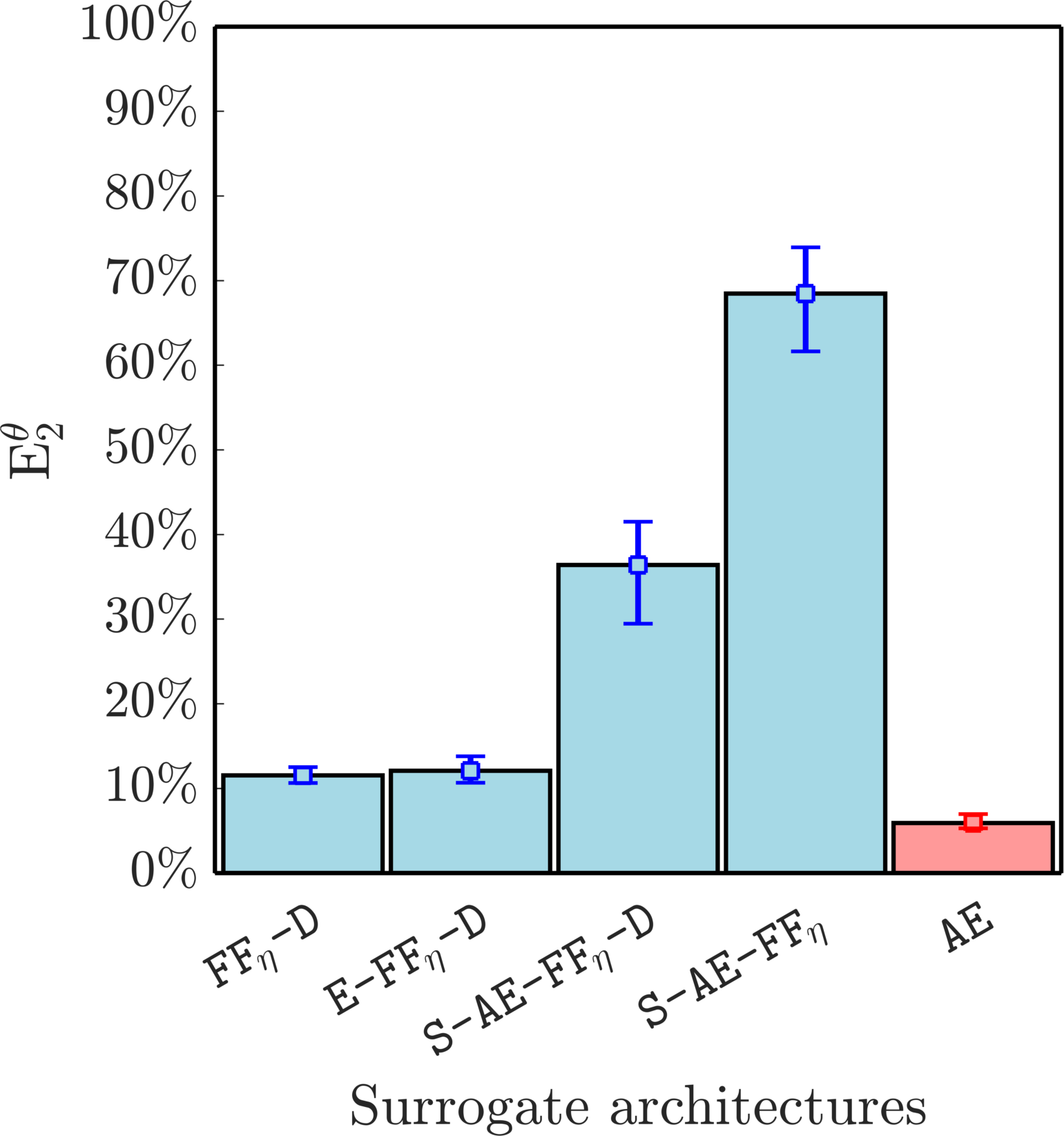}}
	\hspace{5pt}
	\subfigure[Relative mean absolute error]{\includegraphics[height=0.45\textwidth]{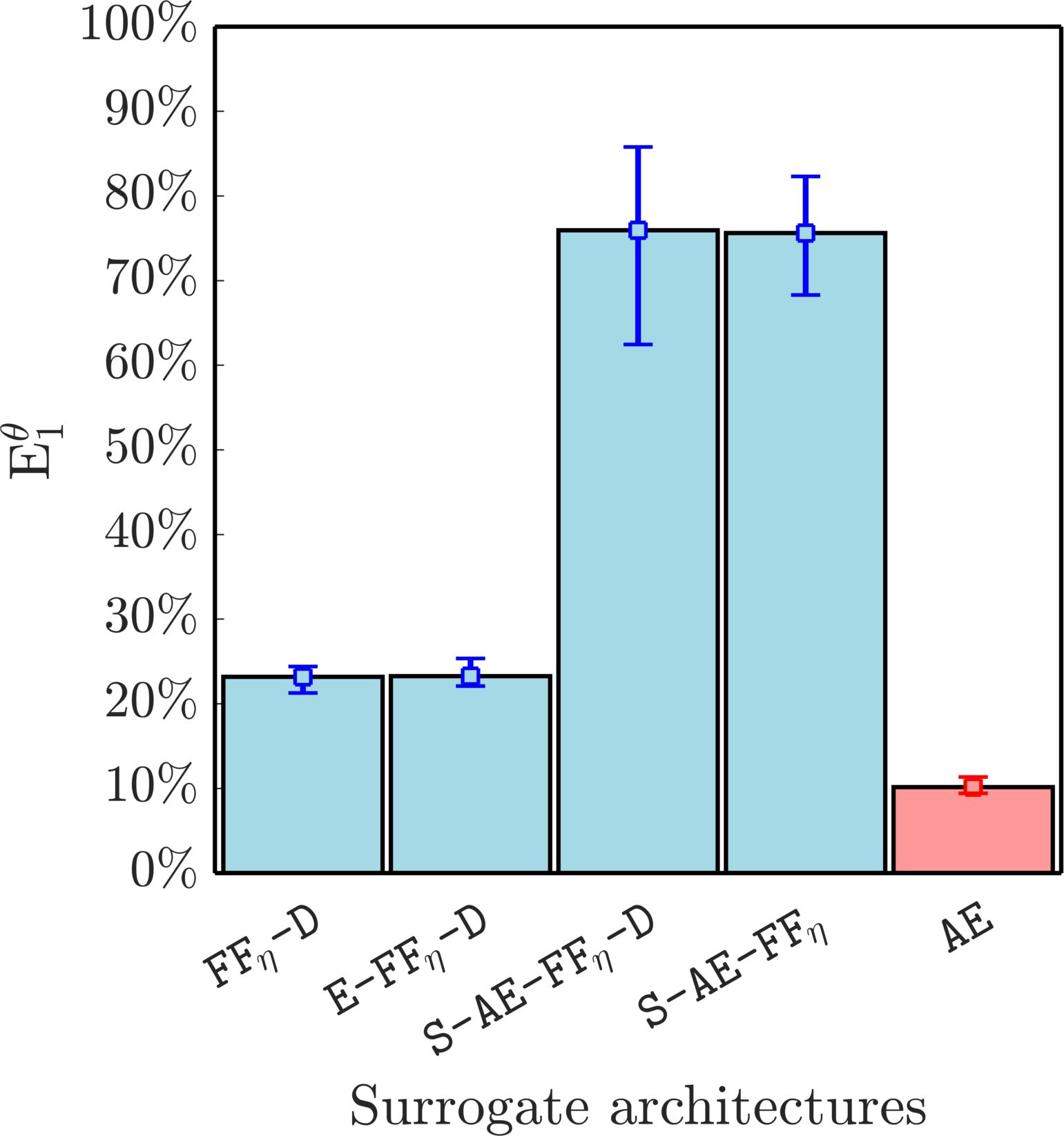}}
	
	\caption{Cross-validation errors ((a): rMSE; (b): rMAE) computed for different surrogate models as the average over the test cases in the datasets $\Aset_k, \ k=1,\ldots,5$.}
	\label{fig:CrossValErr}
\end{figure}

Figure~\ref{fig:CrossValErr} reports the rMSE and the rMAE averaged over the datasets $\Aset_k, \ k=1,\ldots,5$,  for the different proposed architectures of the surrogate models. The results are compared to the target optimal performance of the autoencoder $\texttt{AE}$ (only compressing and decompressing the \emph{ground truth} input data $\bthetaIn$ during the online phase) with the same layers and neurons of the encoder and decoder blocks defined in the surrogate models (see~\ref{sc:appNets} for the details).
The $\texttt{FF}_{\!\eta}\texttt{-D}$ and $\texttt{E-FF}_{\!\eta}\texttt{-D}$ surrogate models present similar accuracies, with average rMSE of $11.55\%$ and $12.10\%$ and average rMAE of $23.18\%$ and $23.26\%$, whereas the baseline $\texttt{AE}$ achieves an average error of $5.91\%$ for rMSE and $10.17\%$ for rMAE.
It is worth noticing that, during the online phase,  the $\texttt{E-FF}_{\!\eta}\texttt{-D}$ net receives as input data both the parameters $\bEta$ and the \emph{ground truth} topology $\bthetaIn$, whereas $\texttt{FF}_{\!\eta}\texttt{-D}$ attains comparable accuracy starting exclusively from the input parameters $\bEta$.
Both architectures clearly outperform the staggered approaches $\texttt{S-AE-FF}_{\!\eta}\texttt{-D}$ and $\texttt{S-AE-FF}_{\!\eta}$ which attain average values of rMSE of $36.41\%$ and $68.45\%$, respectively,  and average rMAE of approximately $76\%$.
Moreover, the variation between the maximum and minimum value of rMSE and rMAE using different validation datasets $\Aset_k, \ k=1,\ldots,5$ is comprised between $1.89\%$ and $3.27\%$ for $\texttt{FF}_{\!\eta}\texttt{-D}$ and $\texttt{E-FF}_{\!\eta}\texttt{-D}$, whereas it grows up to $23.33\%$ for the staggered strategies.

It is worth noticing that whilst $\texttt{FF}_{\!\eta}\texttt{-D}$ and $\texttt{E-FF}_{\!\eta}\texttt{-D}$ provide similar accuracy, their complexity is significantly different. Indeed,  in all the proposed architectures, the majority of trainable parameters is concentrated in the encoder block $\texttt{E}$ and in the decoder block $\texttt{D}$.  Since $\texttt{FF}_{\!\eta}\texttt{-D}$ does not include the encoder block in the architecture,  it features a total of approximately $2.6$ millions of unknowns with respect to the $5.18$ millions of $\texttt{E-FF}_{\!\eta}\texttt{-D}$ (Table~\ref{tab:complexity}).
This directly affects the training time $\TTrain$ of the surrogate models. Denote by $\TTrainAv$ the average training time (measured in minutes) of a model. Table~\ref{tab:complexity} shows that the training of $\texttt{FF}_{\!\eta}\texttt{-D}$ is 2.5 times faster than the one of $\texttt{E-FF}_{\!\eta}\texttt{-D}$ and the speed-up is even more significant when considering the staggered approaches $\texttt{S-AE-FF}_{\!\eta}\texttt{-D}$ and $\texttt{S-AE-FF}_{\!\eta}$.
\begin{table}[!htb]
	\centering
	\resizebox{\linewidth}{!}{%
	\begin{tabular}{| l || c | c | c | c || c |}
	\hline
 	& $\texttt{FF}_{\!\eta}\texttt{-D}$ & $\texttt{E-FF}_{\!\eta}\texttt{-D}$ & $\texttt{S-AE-FF}_{\!\eta}\texttt{-D}$ & $\texttt{S-AE-FF}_{\!\eta}$ & $\texttt{AE}$ \\
	\hline
	$\nNNpar$ & $2,597,025$ & $5,179,850$ & $5,178,425$ ($\nNNpar[,\mathrm{I}]$) & $5,178,425$ ($\nNNpar[,\mathrm{I}]$) & $5,178,425$ \\
	& & & $1,425$ ($\nNNpar[,\mathrm{II}]$) & $1,425$ ($\nNNpar[,\mathrm{II}]$) & \\
	\hline
	$\TTrainAv$ & $60$ & $151$ & $95 (\mathrm{I}) + 70 (\mathrm{II})$ & $95 (\mathrm{I}) + 71 (\mathrm{II})$ & $95$ \\
	\hline
	$\nA/\nA^{\text{max}}$ & $18/25$ & $11.8/25$ & $19.8/25$ & $14.4/25$ & $15.8/25$ \\
	\hline
	\end{tabular}
	}
	
	\caption{Computational performance of the different surrogate models.}
	\label{tab:complexity}
\end{table}

Finally,  note that $\texttt{FF}_{\!\eta}\texttt{-D}$ features an average number of active entries in the latent space $\nA=18$ over the user-defined maximum value $\nA^{\text{max}}=25$.  This entails an effective compression factor $\nT/\nA$ of approximately $711$, slightly below the optimal compression factor of $810$ achieved by the autoencoder $\texttt{AE}$, but higher than the target value of $512$ provided during training.

Given its superior performance with respect to the two staggered approaches $\texttt{S-AE-FF}_{\!\eta}\texttt{-D}$ and $\texttt{S-AE-FF}_{\!\eta}$, henceforth only the $\texttt{FF}_{\!\eta}\texttt{-D}$ model will be considered for the surrogate-based optimisation algorithm.
Moreover, recall that, besides being significantly less efficient than $\texttt{FF}_{\!\eta}\texttt{-D}$ from the computational viewpoint, the $\texttt{E-FF}_{\!\eta}\texttt{-D}$ net is not suitable to be integrated in the optimisation loop because it requires both the parameters $\bEta$ and the \emph{ground truth} $\bthetaIn$ as input data to predict the \emph{quasi-optimal} topology to be employed as \emph{educated} initial guess in algorithm~\ref{alg:topOptSurrogate}.

\hl{
\begin{remark}
The training time reported in Table~\ref{tab:complexity} does not account for the construction of the dataset.
Whilst in the present work the dataset has been synthetically generated by executing the optimisation algorithm~\ref{alg:topOptHiFi} and running the corresponding finite element analyses for all pairs $\bEta \in \Iset$, in realistic scenarios, the presented methodology is meant to exploit existing data of previously optimised topologies.
Of course, if no such previous knowledge is available, the cost of populating a high-quality dataset significantly exceeds the computation of a single instance of the topology optimisation problem~\eqref{eq:elasticityH} and the high-fidelity algorithm~\ref{alg:topOptHiFi} outperforms any surrogate-based approach.
Nonetheless, the results in the following sections will show that, starting from a dataset of optimised designs, algorithm~\ref{alg:topOptSurrogate} is capable of extracting relevant features of the optimised topologies and learn the functional relation between the parameters $\bEta \in \Iset$ and the \emph{quasi-optimal} topology $\bthetaEta$. 
\end{remark}
}

\subsection{Surrogate-based optimisation interpolating \emph{quasi-optimal} topologies}
\label{sc:InterpolationOpti}

In this section, the surrogate model $\texttt{FF}_{\!\eta}\texttt{-D}$ is employed to construct \emph{quasi-optimal} topologies to be employed as \emph{educated} initial guesses for the surrogate-based optimisation algorithm~\ref{alg:topOptSurrogate}.
For each of the five datasets in Figure~\ref{fig:CrossValData}, a random test case is selected according to the pair of parameters $\bEta=(\eta_1,\eta_2)^\top$ reported in Table~\ref{tab:dataOnlineInterp}, with the corresponding position and angle of the external force.
\begin{table}[!htb]
	\centering
	\resizebox{\linewidth}{!}{%
	\begin{tabular}{| c | c || c | c || c | c | c | c |}
	\hline
	\multirow{2}{*}{Case} & \multirow{2}{*}{Dataset}& \multicolumn{2}{c||}{Parameters} & \multicolumn{4}{c|}{Force} \\
	\cline{3-8}
	& & $\eta_1$ & $\eta_2$ & Boundary & Position ($x$) & Position ($y$) & Angle \\
	\hline
	 1 & $\Aset_1$ & $0.45$ & $55$ & $\Ga{b}$ & $[0.75 ,0.85]$ & $0$ & $-75^{\circ}$ \\
	\hline
	 2 & $\Aset_2$ & $1.10$ & $7$ & $\Ga{r}$ & $1$ & $[0.40 ,0.50]$ & $69^{\circ}$ \\
	\hline
	 3 & $\Aset_3$ & $2.30$ & $0$ & $\Ga{t}$ & $[0.30 ,0.40]$ & $1$ & $90^{\circ}$ \\
	\hline
	 4 & $\Aset_4$ & $1.55$ & $59$ & $\Ga{r}$ & $1$ & $[0.85 ,0.95]$ & $-87^{\circ}$ \\
	\hline
	 5 & $\Aset_5$ & $2.00$ & $29$ & $\Ga{t}$ & $[0.60 ,0.70]$ & $1$ & $3^{\circ}$ \\
	\hline
	\end{tabular}
	}

	\caption{Interpolation configurations for online evaluations of the surrogate-based optimisation strategy.}
	\label{tab:dataOnlineInterp}
\end{table}

\begin{figure}[!htb]
	\centering
	\subfigure[Case 1: $\bthetaRef_{\text{opt}}$]{\includegraphics[width=0.3\textwidth]{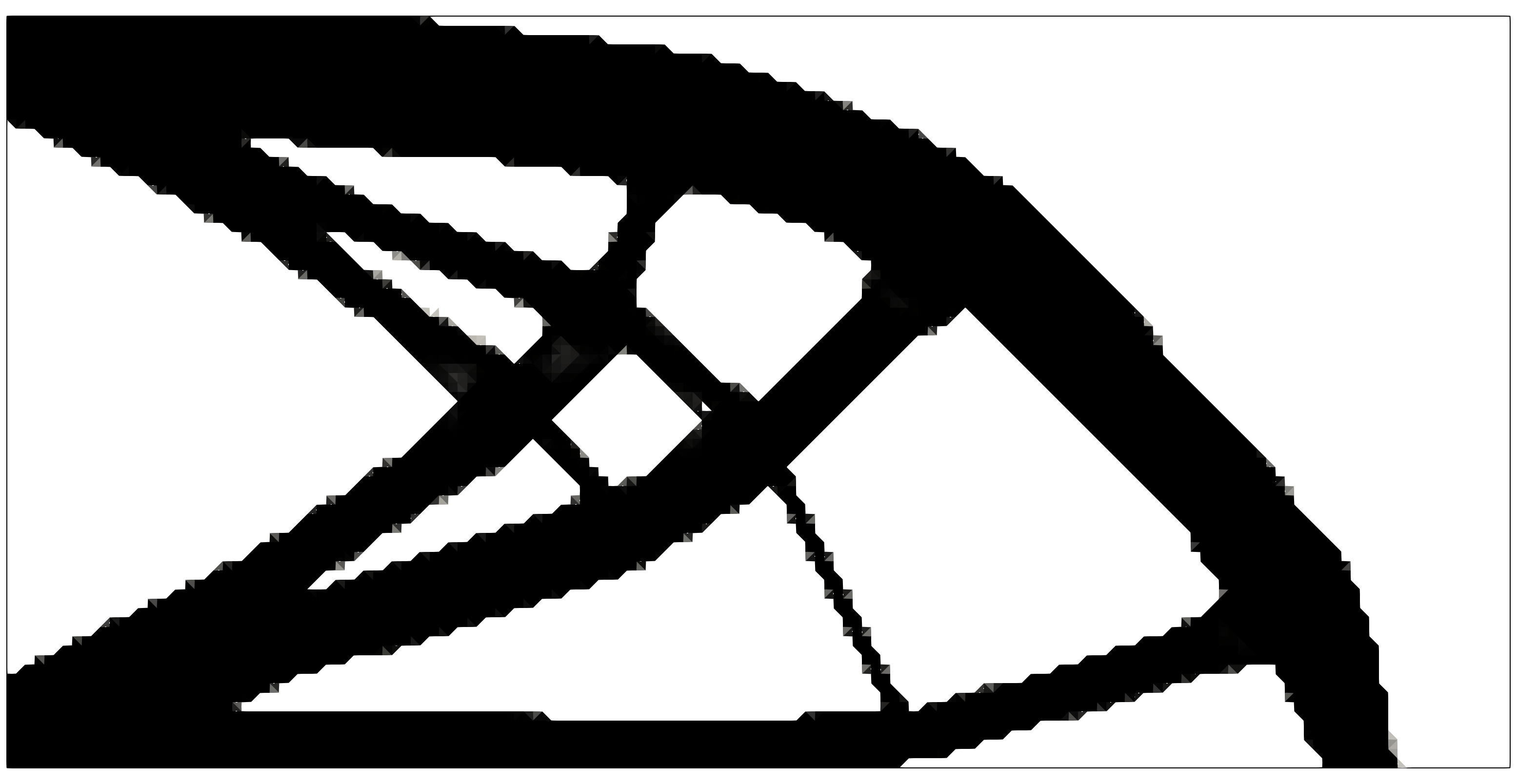}}
	\hspace{5pt}
	\subfigure[Case 1: $\bthetaEta$]{\includegraphics[width=0.3\textwidth]{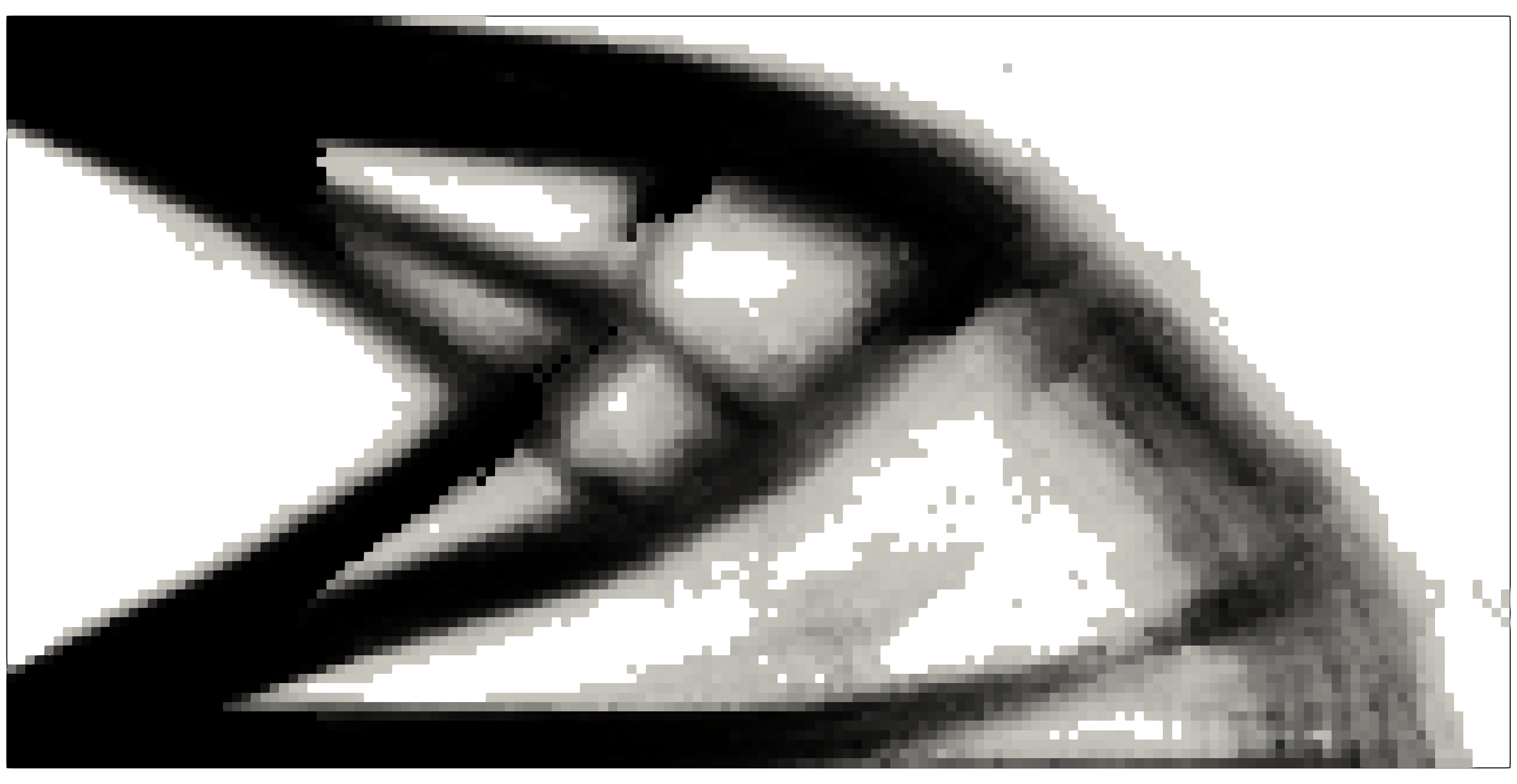}}
	\hspace{5pt}
	\subfigure[Case 1: $\bthetaEta_{\text{opt}}$]{\includegraphics[width=0.3\textwidth]{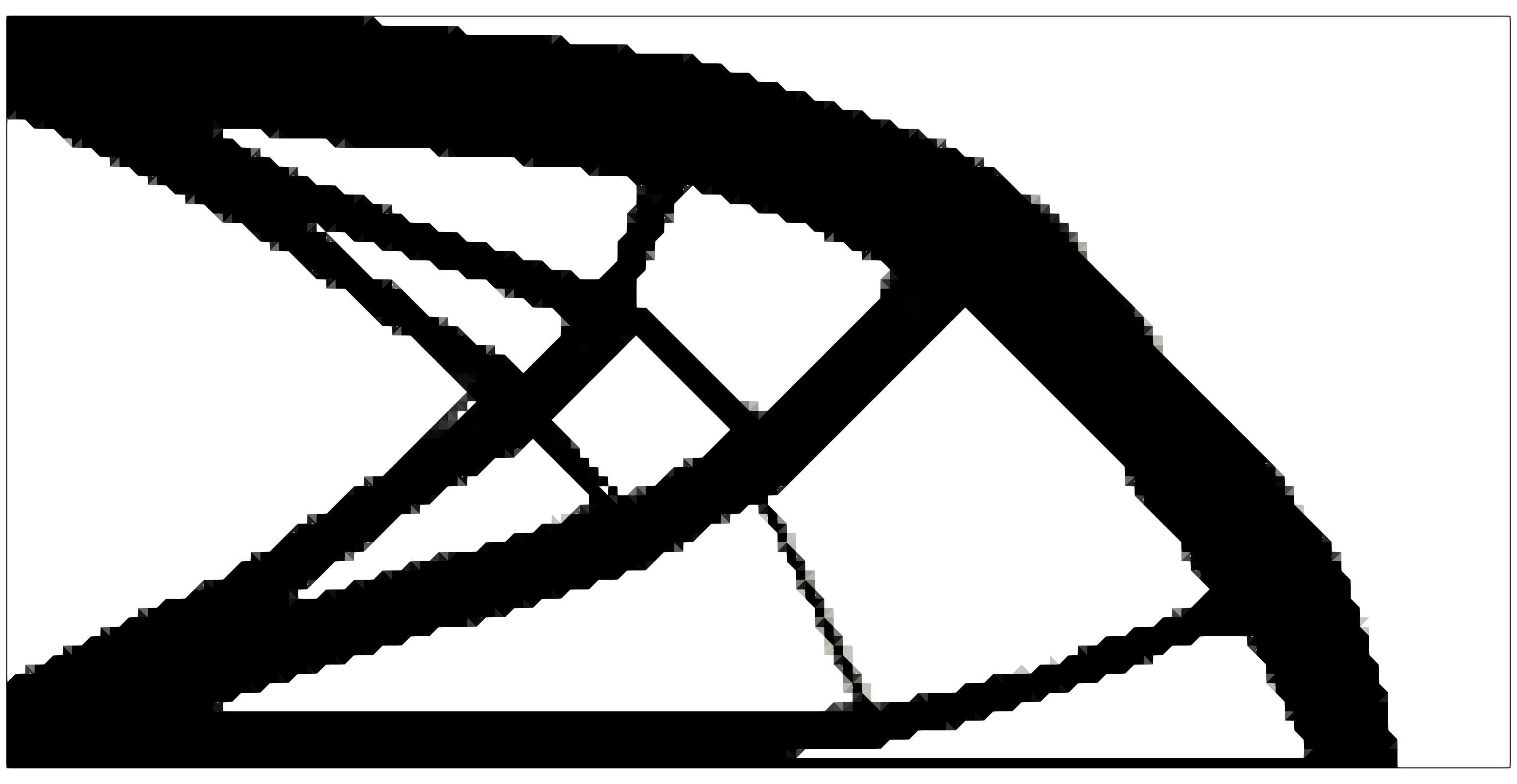}}
	
	\subfigure[Case 2: $\bthetaRef_{\text{opt}}$]{\includegraphics[width=0.3\textwidth]{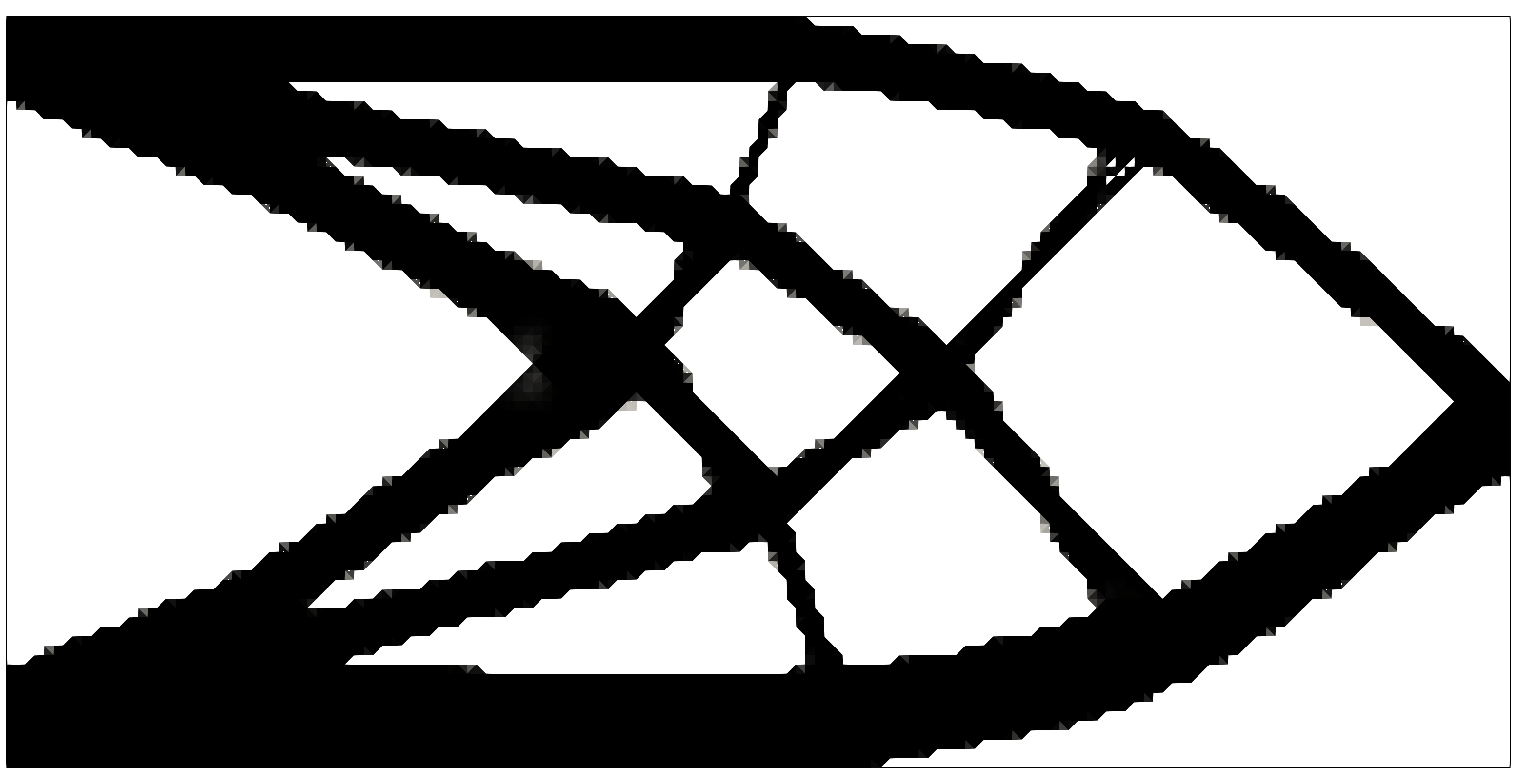}}
	\hspace{5pt}
	\subfigure[Case 2: $\bthetaEta$]{\includegraphics[width=0.3\textwidth]{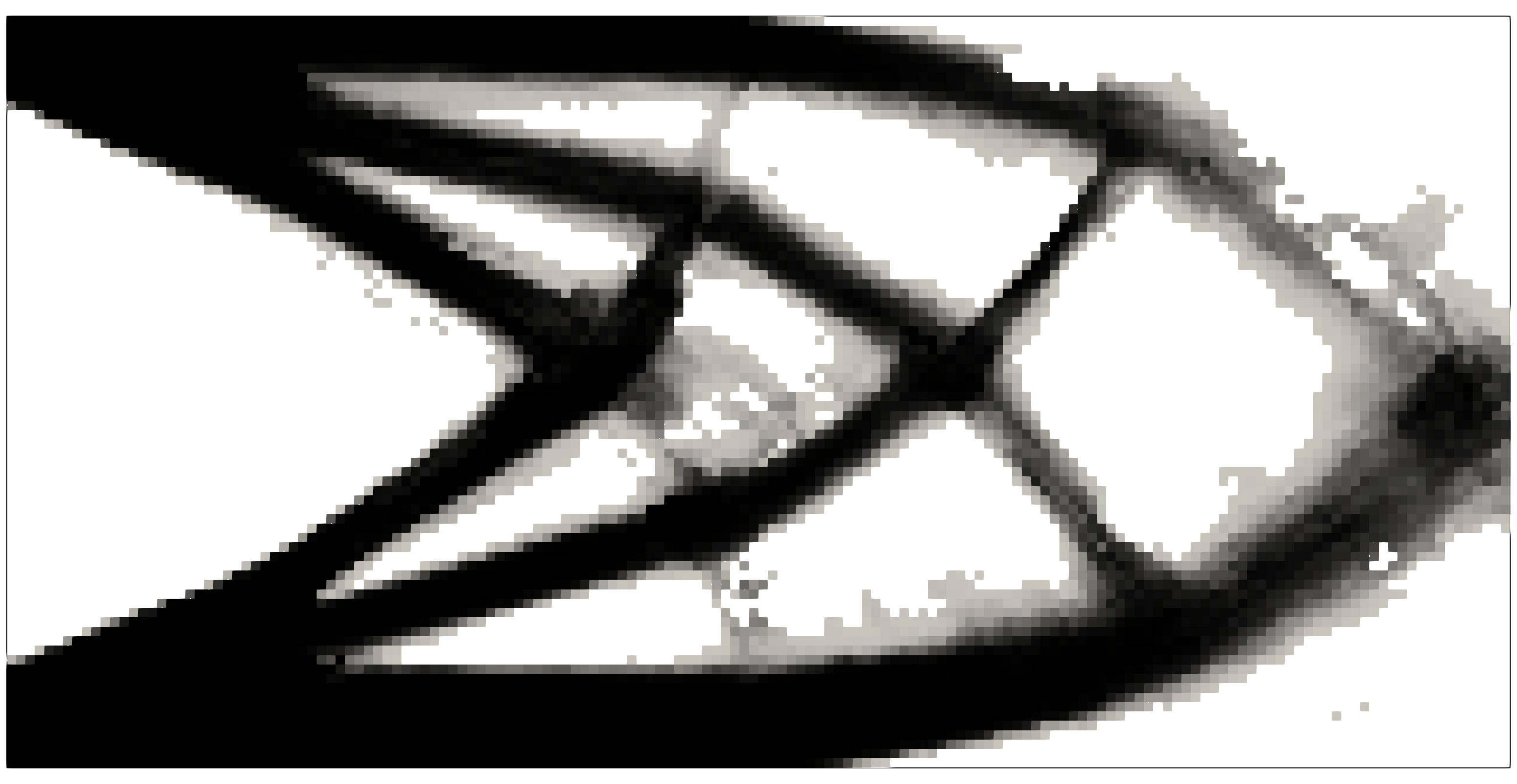}}
	\hspace{5pt}
	\subfigure[Case 2: $\bthetaEta_{\text{opt}}$]{\includegraphics[width=0.3\textwidth]{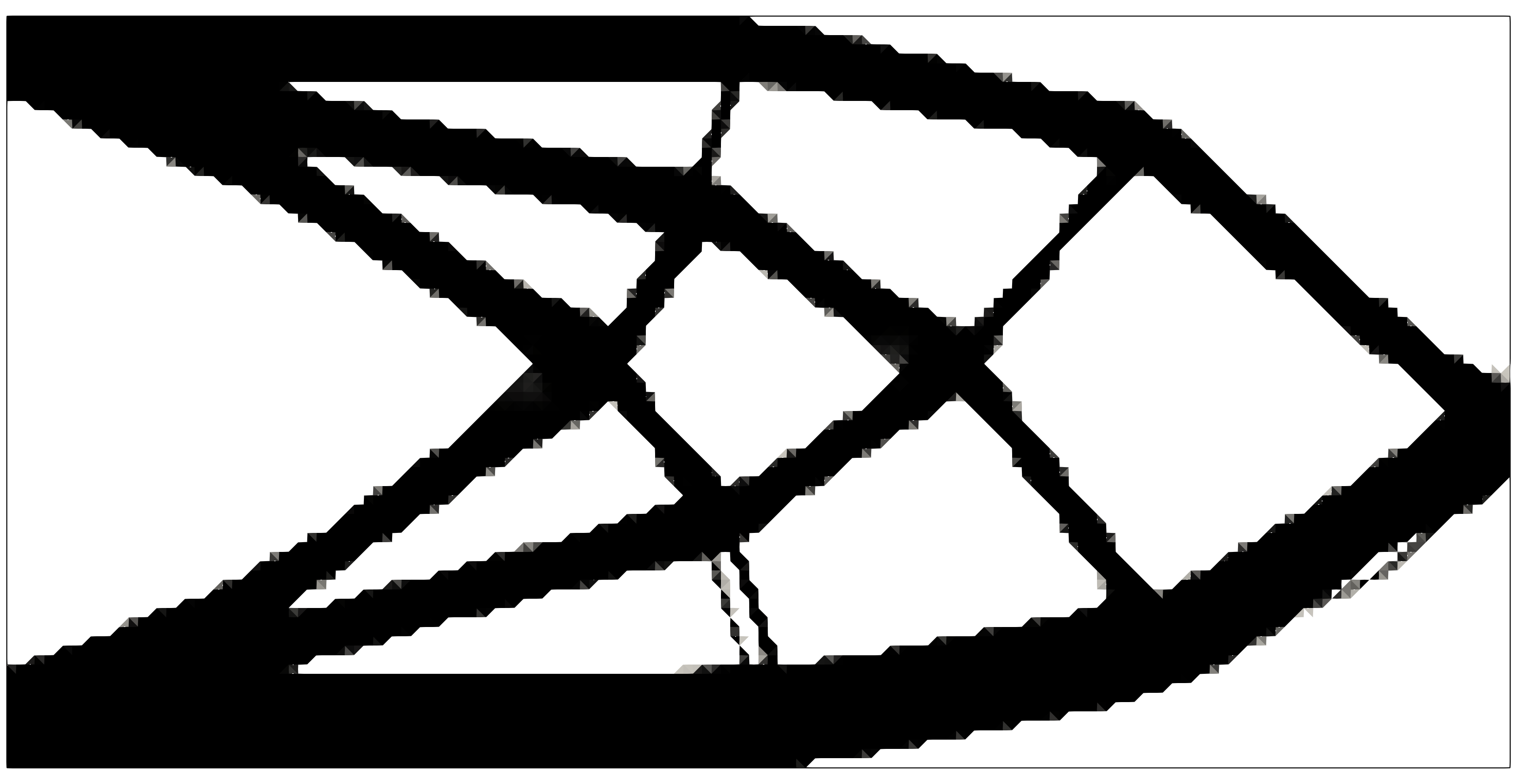}}
	
	\subfigure[Case 3: $\bthetaRef_{\text{opt}}$]{\includegraphics[width=0.3\textwidth]{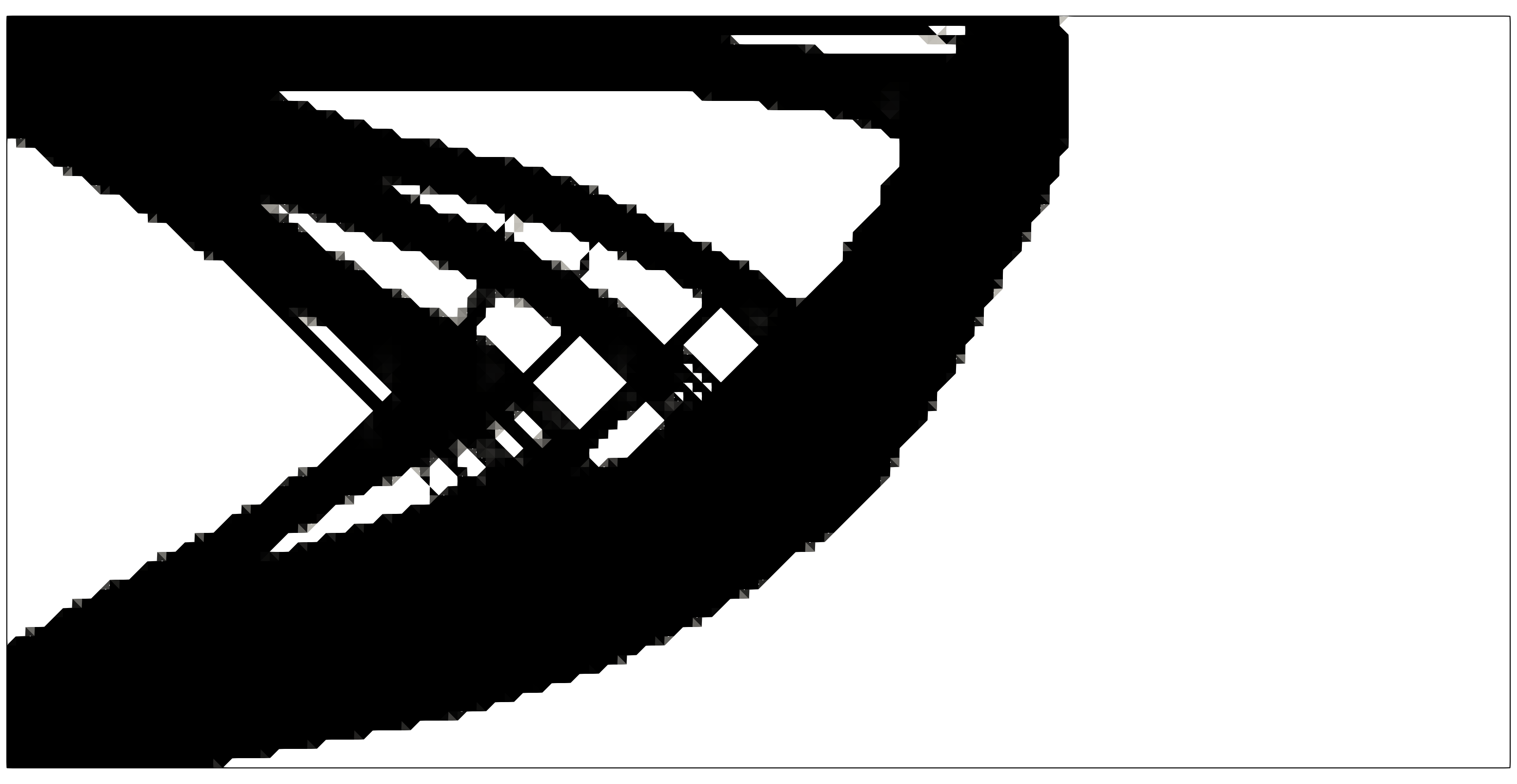}}
	\hspace{5pt}
	\subfigure[Case 3: $\bthetaEta$]{\includegraphics[width=0.3\textwidth]{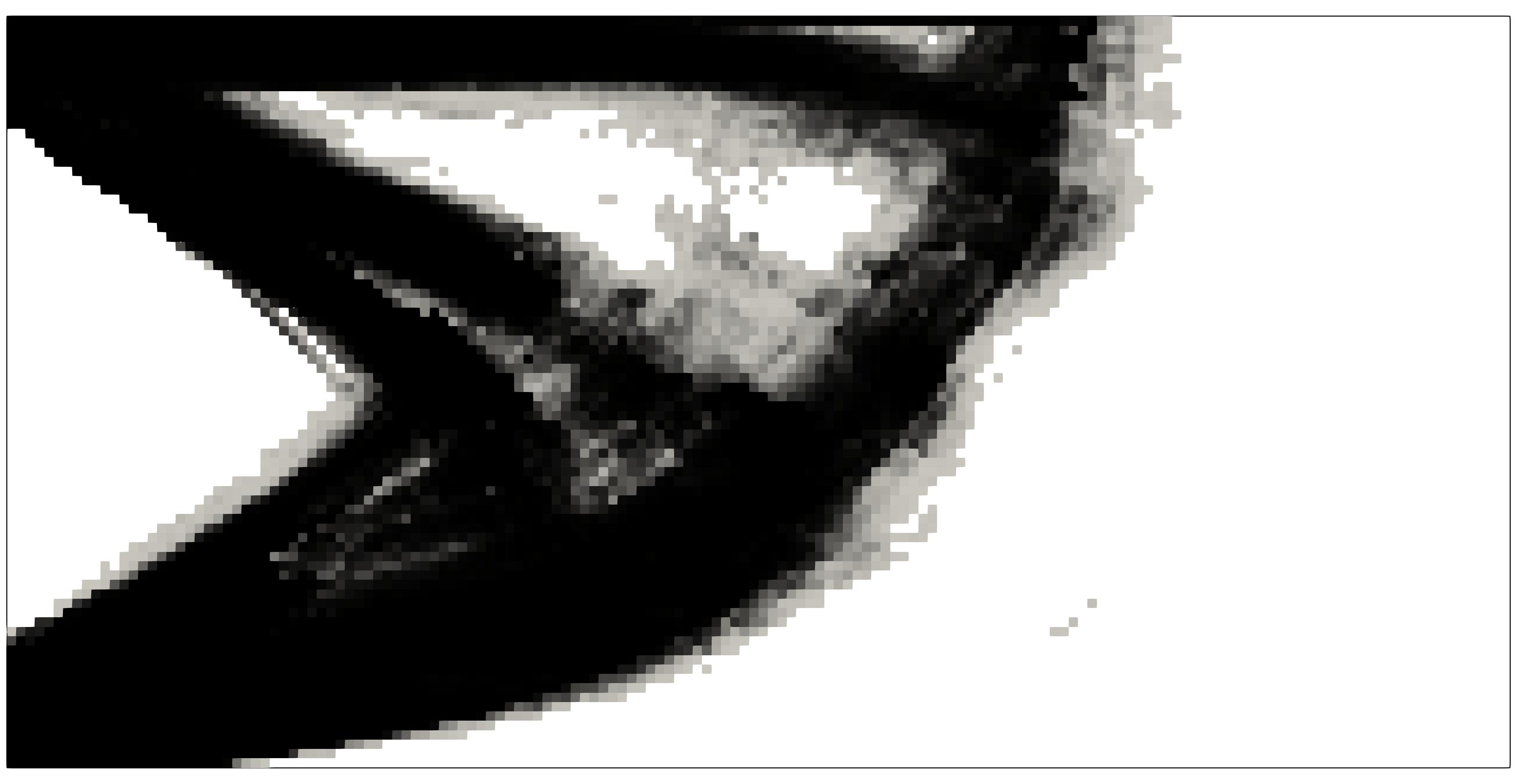}}
	\hspace{5pt}
	\subfigure[Case 3: $\bthetaEta_{\text{opt}}$]{\includegraphics[width=0.3\textwidth]{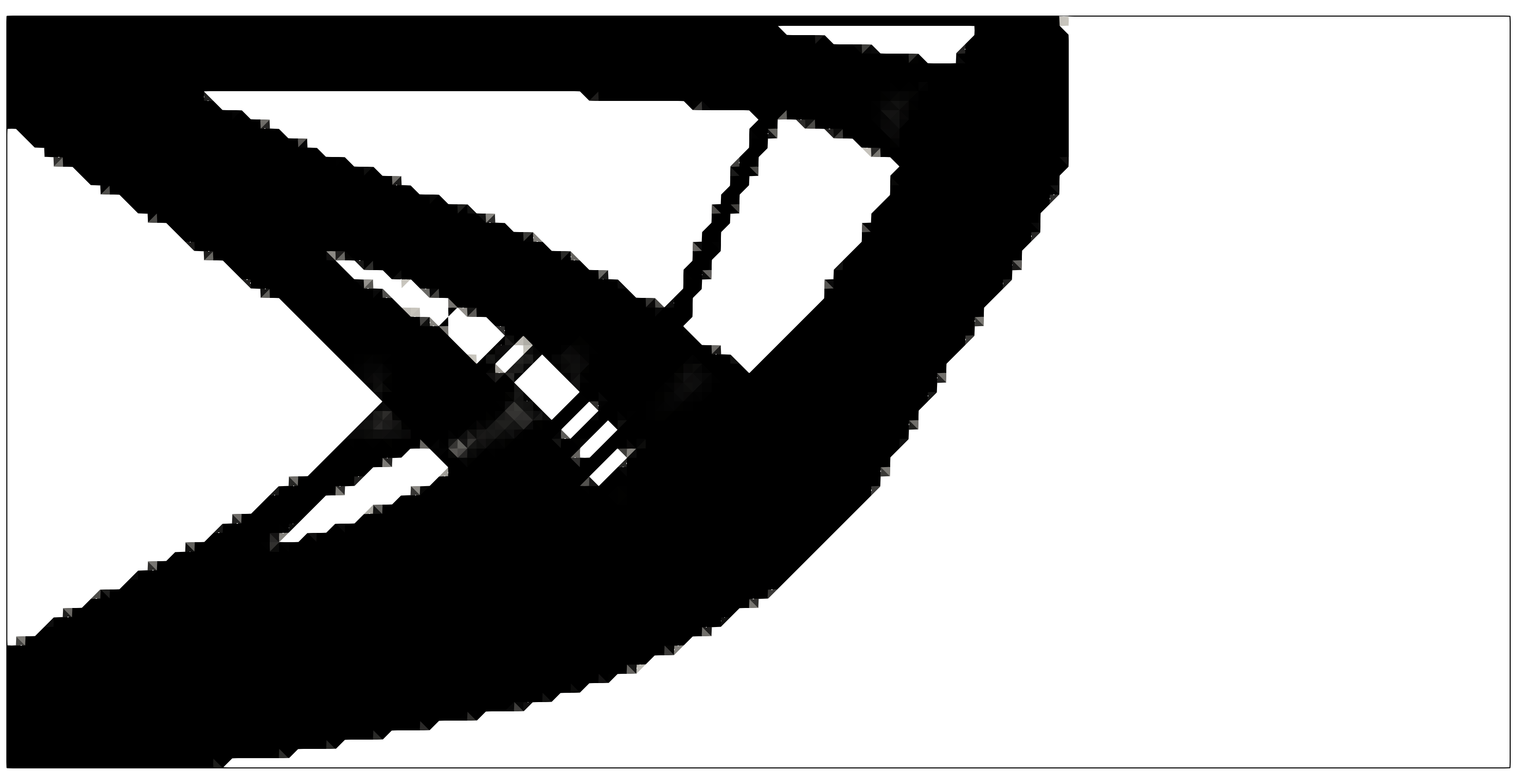}}
	
	\subfigure[Case 4: $\bthetaRef_{\text{opt}}$]{\includegraphics[width=0.3\textwidth]{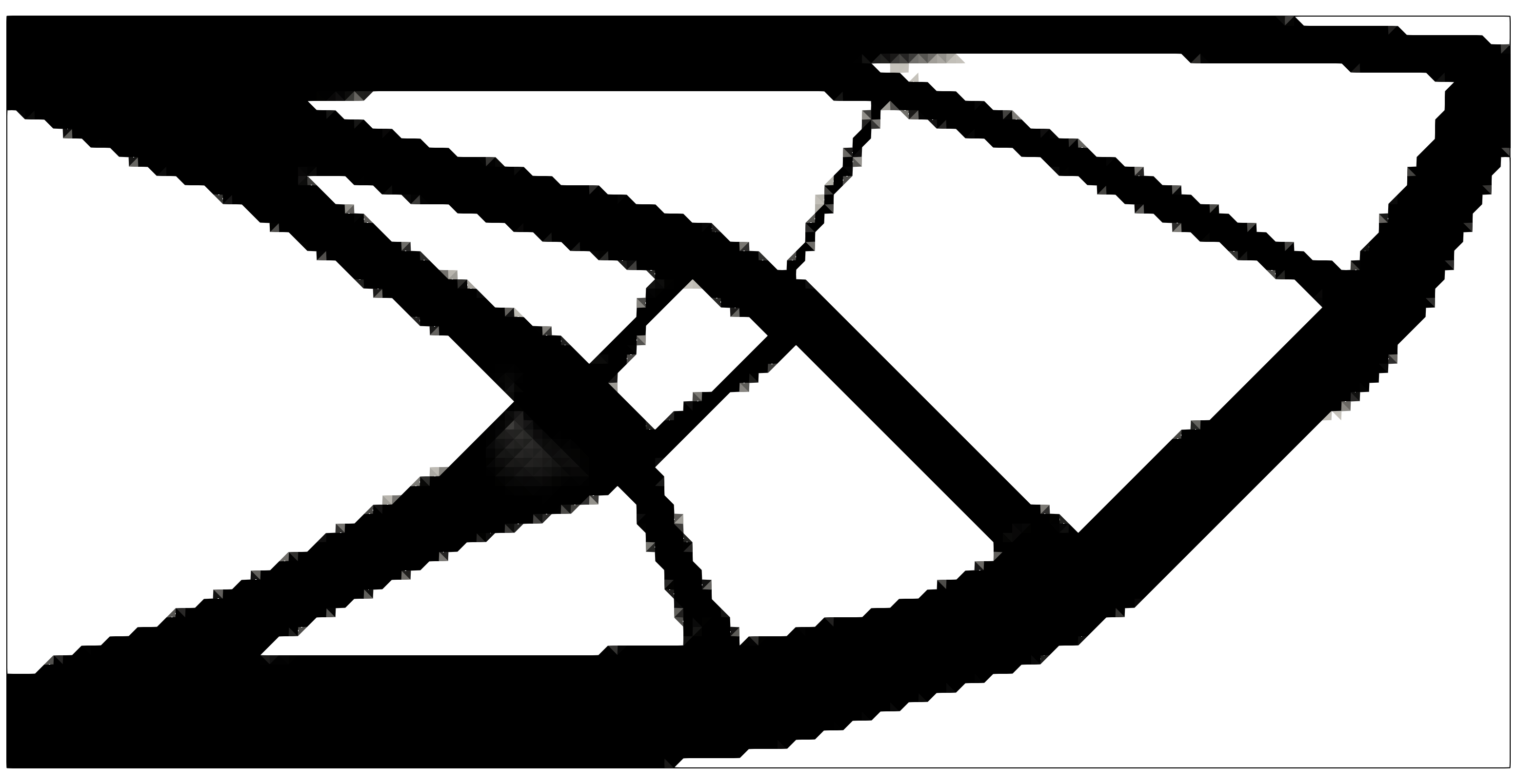}}
	\hspace{5pt}
	\subfigure[Case 4: $\bthetaEta$]{\includegraphics[width=0.3\textwidth]{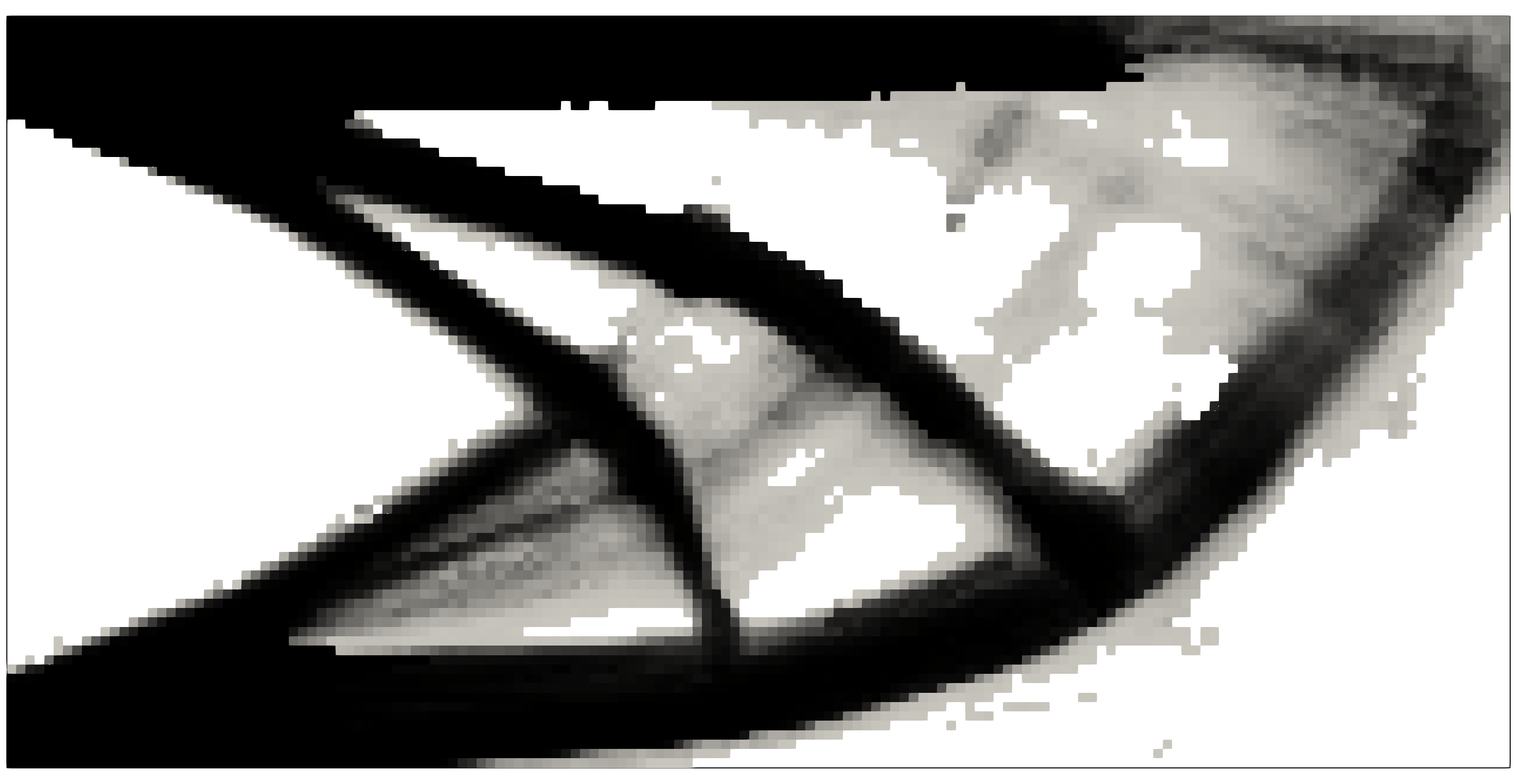}}
	\hspace{5pt}
	\subfigure[Case 4: $\bthetaEta_{\text{opt}}$]{\includegraphics[width=0.3\textwidth]{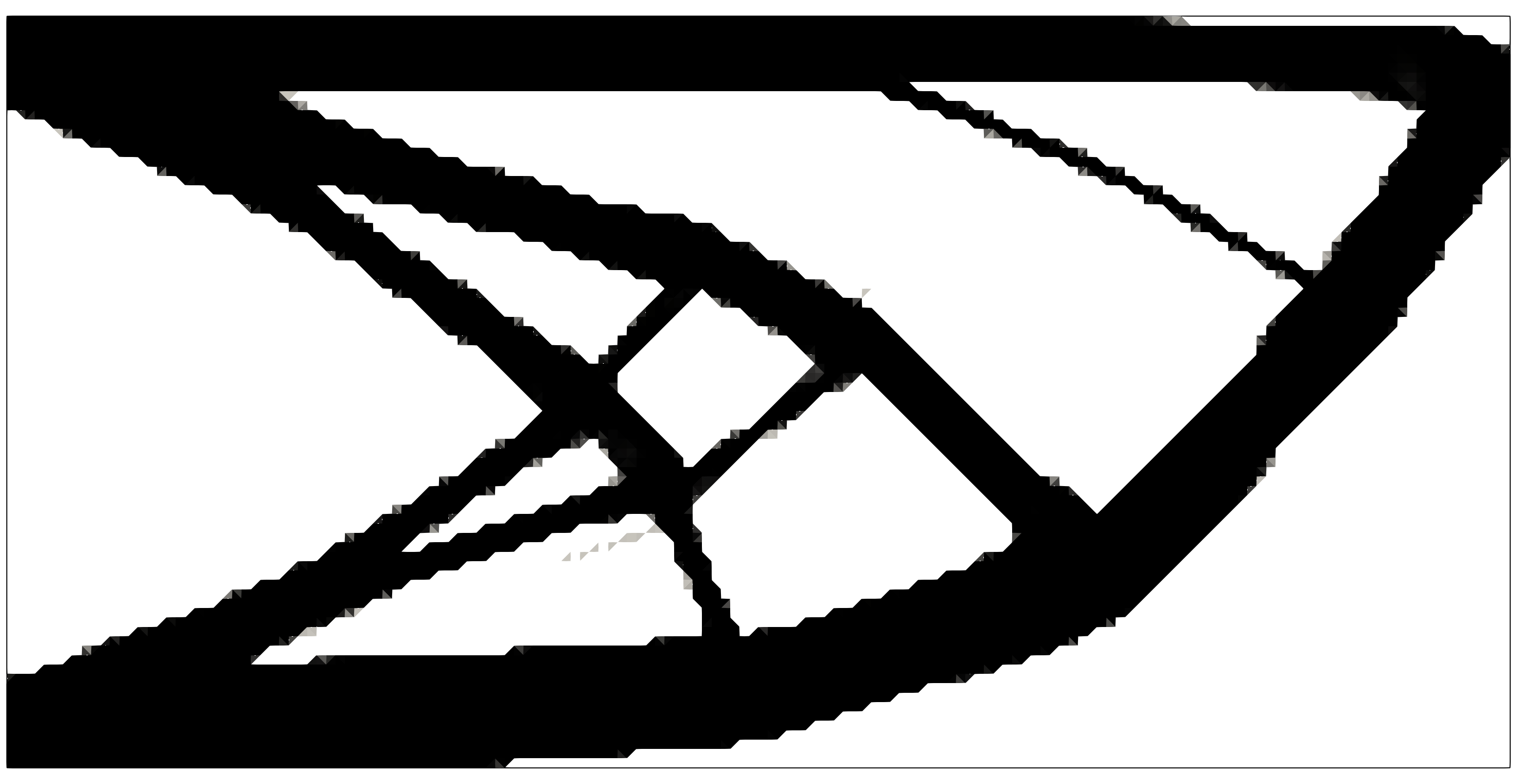}}
	
	\subfigure[Case 5: $\bthetaRef_{\text{opt}}$]{\includegraphics[width=0.3\textwidth]{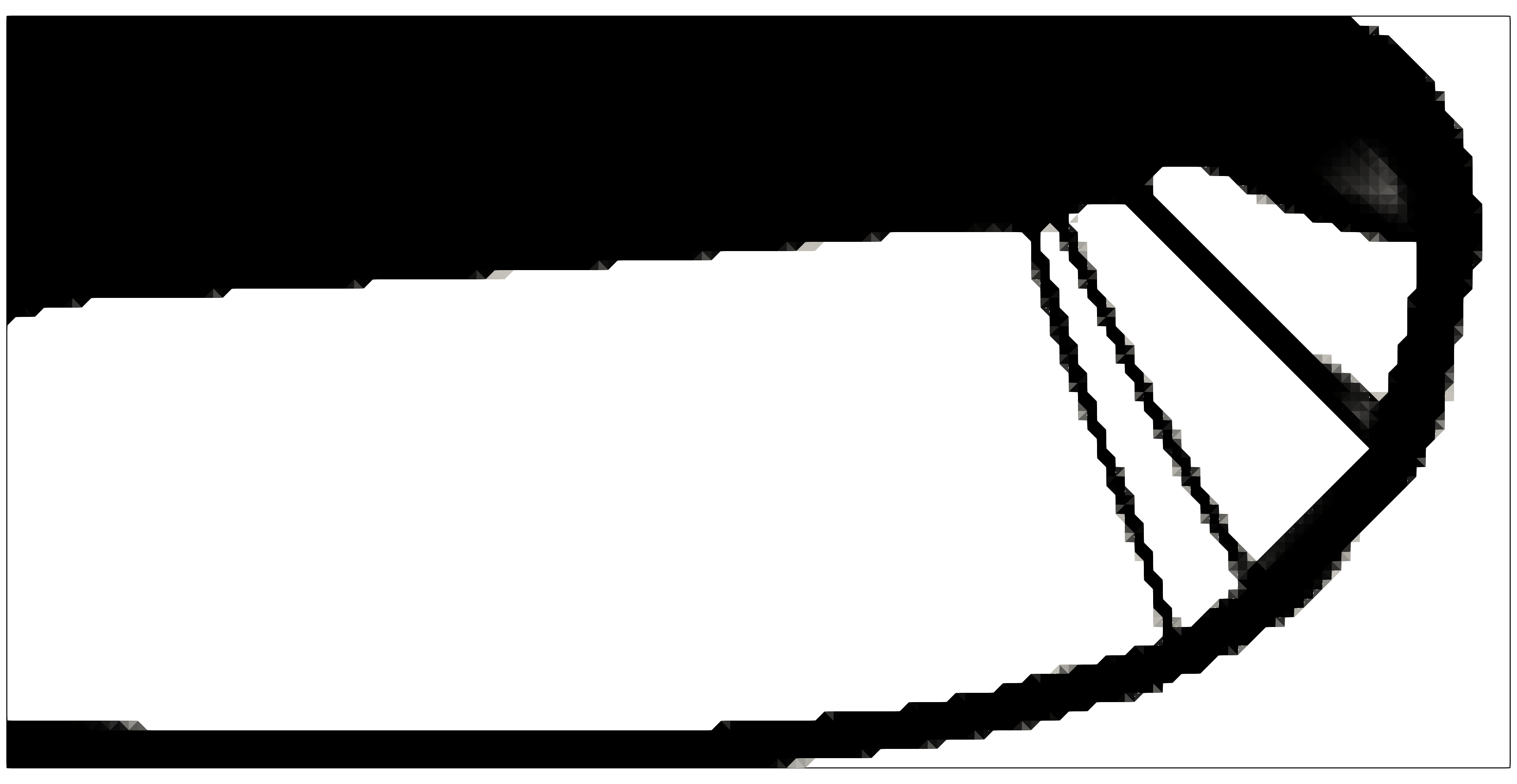}}
	\hspace{5pt}
	\subfigure[Case 5: $\bthetaEta$]{\includegraphics[width=0.3\textwidth]{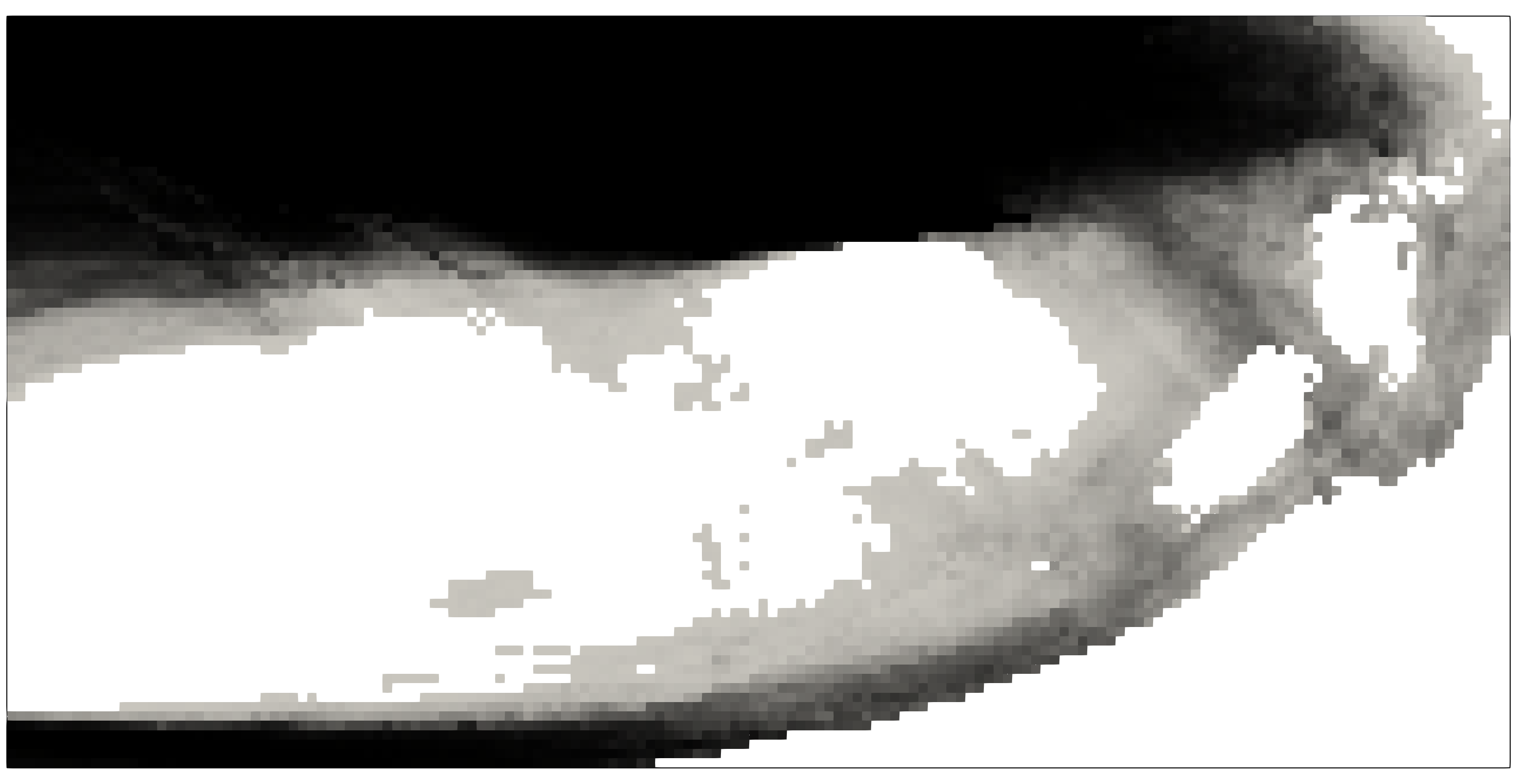}}
	\hspace{5pt}
	\subfigure[Case 5: $\bthetaEta_{\text{opt}}$]{\includegraphics[width=0.3\textwidth]{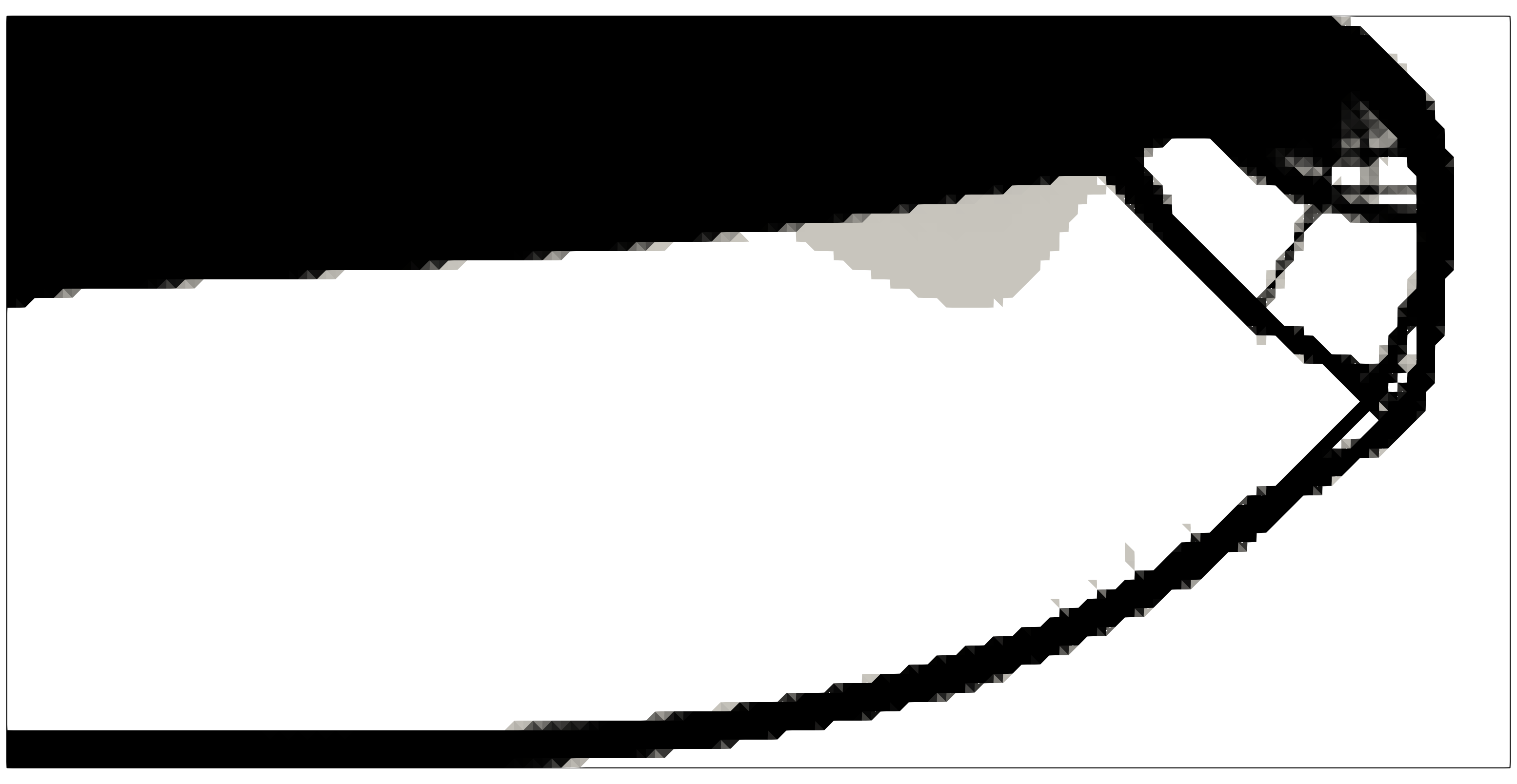}}
	
	\caption{Surrogate-based optimisation for interpolated unseen cases. Left: \emph{ground truth} optimised topology. Centre: \emph{quasi-optimal} topology provided by the surrogate model. Right: topology optimised starting from the \emph{educated} initial guess in the centre.}
	\label{fig:SeedSol}
\end{figure}

Figure~\ref{fig:SeedSol} displays, for the five cases, the \emph{ground truth} solution $\bthetaRef_{\text{opt}}$ computed using algorithm~\ref{alg:topOptHiFi} and not seen by the surrogate model during training, the prediction $\bthetaEta$ performed by $\texttt{FF}_{\!\eta}\texttt{-D}$ for the given set of parameters $(\eta_1,\eta_2)^\top$, and the optimal topology $\bthetaEta_{\text{opt}}$ obtained by means of algorithm~\ref{alg:topOptSurrogate} upon setting the output of the surrogate model as initial condition.

\begin{remark}
Note that the accuracy of the prediction of the surrogate model depends on the precision required during the training phase and on the reliability of the data employed for training.
More precisely, since the \emph{ground truth} optimal topologies have been generated using an optimisation procedure based on a finite element solver (see Section~\ref{sc:DataHomog}), data are polluted by an intrinsic error proportional to the size $h=0.8 \times 10^{-2}$ of the computational mesh.
Hence,  for visualisation purposes, the outputs of the high-fidelity optimisation algorithm, of the surrogate model, and of the surrogate-based optimisation strategy reported in Figure~\ref{fig:SeedSol} are filtered, with all values of material density below $10^{-2}$ being visualised as pure void, (i.e., white). 
\end{remark}

The results in the second column of Figure~\ref{fig:SeedSol} showcase that the surrogate model is able to accurately predict the main mechanical features of the structure, correctly identifying the portions of the boundary where Dirichlet conditions are applied. Despite the position of the Neumann boundary has not been seen during the training procedure, $\texttt{FF}_{\!\eta}\texttt{-D}$ leverages the information on the neighbouring cases explored during training and validation to interpolate the location of the force. Similarly, the most relevant features of the topology are identified by the surrogate model which lacks sufficient resolution to fully describe the structure.

It is worth noticing that the distributions $\bthetaEta$ are qualitatively \emph{close} to the reference $\bthetaRef_{\text{opt}}$.  Indeed,  the third column of Figure~\ref{fig:SeedSol} displays that fully-resolved optimal topologies can be achieved by the surrogate-based algorithm when $\bthetaEta$ is employed as \emph{educated} initial condition. The results show excellent agreement with the reference topologies, accurately reproducing the main mechanical features, while eliminating most of the spurious intermediate material distributions represented in grey.

\begin{figure}[!htb]
	\centering
	\subfigure[High-fidelity compliance]{\includegraphics[width=0.48\textwidth]{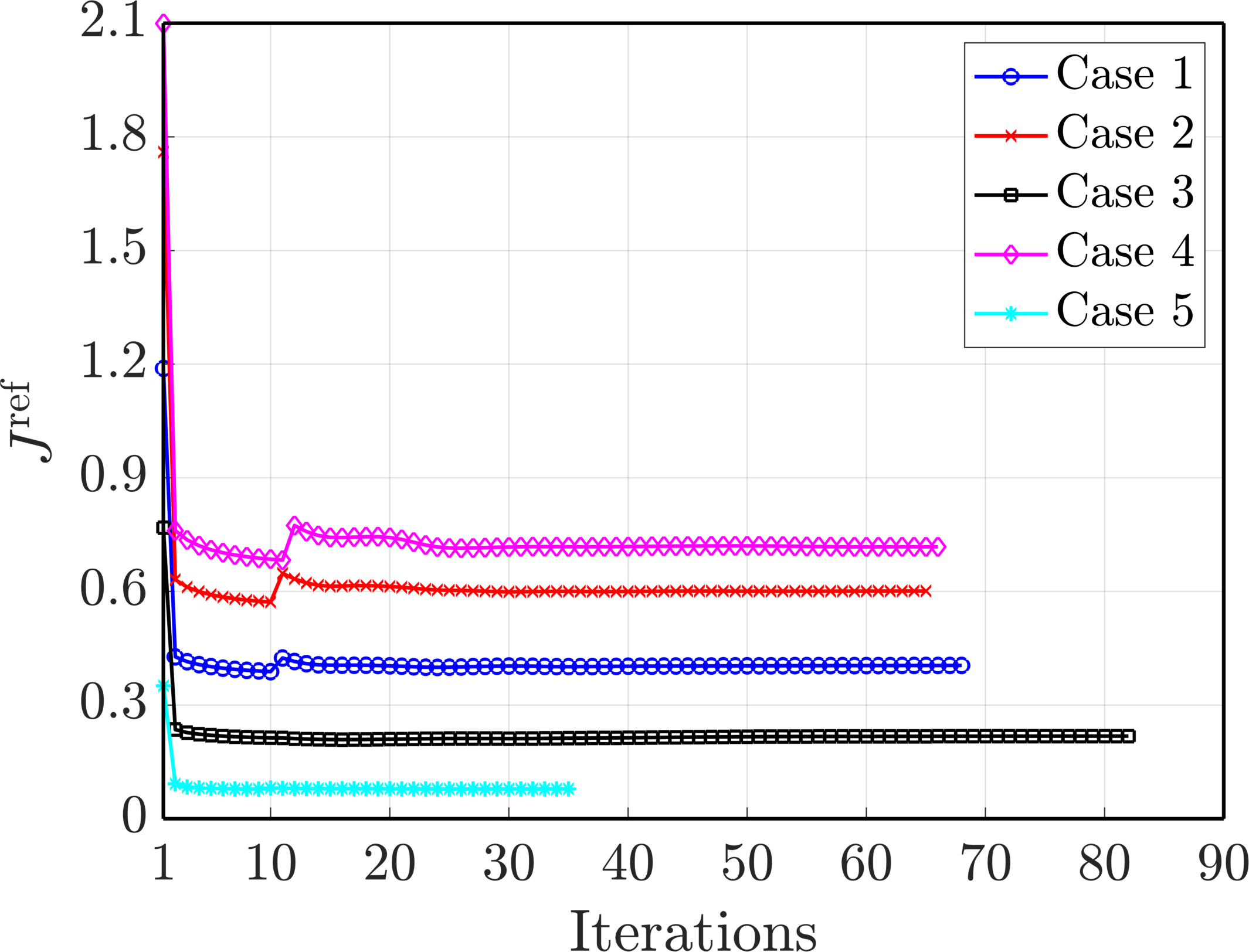}\label{fig:SeedEvolutionA}}
	\hspace{5pt}
	\subfigure[Surrogate compliance]{\includegraphics[width=0.48\textwidth]{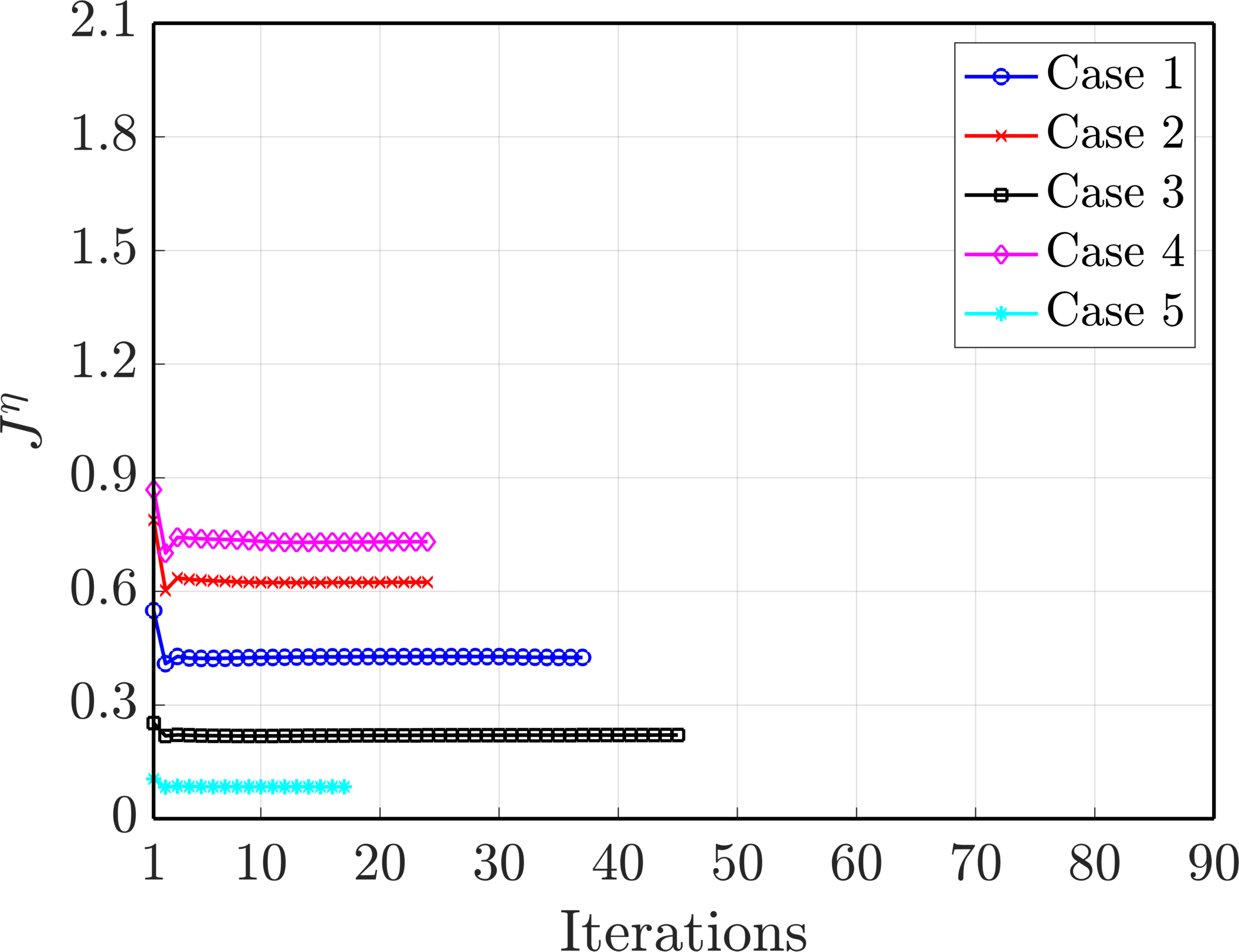}\label{fig:SeedEvolutionB}}
	
	\subfigure[High-fidelity volume fraction]{\includegraphics[width=0.48\textwidth]{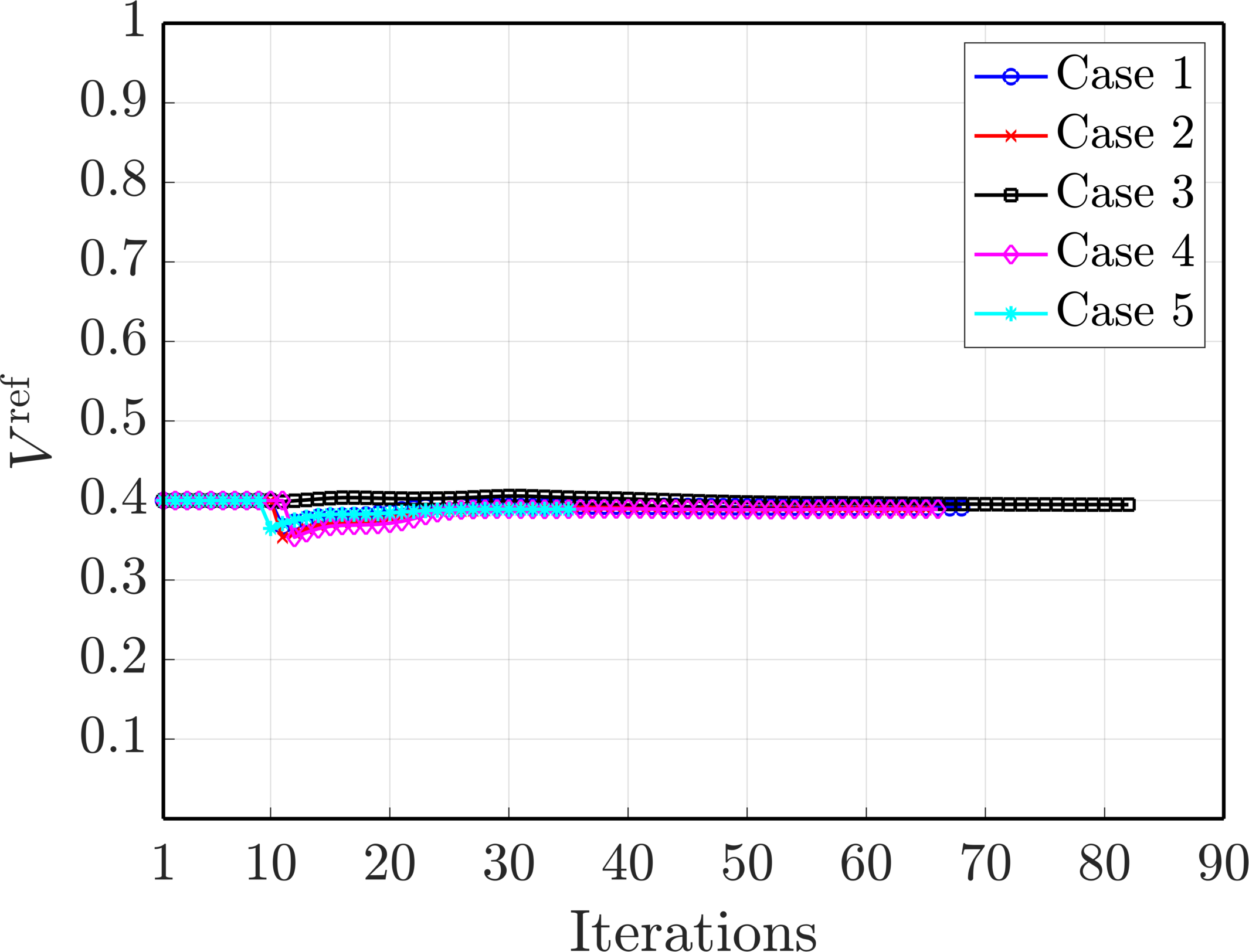}\label{fig:SeedEvolutionC}}
	\hspace{5pt}
	\subfigure[Surrogate volume fraction]{\includegraphics[width=0.48\textwidth]{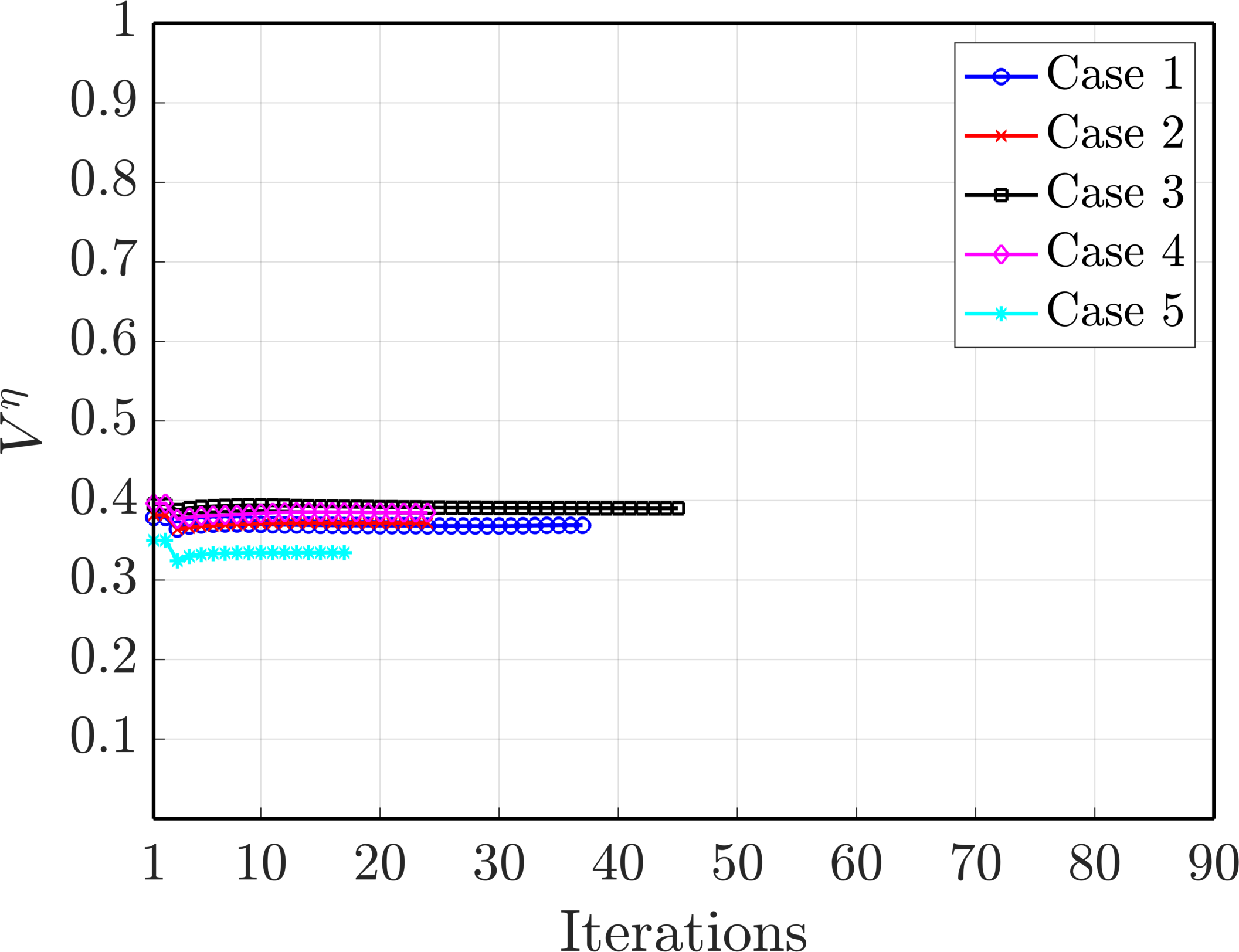}\label{fig:SeedEvolutionD}}
	
	\caption{Evolution of the compliance (top) and the volume fraction (bottom) computed using the high-fidelity (left) and surrogate-based (right) algorithms for the interpolated unseen cases. }
	\label{fig:SeedEvolution}
\end{figure}

To quantitatively assess the accuracy and reliability of the newly optimised topologies $\bthetaEta_{\text{opt}}$, Figure~\ref{fig:SeedEvolution} reports the evolution of the compliance~\eqref{eq:compliance} and the volume fraction~\eqref{eq:volFrac} using algorithm~\ref{alg:topOptSurrogate} with the \emph{educated} initial guess, compared to the high-fidelity algorithm~\ref{alg:topOptHiFi} starting with a uniform material distribution.
The initial guess provided by the surrogate model $\texttt{FF}_{\!\eta}\texttt{-D}$ yields initial values of the compliance significantly closer to the optimum than the standard high-fidelity ones, as visible by comparing Figure~\ref{fig:SeedEvolutionA} and~\ref{fig:SeedEvolutionB}.  Moreover, it is straightforward to observe that the number of iterations is significantly reduced when the \emph{educated} initial guess is employed.
Concerning the volume fraction of the structure (Figure~\ref{fig:SeedEvolutionC} and~\ref{fig:SeedEvolutionD}), it appears that the surrogate model $\texttt{FF}_{\!\eta}\texttt{-D}$ introduces a small error at convergence,  slightly devitating from the target value of $0.4$.

Table~\ref{tab:interpJV} reports the values $\Jeta[0]$ and $\Veta[0]$ of the compliance and the volume fraction provided as output by the evaluation of the surrogate model $\texttt{FF}_{\!\eta}\texttt{-D}$ and the corresponding optimal values $\Jeta[\text{opt}]$ and $\Veta[\text{opt}]$ (respectively, $\Jref[\text{opt}]$ and $\Vref[\text{opt}]$) attained using the surrogate-based (respectively, high-fidelity) topology optimisation algorithm.
\begin{table}[!htb]
	\centering
	\begin{tabular}{| c || c | c | c || c | c | c || c | c |}
	\hline
	Case & $\Jeta[0]$ & $\Jeta[\text{opt}]$ & $\Jref[\text{opt}]$ & $\Veta[0]$ & $\Veta[\text{opt}]$ & $\Vref[\text{opt}]$ & $\niter[\eta]$ & $\niter[\text{ref}]$ \\
	\hline
	 1 & $0.550$ & $0.426$ & $0.405$ & $0.378$ & $0.369$ & $0.391$& $37$ & $68$\\
	\hline
	 2 & $0.788$ & $0.624$ & $0.601$ & $0.381$ & $0.372$ & $0.389$ & $24$ & $65$ \\
	\hline
	 3 & $0.253$ & $0.221$ & $0.218$ & $0.395$ & $0.390$ & $0.394$ & $45$ & $82$ \\
	\hline
	 4 & $0.868$ & $0.731$ & $0.718$ & $0.396$ & $0.384$ & $0.389$ & $24$ & $66$ \\
	\hline
	 5 & $0.105$ & $0.084$ & $0.079$ & $0.350$ & $0.334$ & $0.389$ & $17$ & $35$ \\
	\hline
	\end{tabular}

	\caption{Initial, optimal, and reference values of the compliance and the volume fraction computed using the high-fidelity and surrogate-based optimisation algorithms and corresponding number of iterations for the interpolated unseen cases.}
	\label{tab:interpJV}
\end{table}
The results show that the surrogate model initially features an error with respect to the reference compliance $\Jref[\text{opt}]$ between $15.78\%$ and $35.86\%$,  which quickly decreases until attaining accuracies ranging from $1.25\%$ to $7.12\%$.
Moreover,  the number of iterations required for convergence is significantly lowered, with an average reduction of $54\%$ of iterations and peak computational gains of $64\%$ with respect to the high-fidelity strategy.

\subsection{Generalisation capabilities of the surrogate models}
\label{sc:Generalisation}

\begin{figure}[!htb]
	\centering	
	\subfigure[Data $\Eset_1$]{\includegraphics[width=0.14\textwidth]{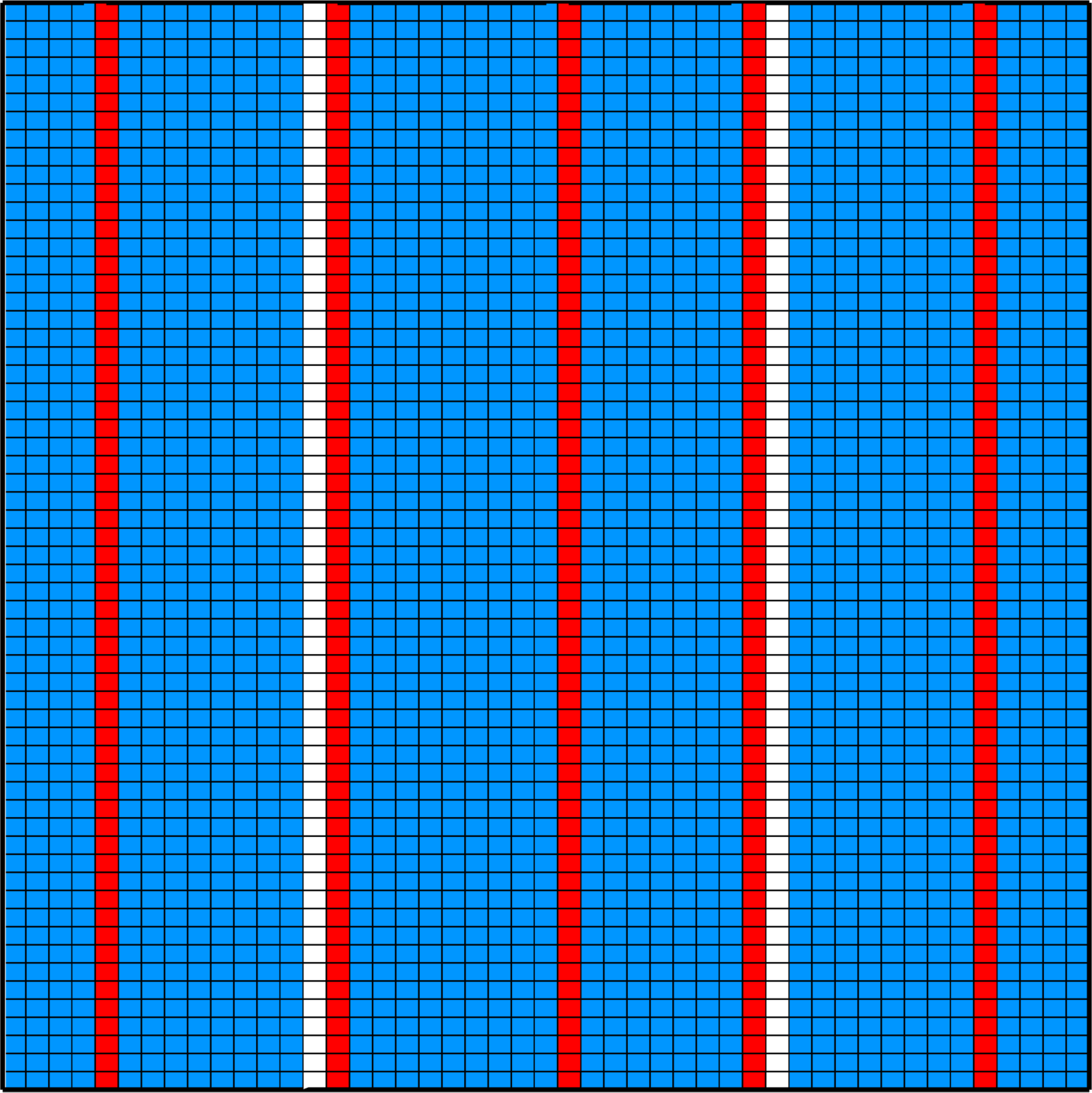}}
	\hspace{2pt}
	\subfigure[Data $\Eset_2$]{\includegraphics[width=0.14\textwidth]{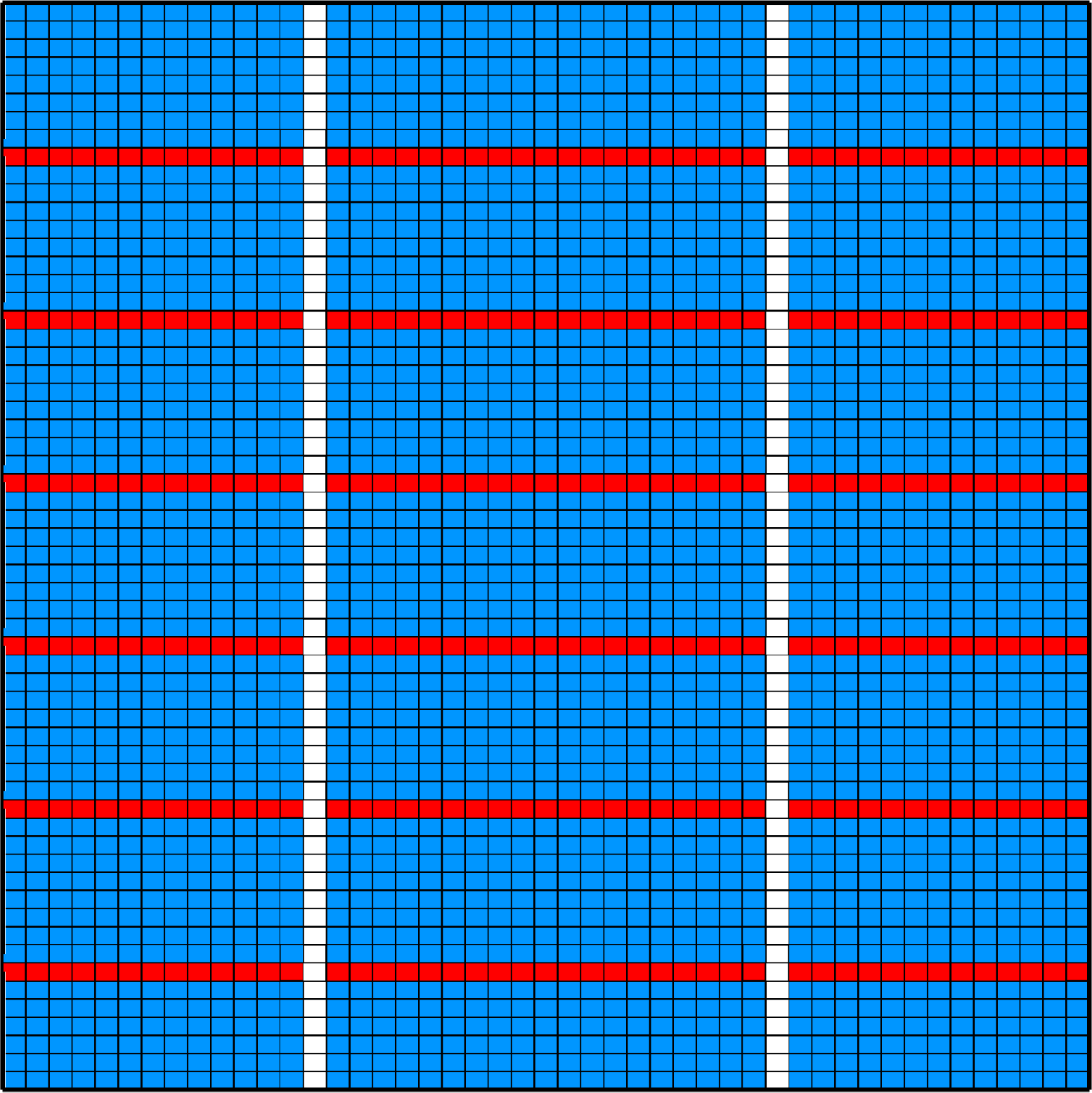}}
	\hspace{2pt}
	\subfigure[Data $\Eset_3$]{\includegraphics[width=0.14\textwidth]{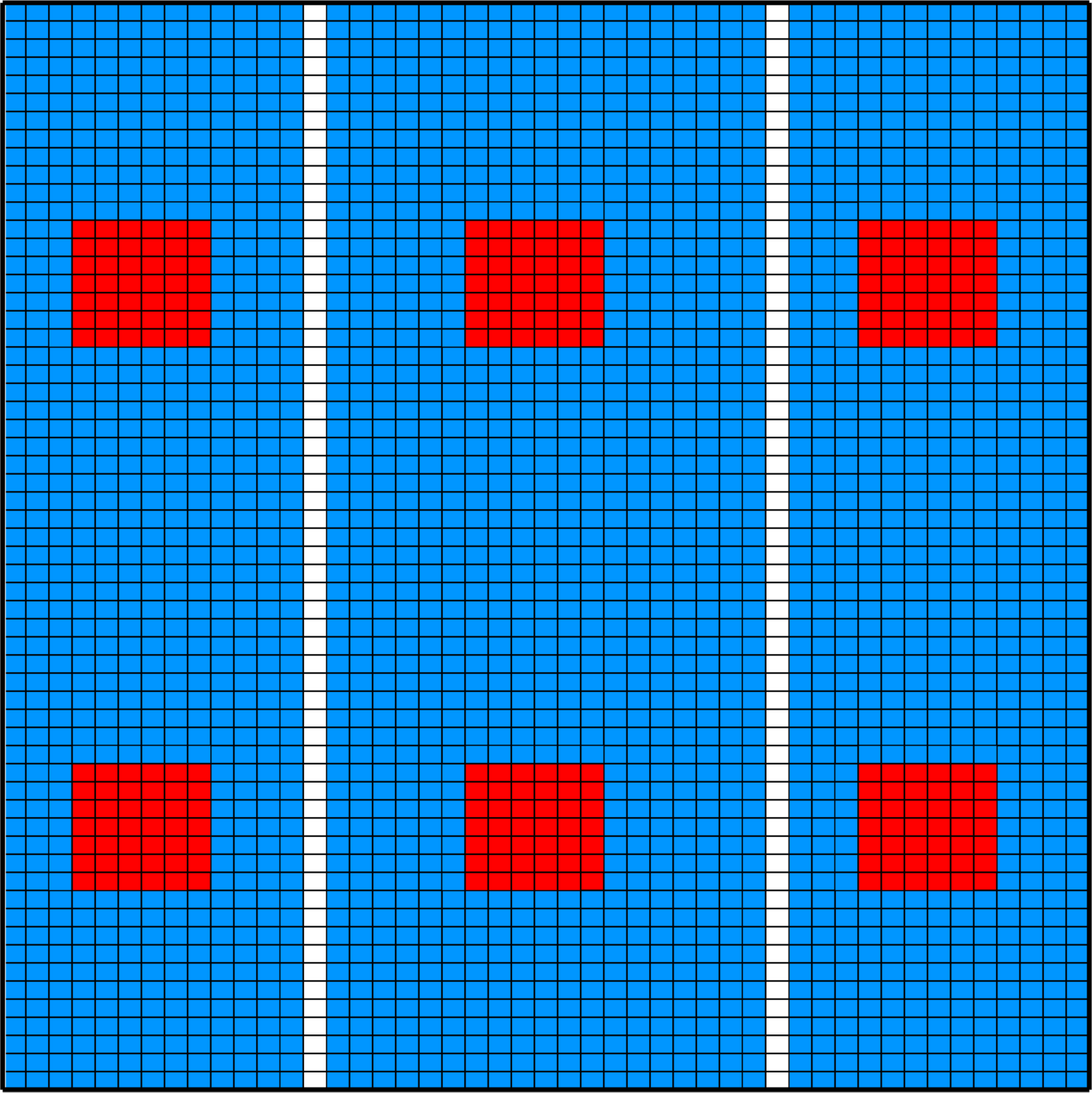}}
	\hspace{2pt}
	\subfigure[Data $\Eset_4$]{\includegraphics[width=0.14\textwidth]{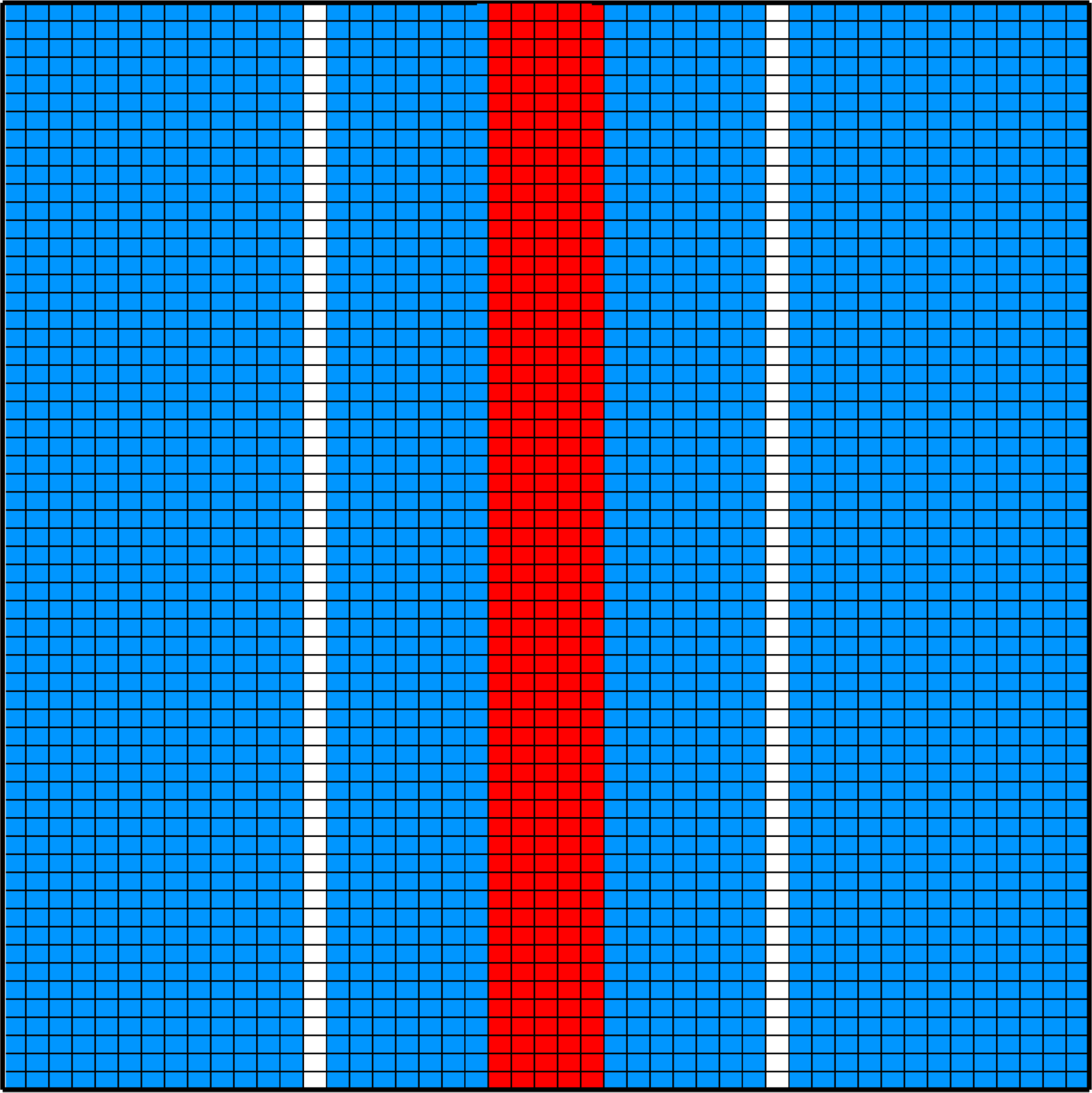}}
	\hspace{2pt}
	\subfigure[Data $\Eset_5$]{\includegraphics[width=0.14\textwidth]{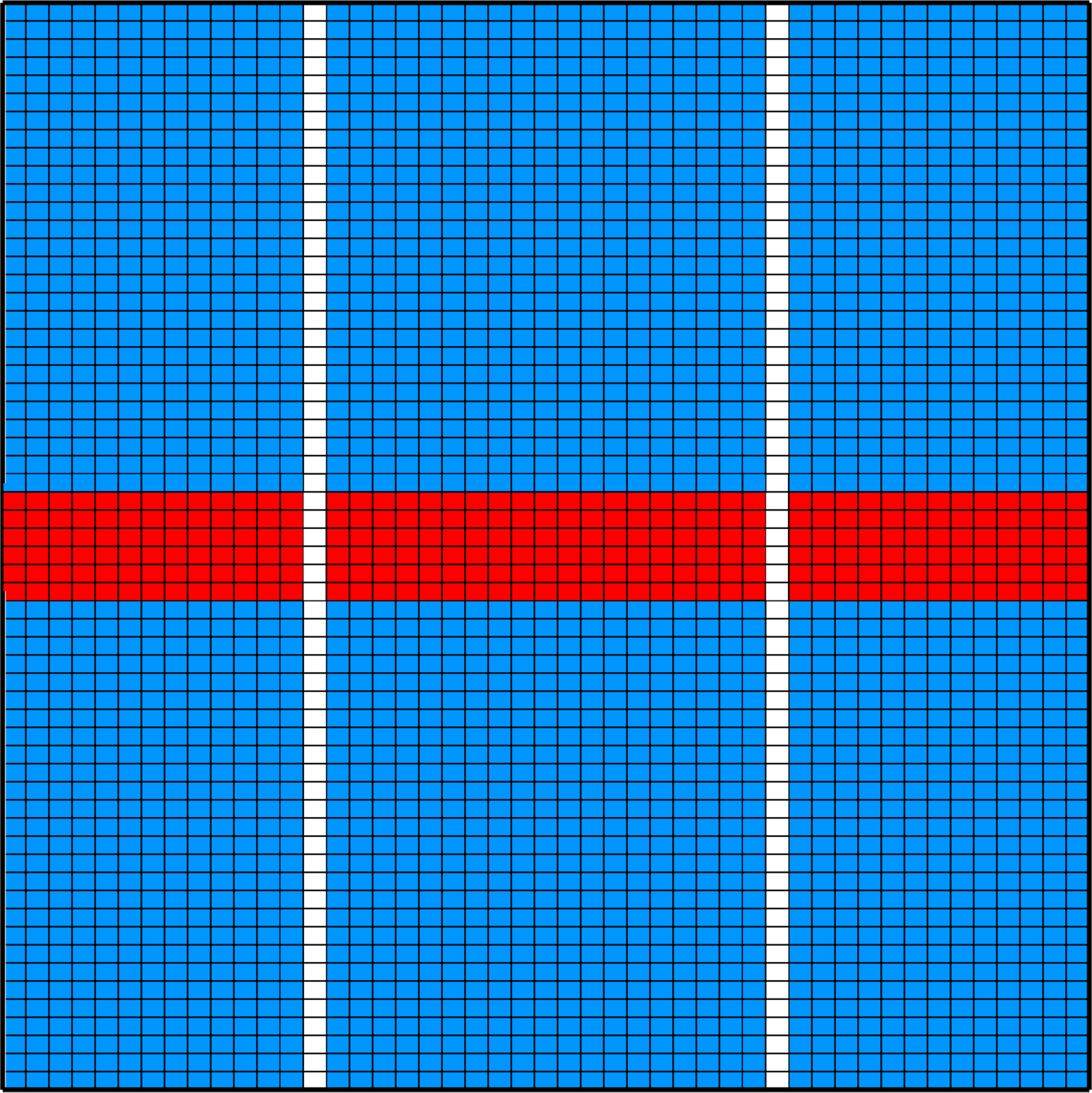}}
	\hspace{2pt}
	\subfigure[Data $\Eset_6$]{\includegraphics[width=0.14\textwidth]{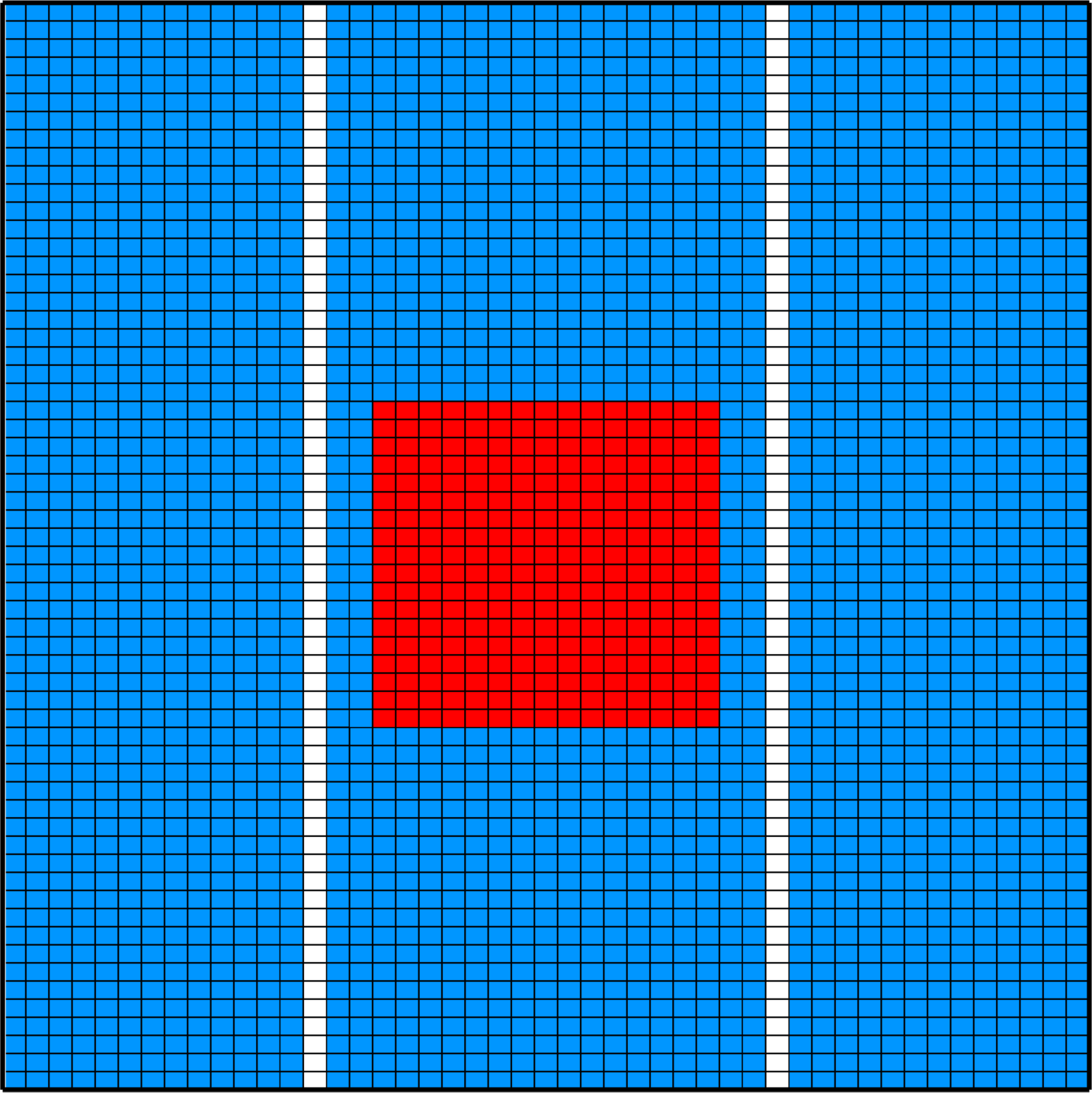}}
	
	\subfigure[Data $\Eset_7$]{\includegraphics[width=0.14\textwidth]{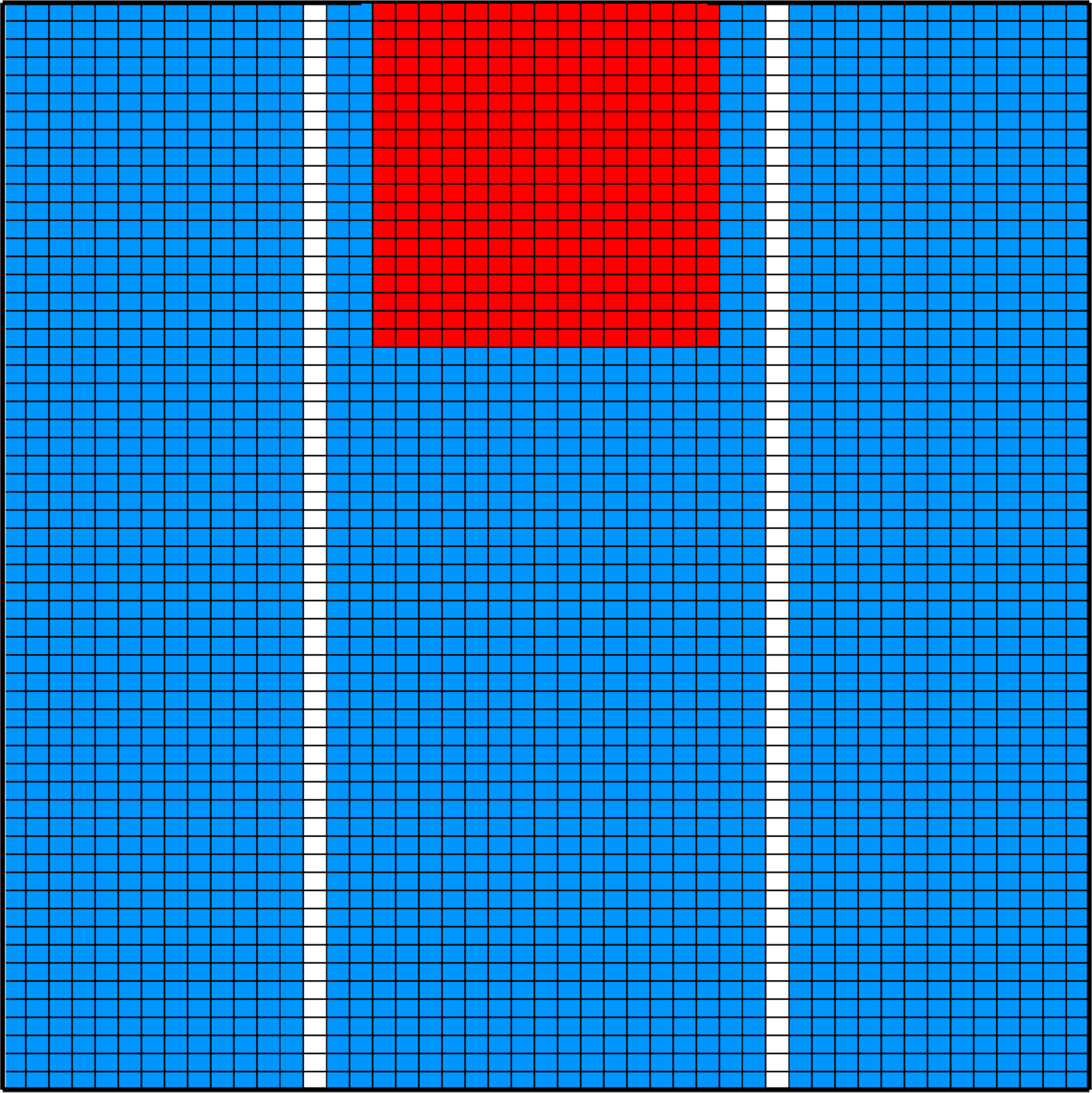}}
	\hspace{2pt}
	\subfigure[Data $\Eset_8$]{\includegraphics[width=0.14\textwidth]{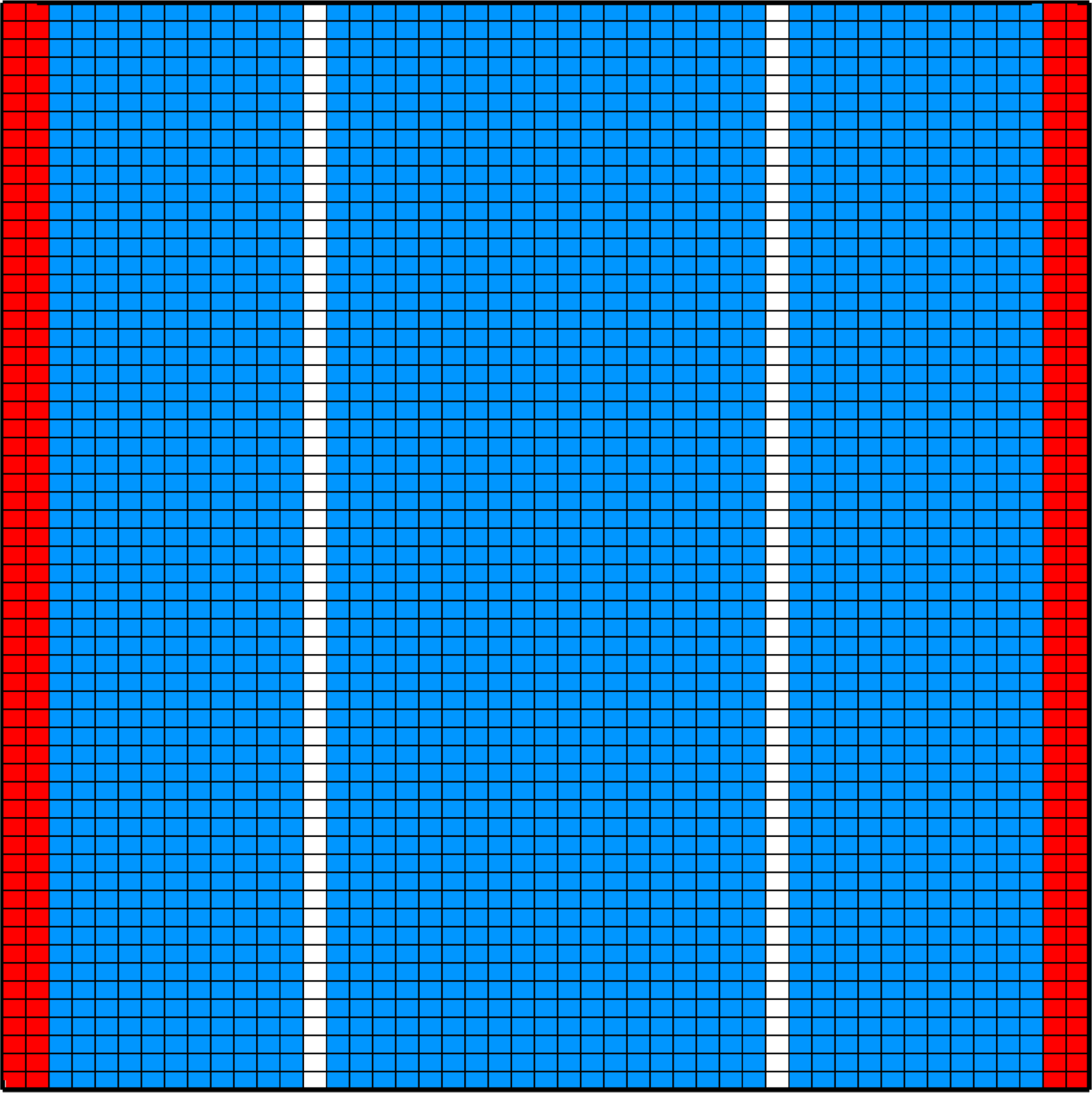}}
	\hspace{2pt}
	\subfigure[Data $\Eset_9$]{\includegraphics[width=0.14\textwidth]{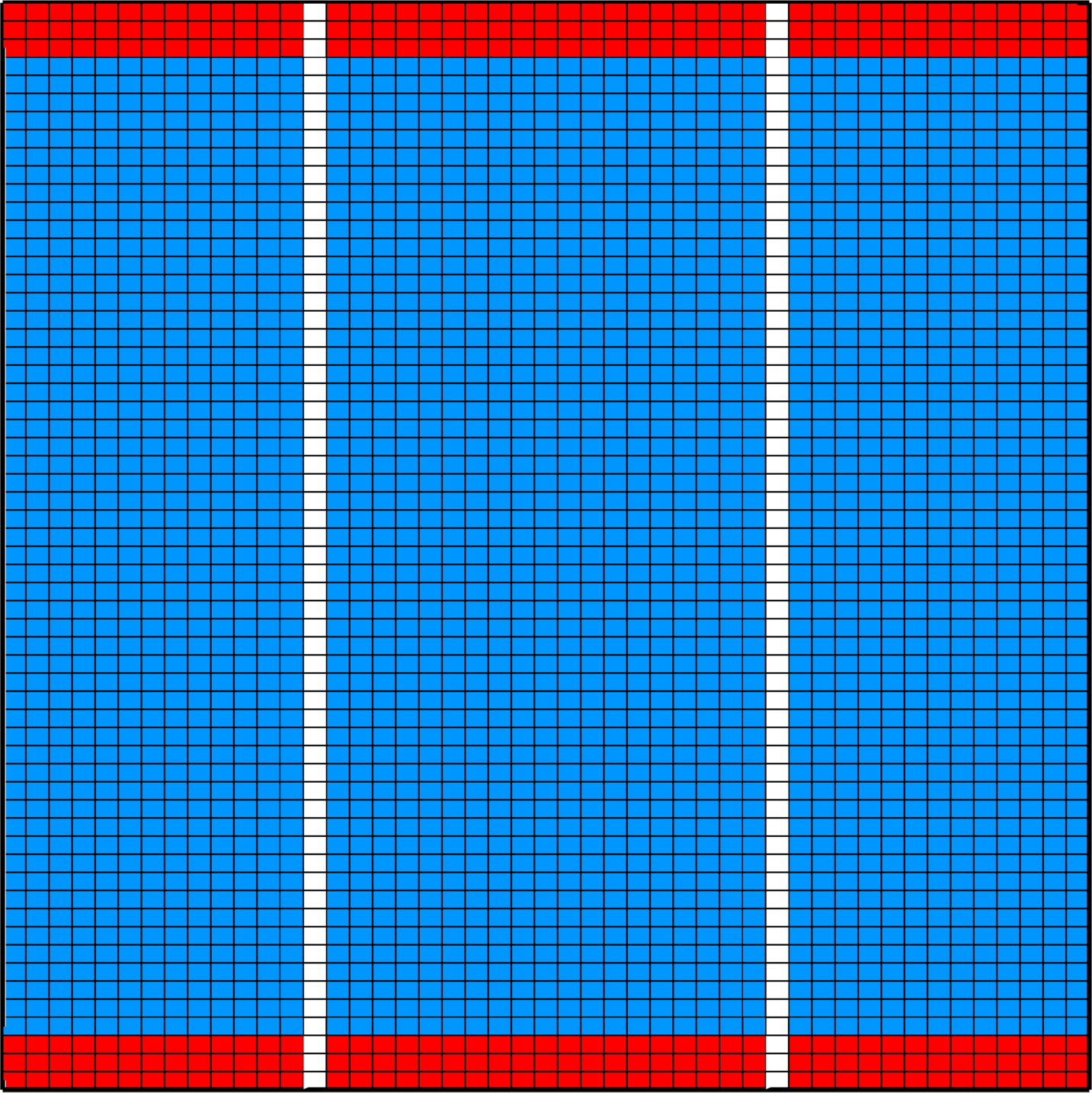}}
	\hspace{2pt}
	\subfigure[Data $\Eset_{10}$]{\includegraphics[width=0.14\textwidth]{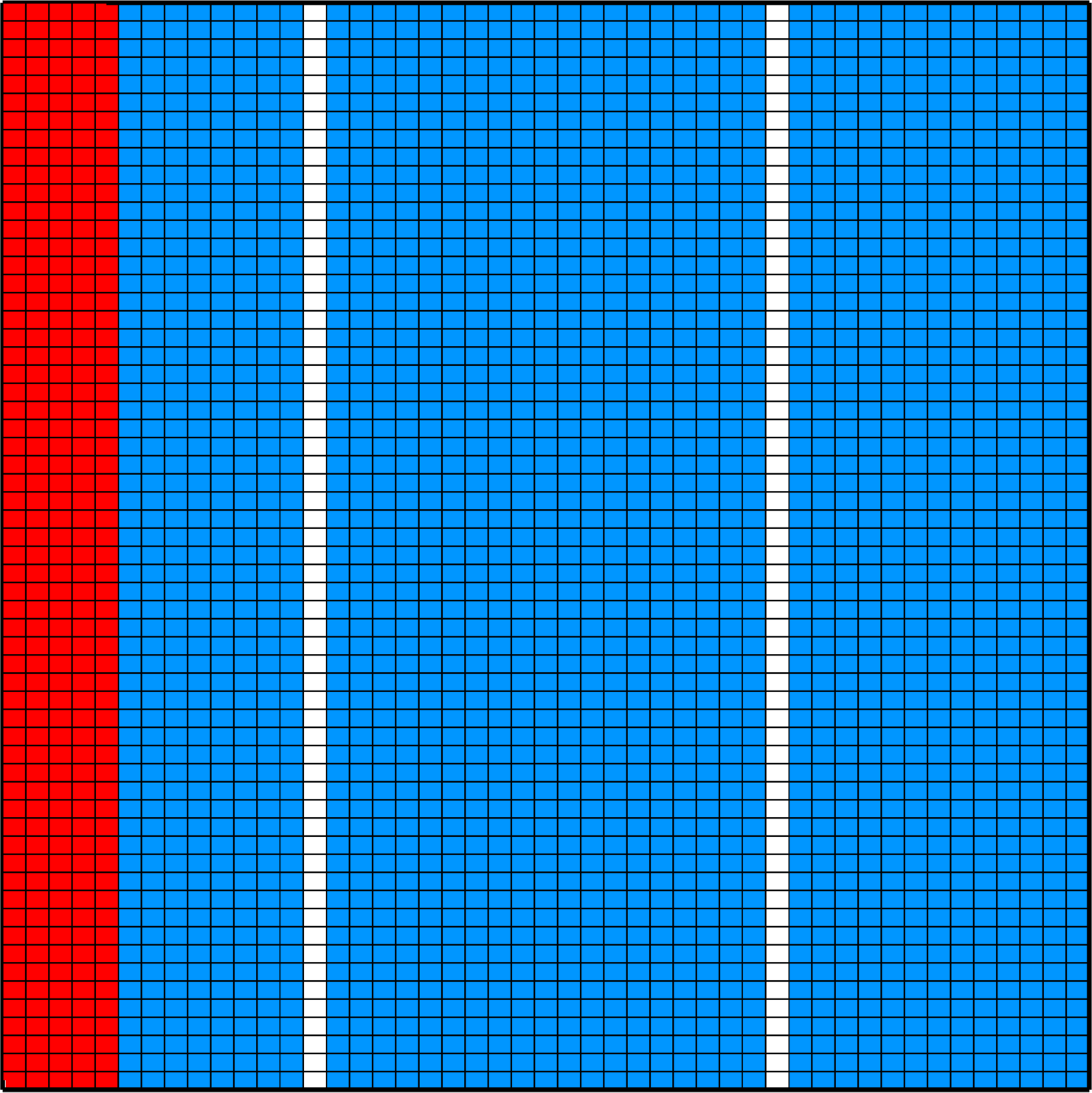}}
	\hspace{2pt}
	\subfigure[Data $\Eset_{11}$]{\includegraphics[width=0.14\textwidth]{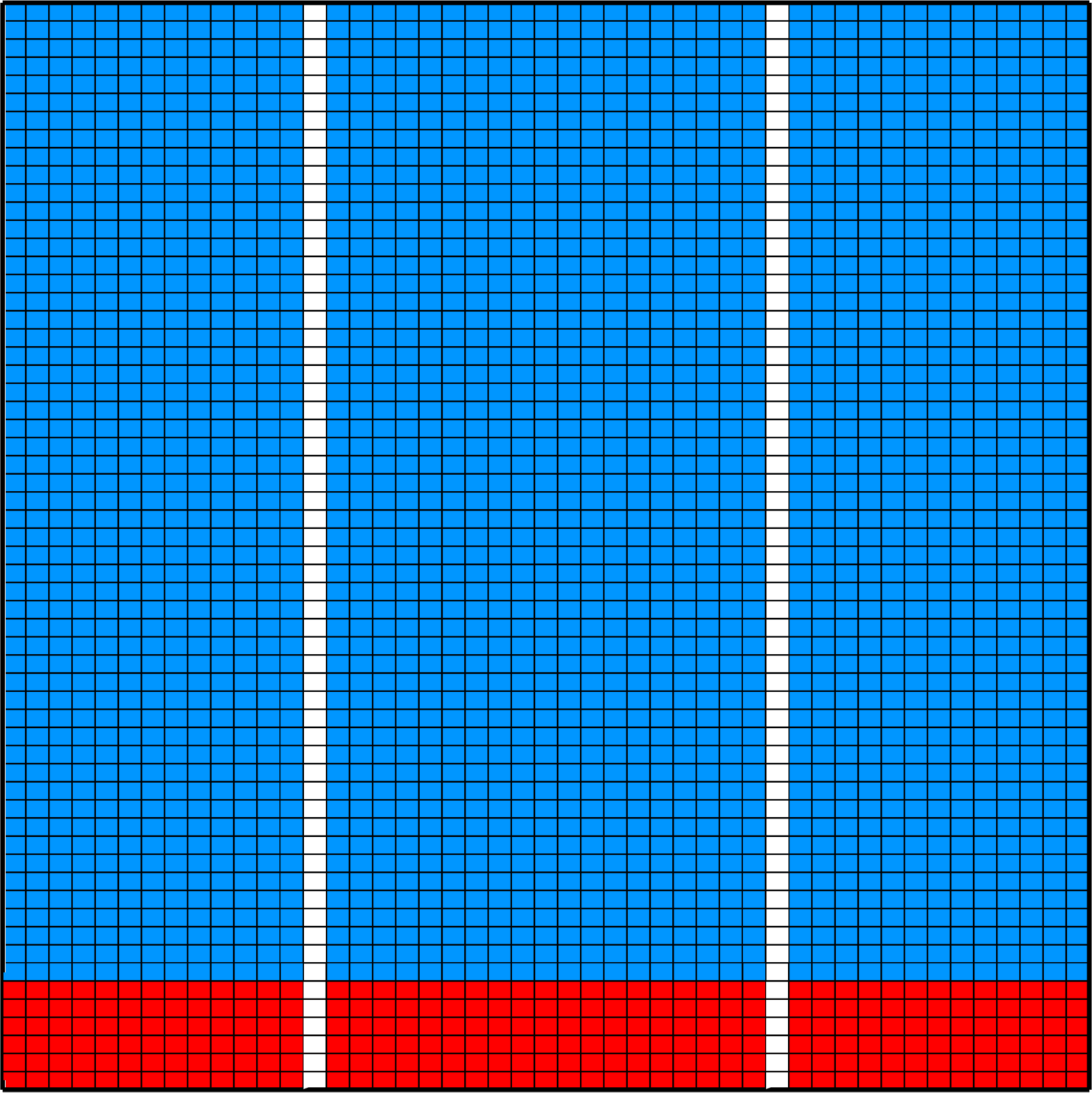}}
	\hspace{2pt}
	\subfigure[Data $\Eset_{12}$]{\includegraphics[width=0.14\textwidth]{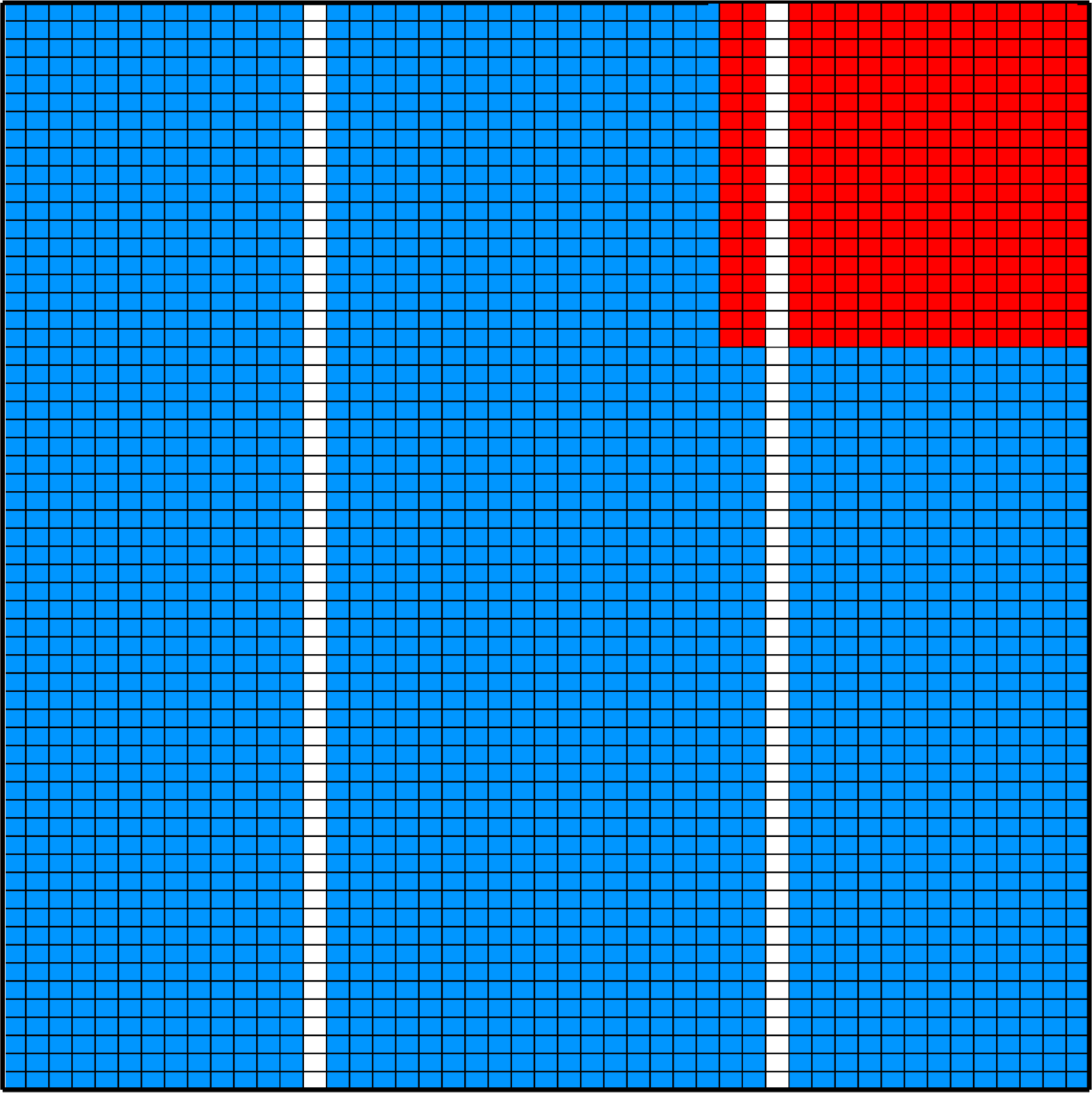}}
	
	\caption{Datasets for extrapolation. Twelve splittings of the dataset of $2,700$ optimal topologies into training and validation cases (blue) and test cases (red).}
	\label{fig:ExtrapolData}
\end{figure}

This section addresses the more complex framework, often overlooked in the scientific machine learning literature, of employing a surrogate model (here, $\texttt{FF}_{\!\eta}\texttt{-D}$) to perform extrapolation outside the training and validation datasets.
The datasets $\Eset_k, \ k=1,\ldots,12$ in Figure~\ref{fig:ExtrapolData} are defined for this task, using the cases displayed in blue for training and validation and the red ones for testing.
The datasets are obtained by completely removing from the training procedure all angles for a set of force locations (datasets $\Eset_1$,  $\Eset_4$, $\Eset_8$, and $\Eset_{10}$),  all locations for given values of the angle (datasets $\Eset_2$,  $\Eset_5$, $\Eset_9$, and $\Eset_{11}$), and entire regions of the parametric space $\Iset$ for datasets $\Eset_3$,  $\Eset_6$, $\Eset_7$, and $\Eset_{12}$.

\begin{figure}[!htb]
	\centering
%

	\subfigure[rMSE of $\texttt{FF}_{\!\eta}\texttt{-D}$ ]{\includegraphics[height=0.45\textwidth]{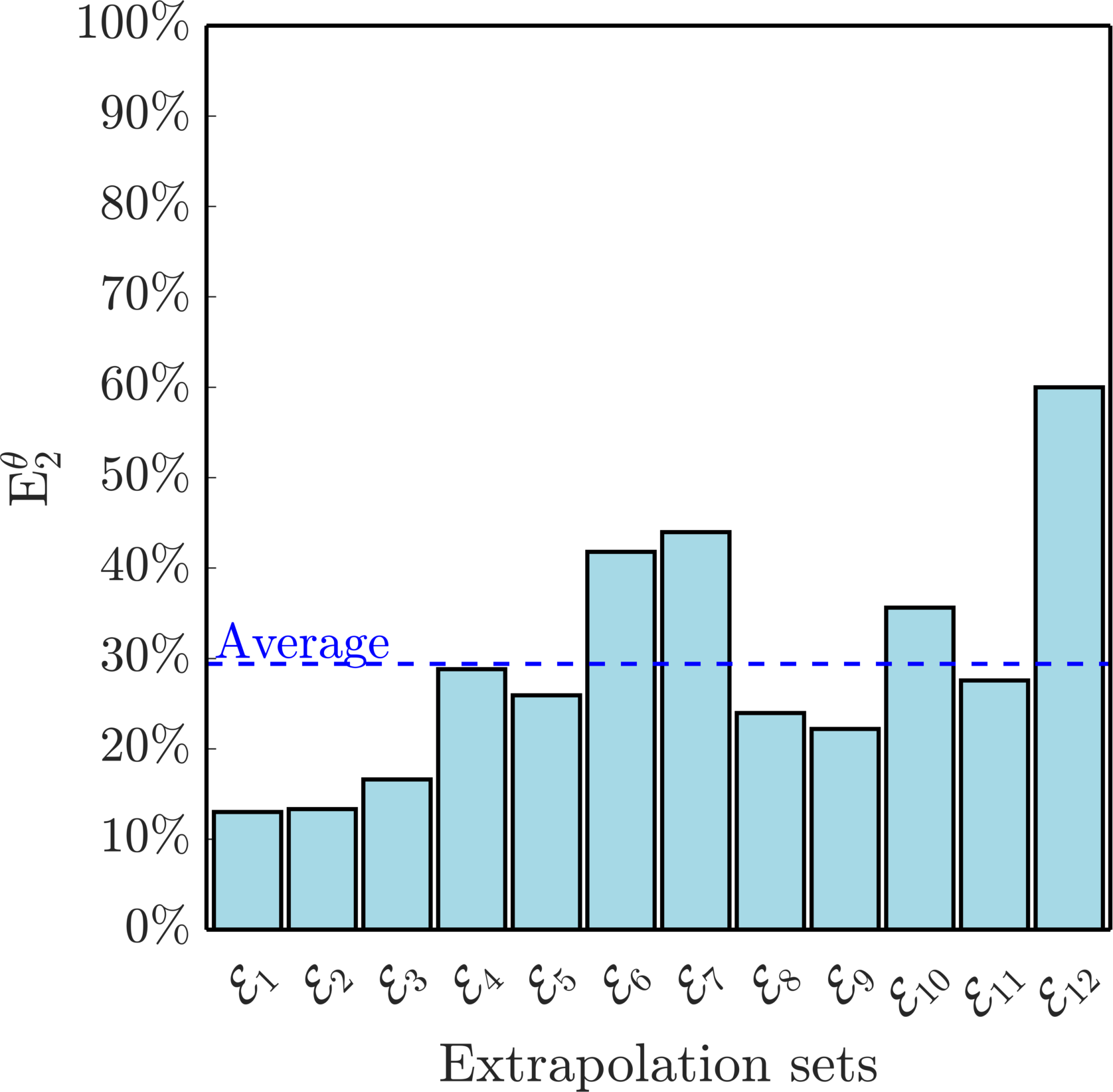}}
	\hspace{5pt}
	\subfigure[rMSE of $\texttt{AE}$]{\includegraphics[height=0.45\textwidth]{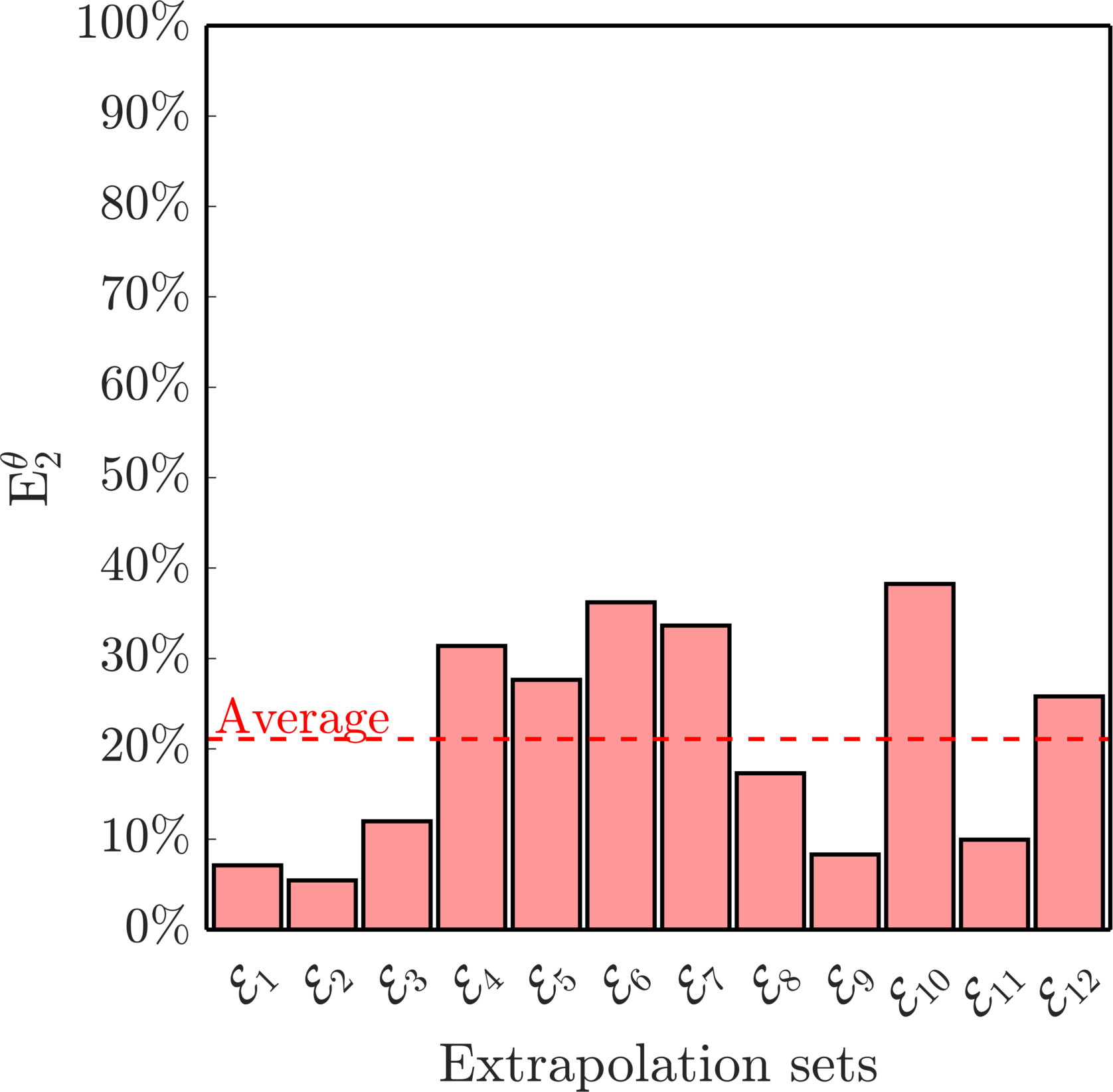}}
	
	\subfigure[rMAE of $\texttt{FF}_{\!\eta}\texttt{-D}$ ]{\includegraphics[height=0.45\textwidth]{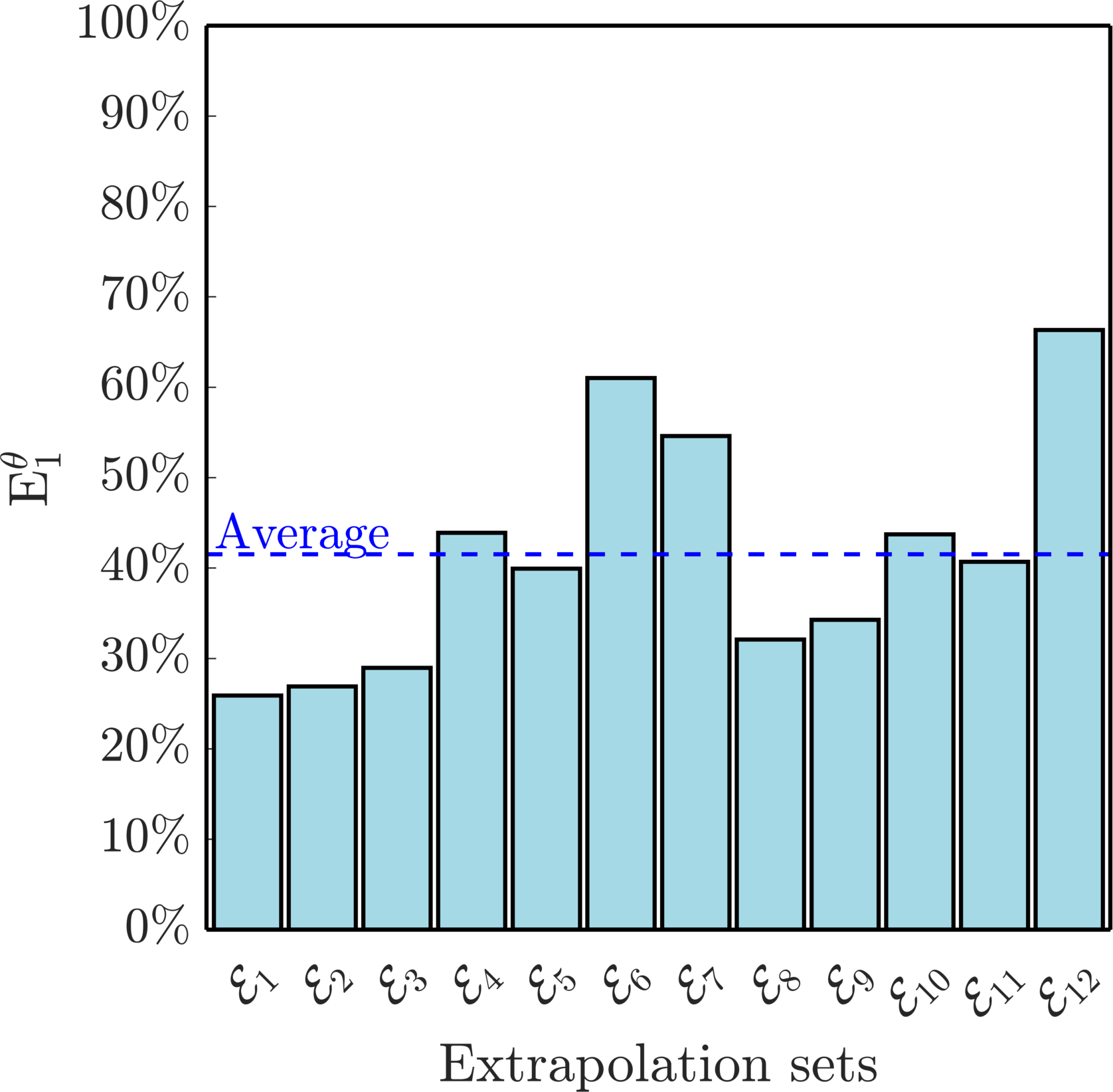}}
	\hspace{5pt}
	\subfigure[rMAE of $\texttt{AE}$]{\includegraphics[height=0.45\textwidth]{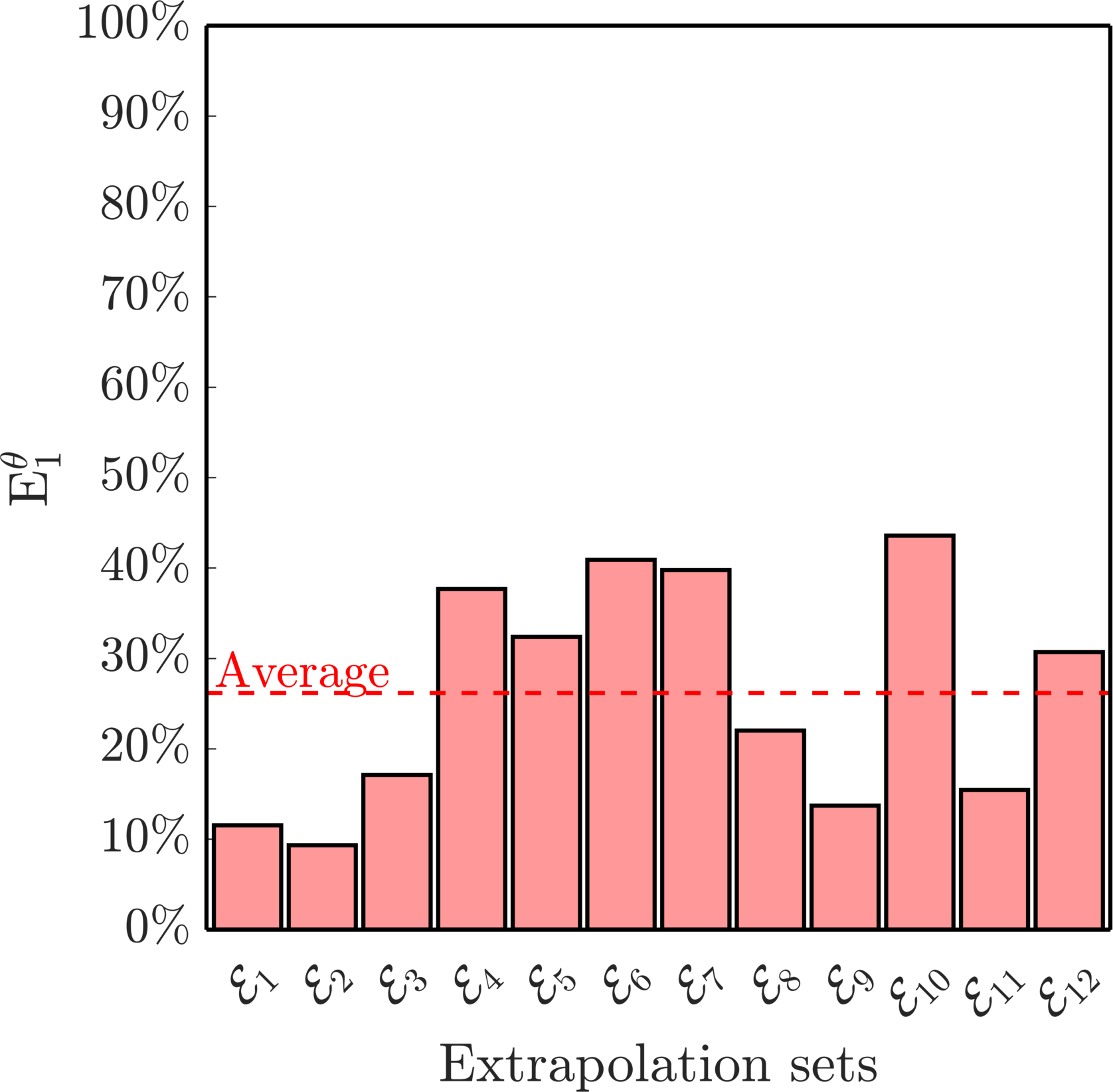}}
	
	\caption{Extrapolation error (top: rMSE; bottom: rMAE) computed over the test cases in the datasets $\Eset_k, \ k=1,\ldots,12$. }
	\label{fig:ExtrapolErr}
\end{figure}

Figure~\ref{fig:ExtrapolErr} displays the rMSE and rMAE for the $\texttt{FF}_{\!\eta}\texttt{-D}$ surrogate model when trained and tested on the extrapolation datasets $\Eset_k, \ k=1,\ldots,12$. 
Configurations such as $\Eset_1$ and $\Eset_2$ in which testing cases, although clustered,  are surrounded by multiple entries of the training and validation sets showcase behaviours similar to the ones reported in Section~\ref{sc:Validation} with random splitting of training, validations, and test sets. The corresponding values of rMSE and rMAE (respectively, $13\%$ and $25-26\%$) are indeed comparable to the accuracy achieved during cross-validation, with rMSE being approximately $12\%$ and rMAE $23\%$.

On the contrary,  the less information is provided to the surrogate models on a certain region of the parametric domain during training (e.g., datasets $\Eset_6$, $\Eset_7$, and $\Eset_{12}$), the higher the error in the associated predictions. In these three cases, the rMSE of $\texttt{FF}_{\!\eta}\texttt{-D}$ achieves $41.80\%$, $43.97\%$, and $60.00\%$.
It is worth noticing that these problems are particularly challenging as no prior knowledge is available for \emph{any} angles within a certain interval of values and \emph{any} position of the force. As a matter of fact, the target optimal architecture represented by the autoencoder $\texttt{AE}$ also attains low precisions in these scenarios, with errors spanning from $25.80\%$ to $36.20\%$, while only performing a compression and a reconstruction of the \emph{ground truth} input data $\bthetaIn$ without learning any functional dependence of the solution on the parameters $\bEta$.

A detailed overview of the values of rMSE and rMAE for the testing of the $\texttt{FF}_{\!\eta}\texttt{-D}$ surrogate model is reported in Table~\ref{tab:errExtrapol}, together with the data associated with the baseline autoencoder $\texttt{AE}$.
\begin{table}[!htb]
\centering
\begin{tabular}{| c || c | c || c | c || c | c |}
\hline
\multirow{2}{*}{Dataset}& \multicolumn{2}{c||}{rMSE $\Err{2}$} & \multicolumn{2}{c||}{rMAE $\Err{1}$}  & \multicolumn{2}{c|}{$\nA/\nA^{\text{max}}$}\\
\cline{2-7} 
& $\texttt{FF}_{\!\eta}\texttt{-D}$ & $\texttt{AE}$ & $\texttt{FF}_{\!\eta}\texttt{-D}$ & $\texttt{AE}$ & $\texttt{FF}_{\!\eta}\texttt{-D}$ & $\texttt{AE}$ \\
\hline
$\Eset_1$ & $13.00\%$ & $7.12\%$ & $25.91\%$ & $11.55\%$ & $19/25$ & $14/25$ \\
\hline
$\Eset_2$ & $13.34\%$ & $5.45\%$ & $26.91\%$ & $9.34\%$ & $22/25$ & $18/25$ \\
\hline
$\Eset_3$ & $16.61\%$ & $11.99\%$ & $28.96\%$ & $17.12\%$ & $16/25$ & $13/25$ \\
\hline
$\Eset_4$ & $28.82\%$ & $31.38\%$ & $43.89\%$ & $37.69\%$ & $18/25$ & $18/25$ \\
\hline
$\Eset_5$ & $25.93\%$ & $27.65\%$ & $39.94\%$ & $32.40\%$ & $16/25$ & $15/25$ \\
\hline
$\Eset_6$ & $41.80\%$ & $36.20\%$ & $61.01\%$ & $40.92\%$ & $17/25$ & $21/25$ \\
\hline
$\Eset_7$ & $43.97\%$ & $33.64\%$ & $54.60\%$ & $39.79\%$ & $20/25$ & $12/25$ \\
\hline
$\Eset_8$ & $23.96\%$ & $17.30\%$ & $32.09\%$ & $22.02\%$ & $20/25$ & $14/25$ \\
\hline
$\Eset_9$ & $22.19\%$ & $8.32\%$ & $34.27\%$ & $13.73\%$ & $13/25$ & $9/25$\\
\hline
$\Eset_{10}$ & $35.62\%$ & $38.24\%$ & $43.74\%$ & $43.57\%$ & $17/25$ & $10/25$ \\
\hline
$\Eset_{11}$ & $27.58\%$ & $9.96\%$ & $40.69\%$ & $15.48\%$ & $14/25$ & $14/25$ \\
\hline
$\Eset_{12}$ & $60.00\%$ & $25.80\%$ & $66.33\%$ & $30.70\%$ & $18/25$ & $10/25$ \\
\hline\hline
Average & $29.40\%$ & $21.09\%$ & $41.53\%$ & $26.19\%$ & $17.5/25$ & $14/25$ \\
\hline
\end{tabular}

\caption{Relative mean squared error,  relative mean absolute error, and number of active parameters in the latent space for the twelve extrapolation datasets.}
\label{tab:errExtrapol}
\end{table}
In this context, the case of datasets $\Eset_4$, $\Eset_5$, and $\Eset_{10}$ is of particular interest since $\texttt{FF}_{\!\eta}\texttt{-D}$ slightly outperforms the autoencoder $\texttt{AE}$.
In datasets $\Eset_4$, and $\Eset_{10}$ the set of unseen cases corresponds to all the angles of the force applied on the central portion of the lateral surface $\Ga{r}$ and on the left-most portion of the bottom surface $\Ga{b}$. Conversely,  dataset $\Eset_5$ requires the extrapolation for a set of angles applied on all the admissible positions of the force. The complete lack of information on angles ($\Eset_4$ and $\Eset_{10}$) and on positions ($\Eset_5$) is responsible for relative mean squared errors of the autoencoder $\texttt{AE}$ of $31.38\%$, $27.65\%$, and $38.24\%$.
In these cases, $\texttt{FF}_{\!\eta}\texttt{-D}$ attains errors of $28.82\%$, $25.93\%$, and $35.62\%$, respectively. Whilst the difference is not particularly large, the trend is significant: in the complete absence of information on all angles or all positions during training and validation, $\texttt{FF}_{\!\eta}\texttt{-D}$ is capable of learning relevant information by exploiting the functional dependence of the \emph{ground truth} optimal topology $\bthetaIn$ on the input parameters $\bEta$.  On the contrary, the autoencoder $\texttt{AE}$ being purely data-driven disregards any relation between $\bEta$ and $\bthetaIn$.

Finally, even in the presence of training and validation sets with entire portions of missing data,  $\texttt{FF}_{\!\eta}\texttt{-D}$ is capable of achieving significant compression, with an average of $\nA=17.5$ active entries in the latent space over the user-defined maximum $\nA^{\text{max}}=25$. The corresponding compression factor $\nT/\nA$ is approximately $731$, whereas the autoencoder attains $914$, not learning any functional dependence of the solution on the parameters $\bEta$.

\subsection{Surrogate-based optimisation extrapolating \emph{quasi-optimal} topologies}
\label{sc:ExtrapolationOpti}
 
\begin{figure}[!htb]
	\centering	
	\subfigure[Case 6]{\includegraphics[width=0.3\textwidth]{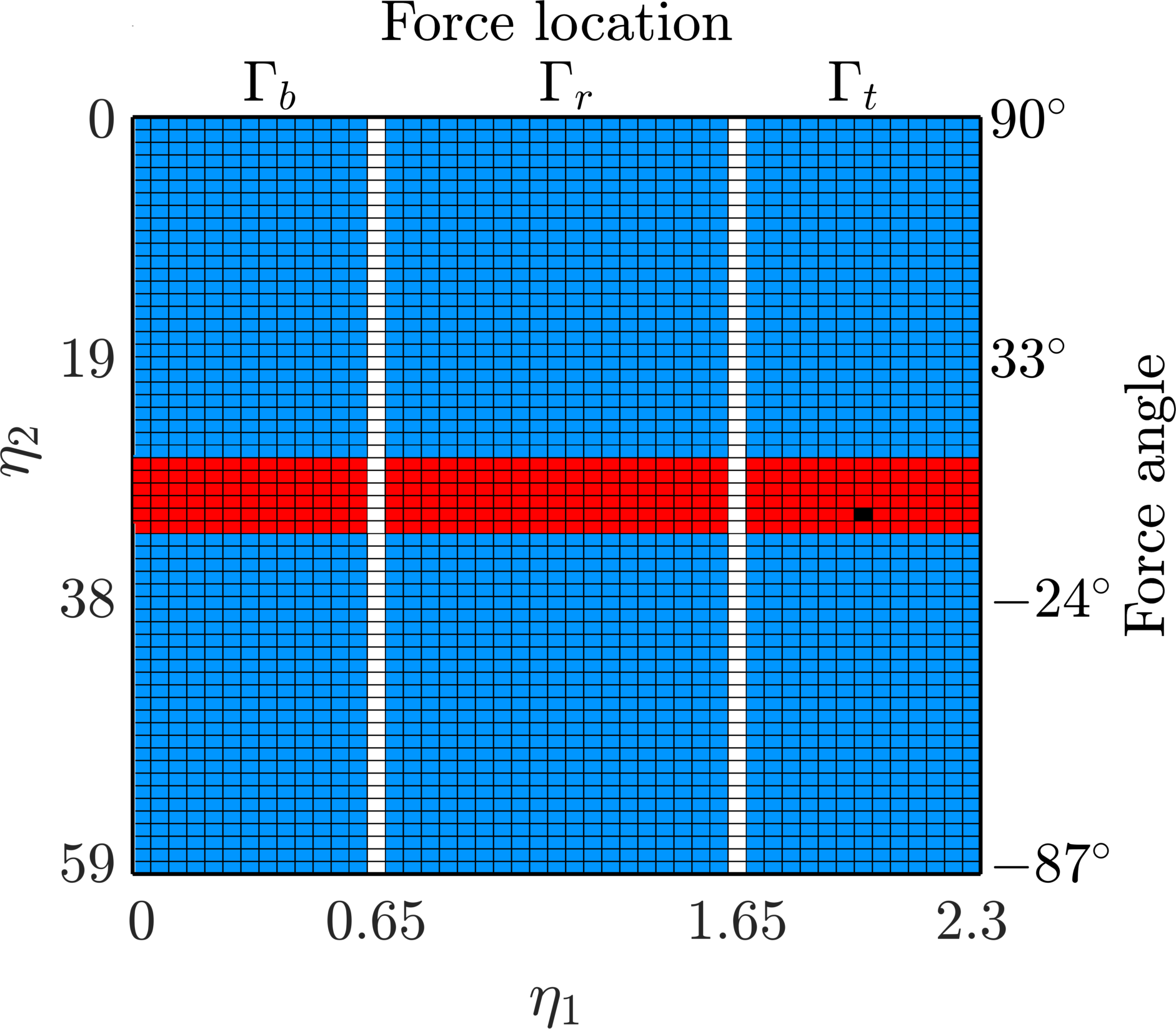}}
	\hspace{3pt}
	\subfigure[Case 7]{\includegraphics[width=0.3\textwidth]{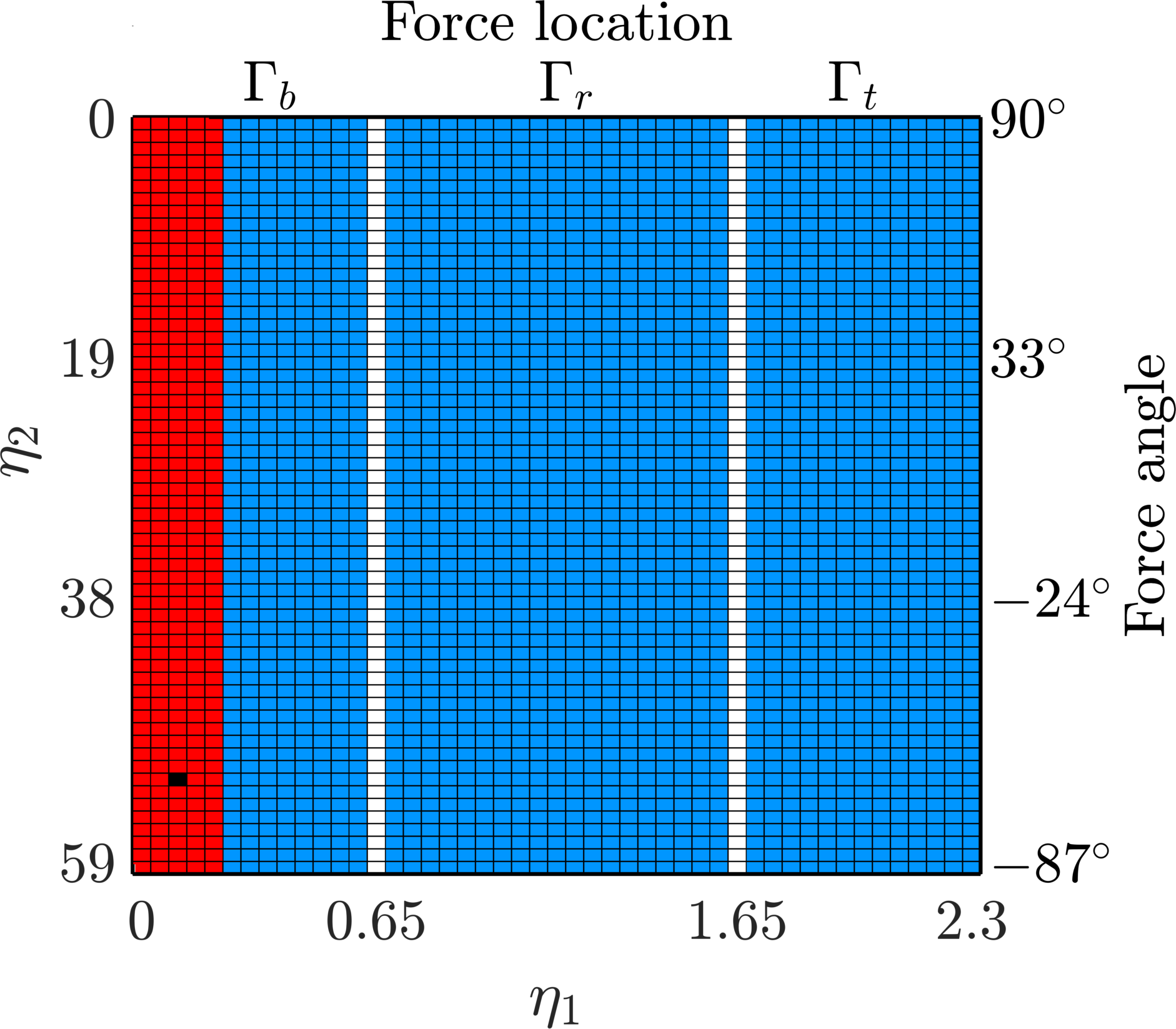}}
	\hspace{3pt}
	\subfigure[Case 8]{\includegraphics[width=0.3\textwidth]{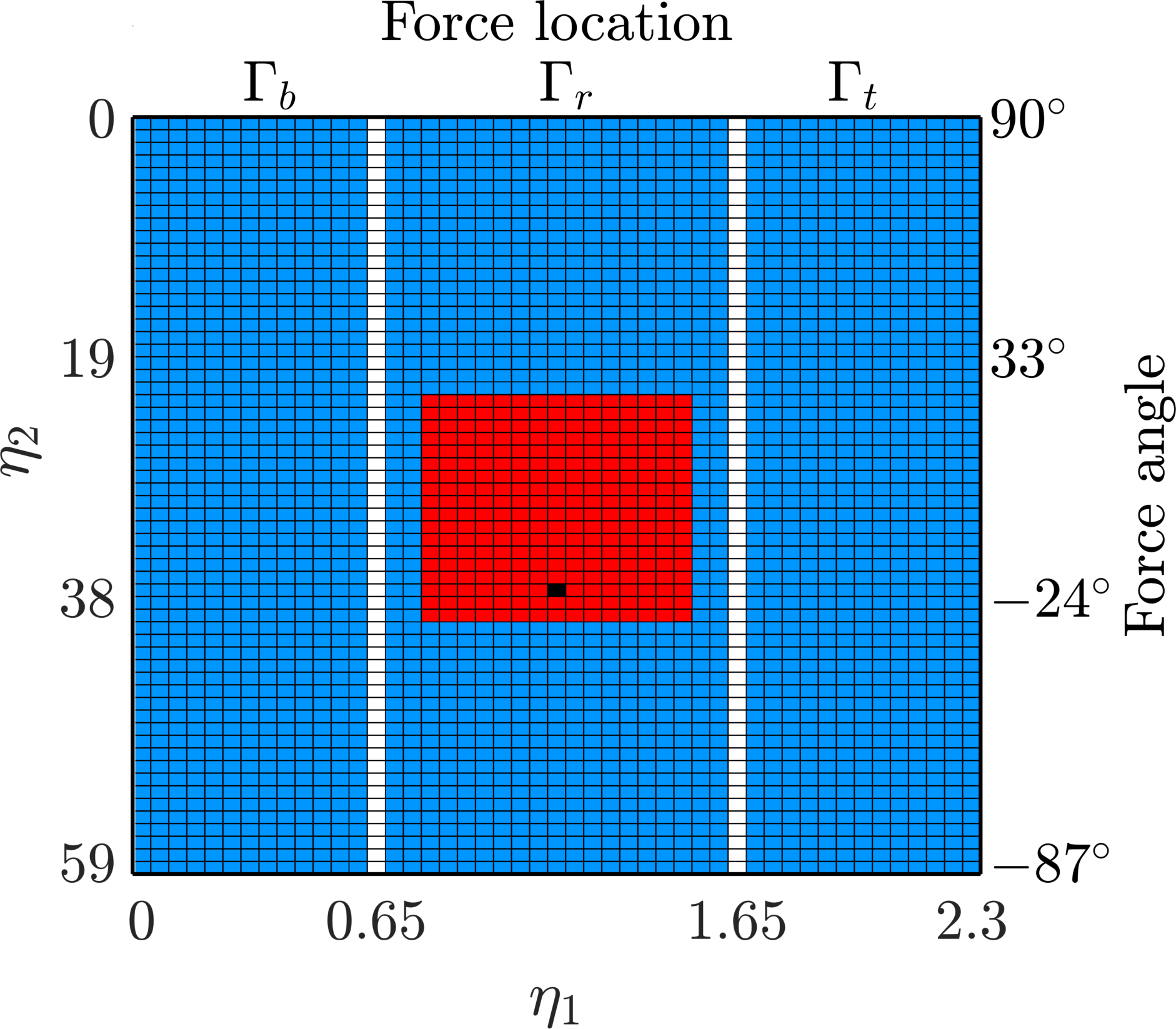}}

	\subfigure[Case 9]{\includegraphics[width=0.3\textwidth]{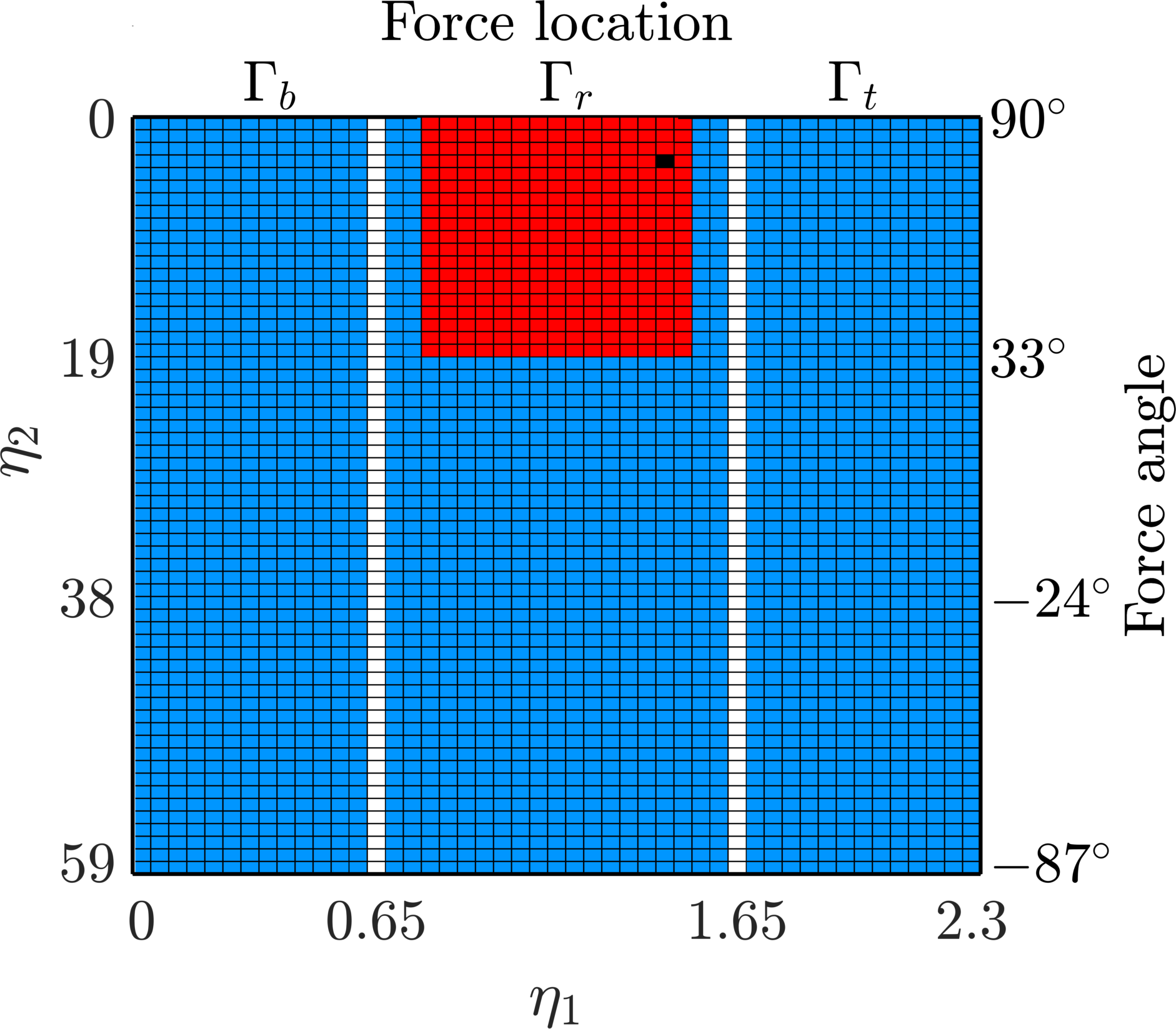}}
	\hspace{1pt}
	\subfigure[Case 10]{\includegraphics[width=0.3\textwidth]{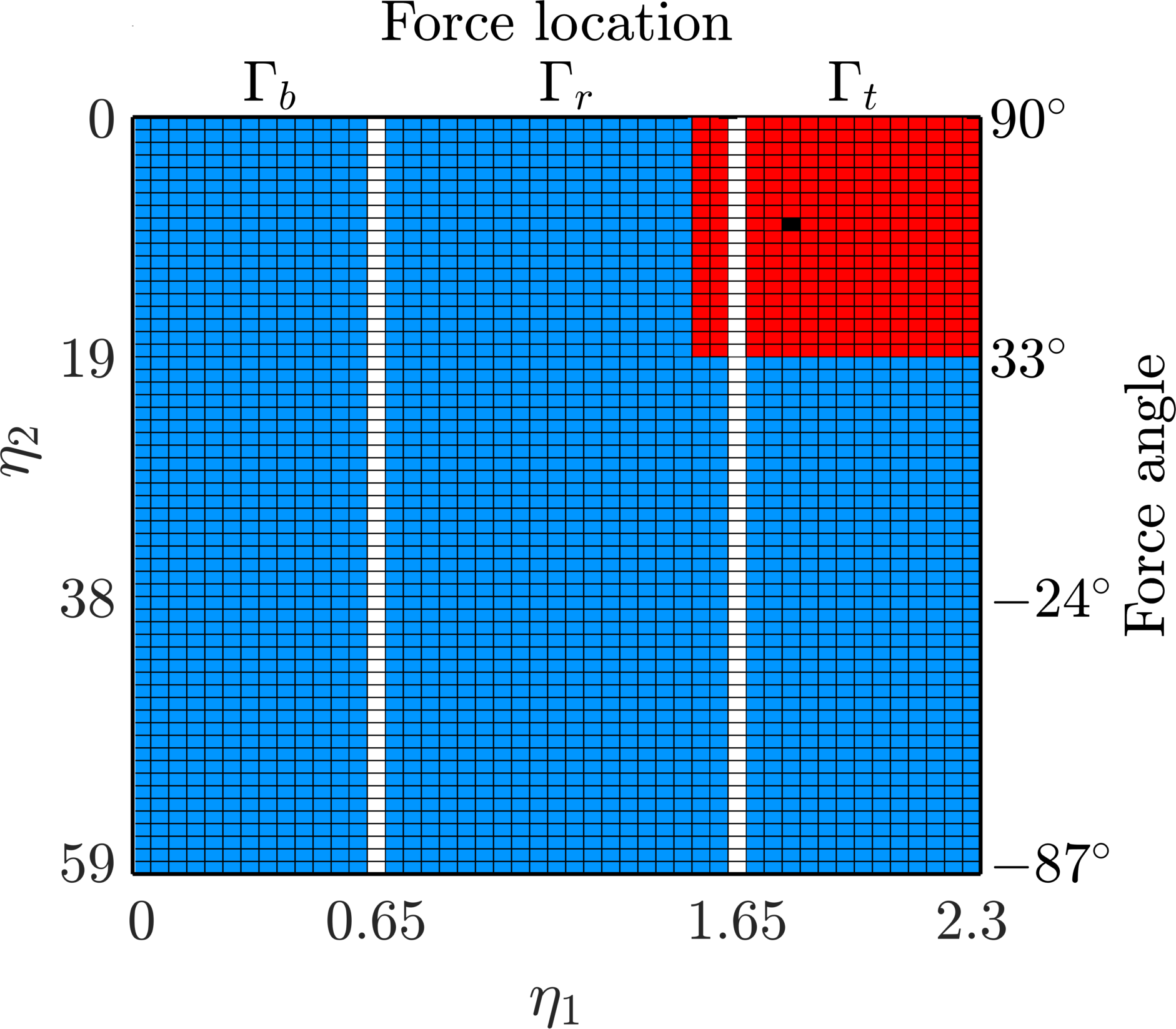}}
	
	\caption{Extrapolation test cases. Blue: training and validation sets. Red: test set. Black: unseen case used for extrapolation test.}
	\label{fig:ExtrapolTest}
\end{figure}

The last section evaluates the suitability of the \emph{quasi-optimal} topologies constructed in Section~\ref{sc:Generalisation} for extrapolated unseen cases as initial conditions for the surrogate-based optimisation algorithm~\ref{alg:topOptSurrogate}.
In particular,  five unseen test cases are selected from the most challenging datasets previously identified, namely $\Eset_5$, $\Eset_{10}$, $\Eset_6$, $\Eset_7$, and $\Eset_{12}$ (Figure~\ref{fig:ExtrapolTest}). The details of the extrapolation data are described in Table~\ref{tab:dataOnlineExtrapol}.
\begin{table}[!htb]
	\centering
	\resizebox{\linewidth}{!}{%
	\begin{tabular}{| c | c || c | c || c | c | c | c |}
	\hline
	\multirow{2}{*}{Case} & \multirow{2}{*}{Dataset}& \multicolumn{2}{c||}{Parameters} & \multicolumn{4}{c|}{Force} \\
	\cline{3-8}
	& & $\eta_1$ & $\eta_2$ & Boundary & Position ($x$) & Position ($y$) & Angle \\
	\hline
	 6 & $\Eset_5$ & $2.00$ & $32$ & $\Ga{t}$ & $[0.60 ,0.70]$ & $1$ & $-6^{\circ}$ \\
	\hline
	 7 & $\Eset_{10}$ & $0.10$ & $52$ & $\Ga{b}$ & $[0.40 ,0.50]$ & $0$ & $-66^{\circ}$ \\
	\hline
	 8 & $\Eset_6$ & $1.15$ & $37$ & $\Ga{r}$ & $1$ & $[0.45 ,0.55]$ & $-21^{\circ}$ \\
	\hline
	 9 & $\Eset_7$ & $1.45$ & $3$ & $\Ga{r}$ & $1$ & $[0.75 ,0.85]$ & $81^{\circ}$ \\
	\hline
	 10 & $\Eset_{12}$ & $1.80$ & $8$ & $\Ga{t}$ & $[0.80 ,0.90]$ & $1$ & $66^{\circ}$ \\
	\hline
	\end{tabular}
	}

	\caption{Extrapolation configurations for online evaluations of the surrogate-based optimisation strategy.}
	\label{tab:dataOnlineExtrapol}
\end{table}

\begin{figure}[!htb]
	\centering
	\subfigure[Case 6: $\bthetaRef_{\text{opt}}$]{\includegraphics[width=0.3\textwidth]{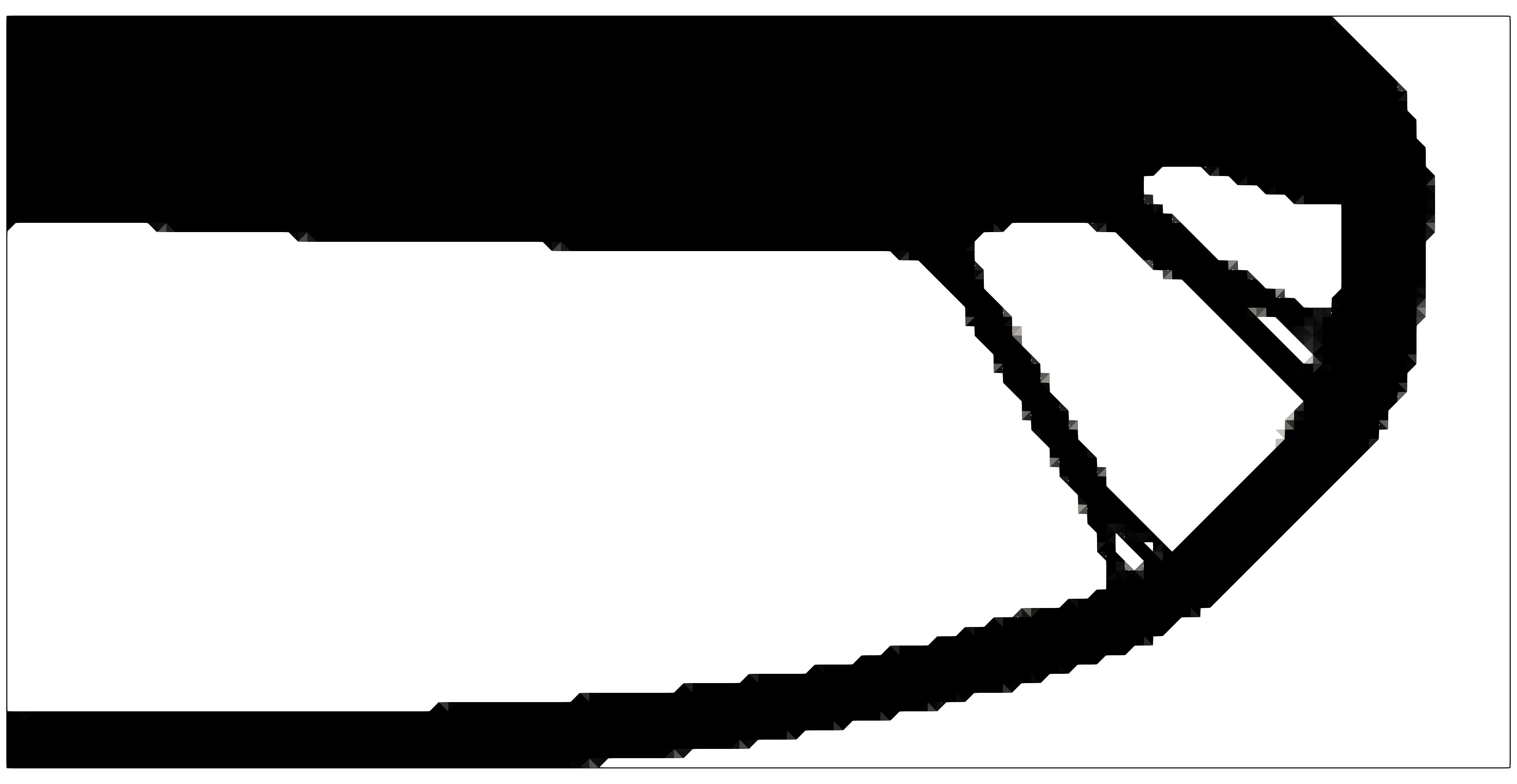}}
	\hspace{5pt}
	\subfigure[Case 6: $\bthetaEta$]{\includegraphics[width=0.3\textwidth]{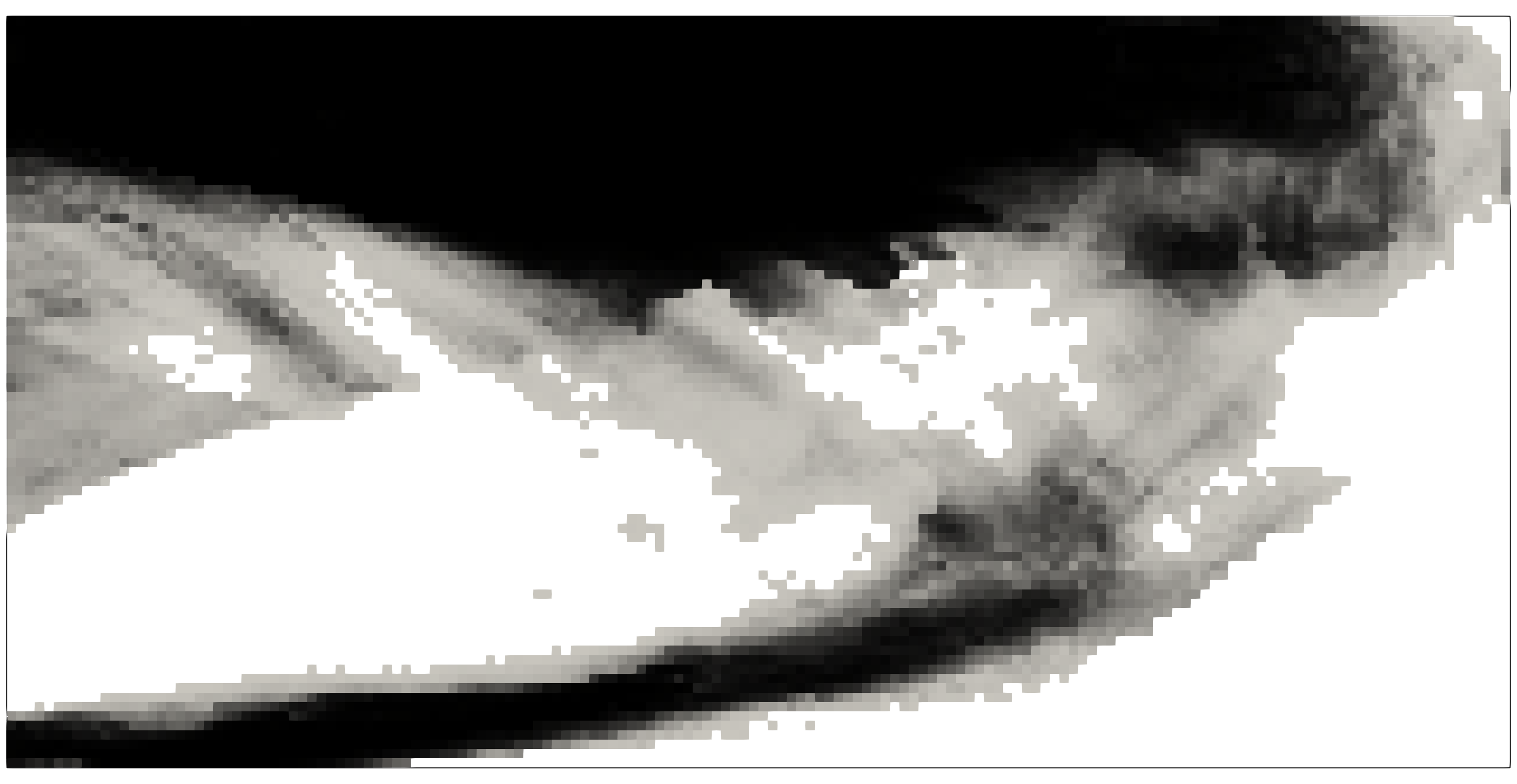}}
	\hspace{5pt}
	\subfigure[Case 6: $\bthetaEta_{\text{opt}}$]{\includegraphics[width=0.3\textwidth]{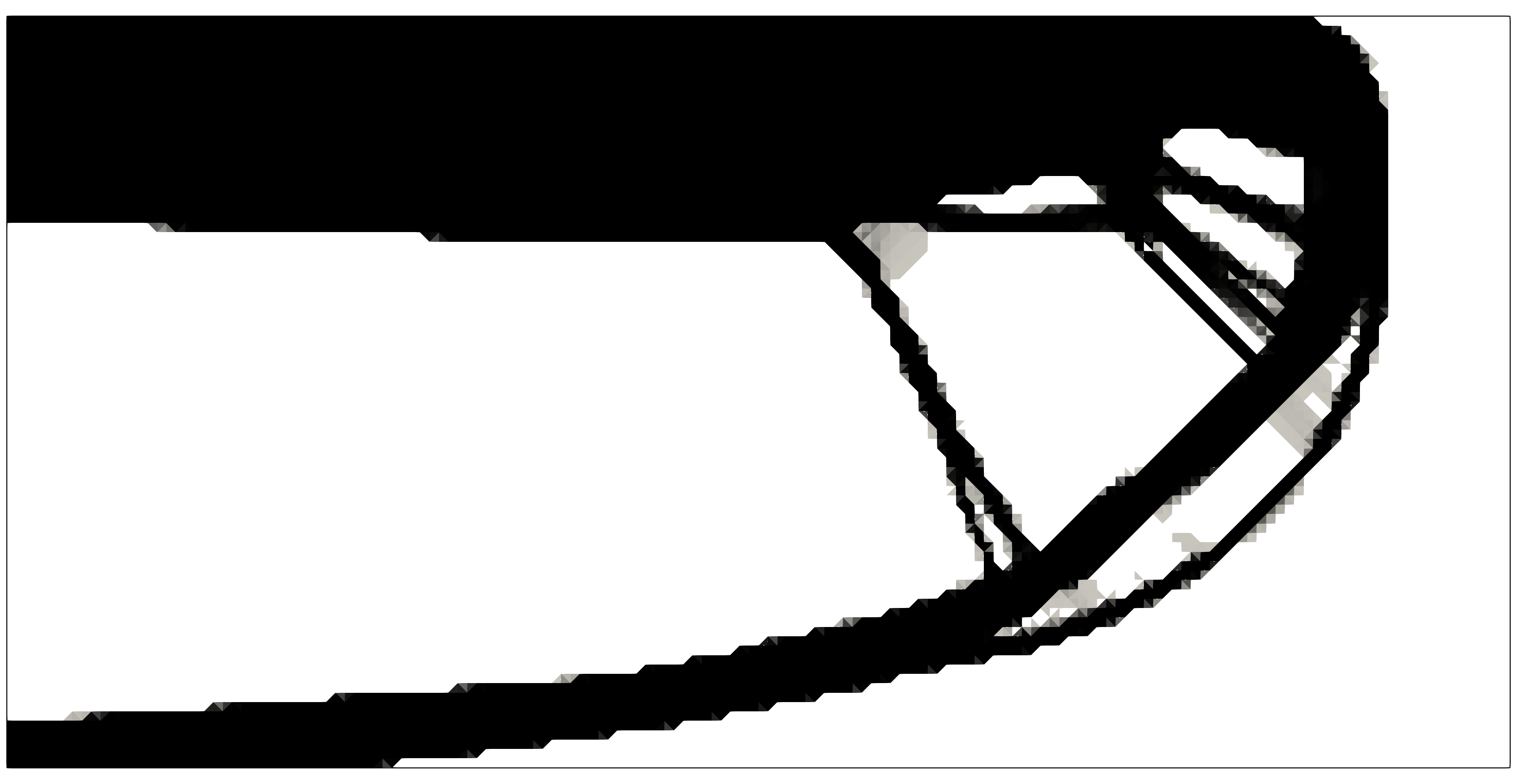}}
	
	\subfigure[Case 7: $\bthetaRef_{\text{opt}}$]{\includegraphics[width=0.3\textwidth]{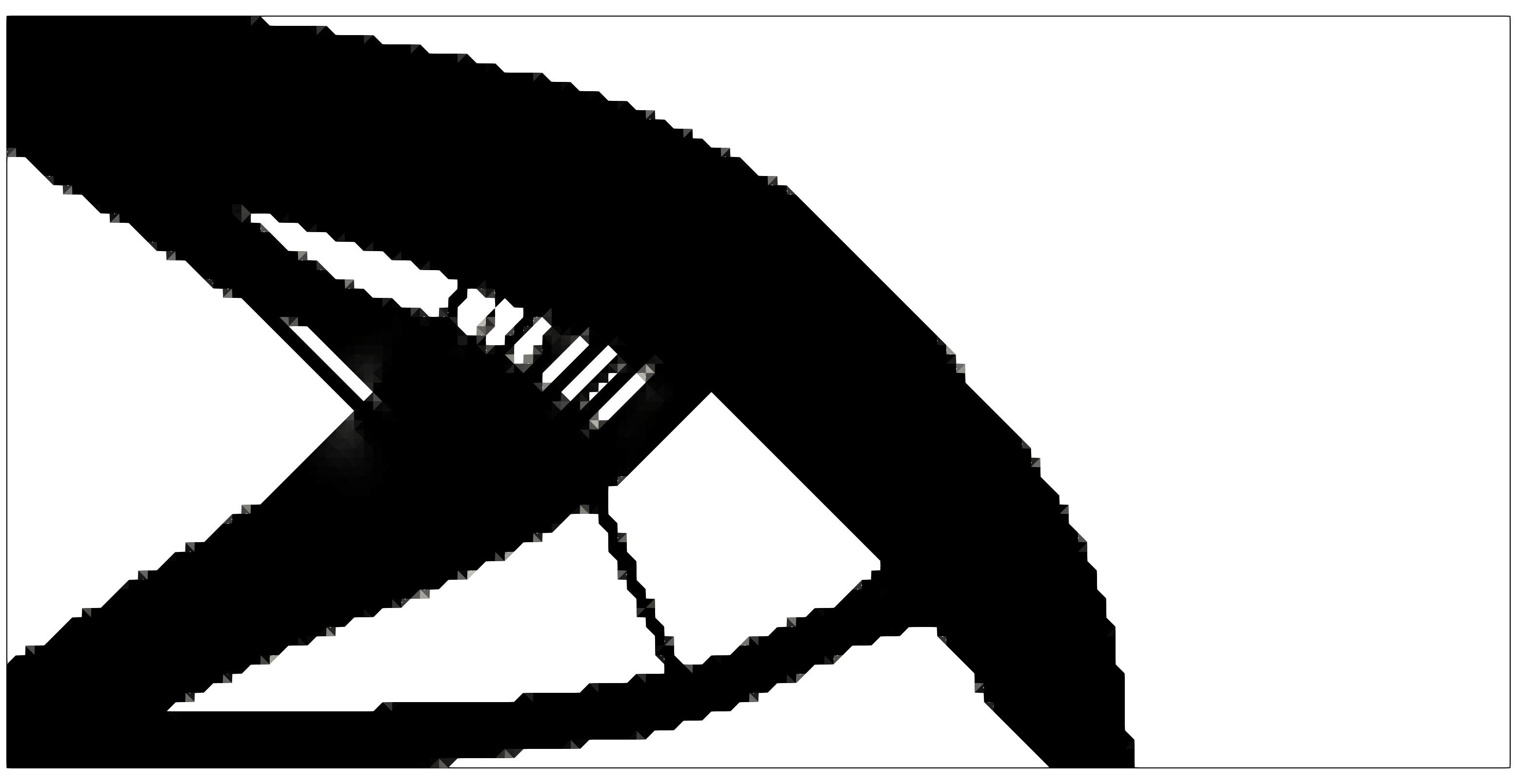}}
	\hspace{5pt}
	\subfigure[Case 7: $\bthetaEta$]{\includegraphics[width=0.3\textwidth]{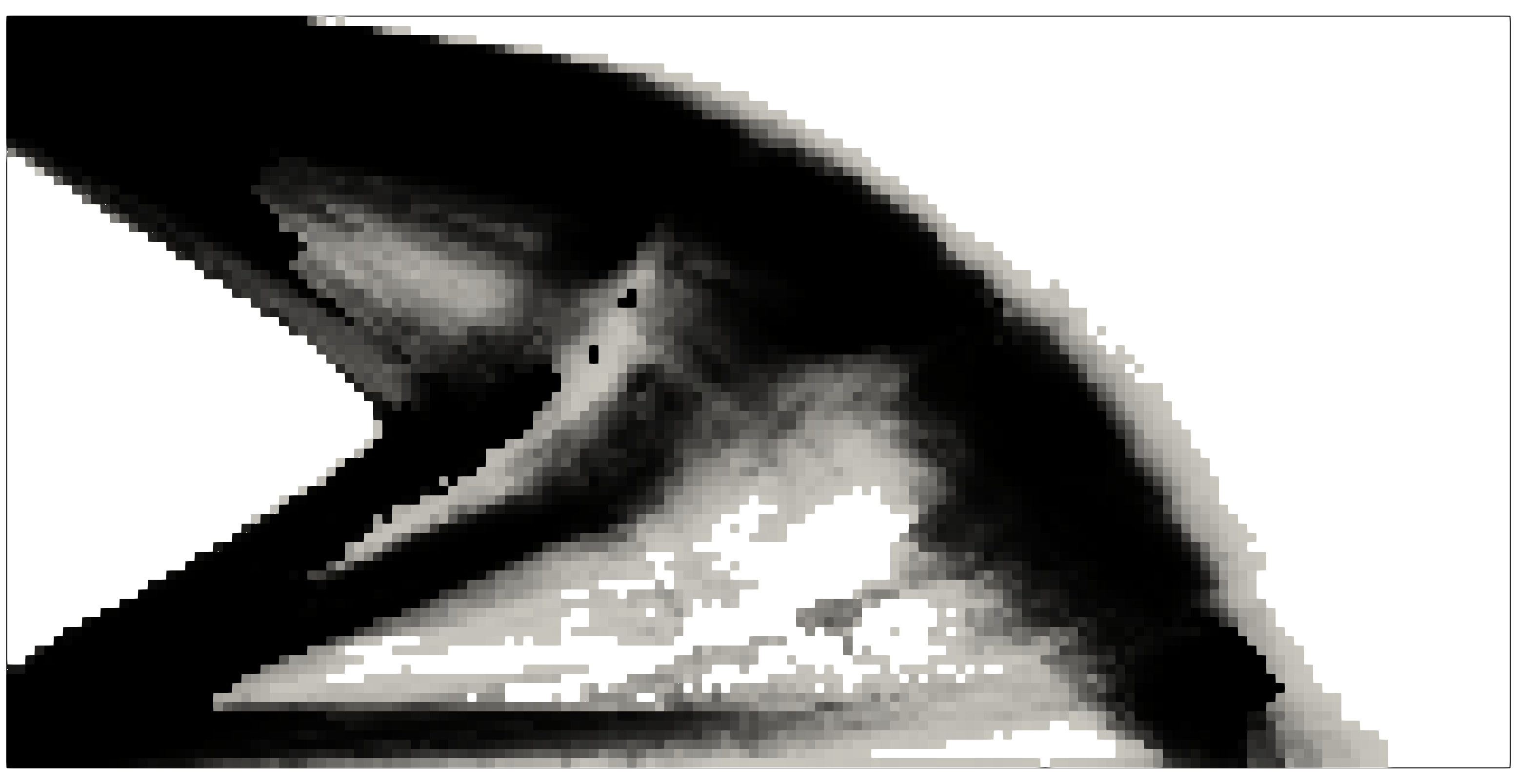}}
	\hspace{5pt}
	\subfigure[Case 7: $\bthetaEta_{\text{opt}}$]{\includegraphics[width=0.3\textwidth]{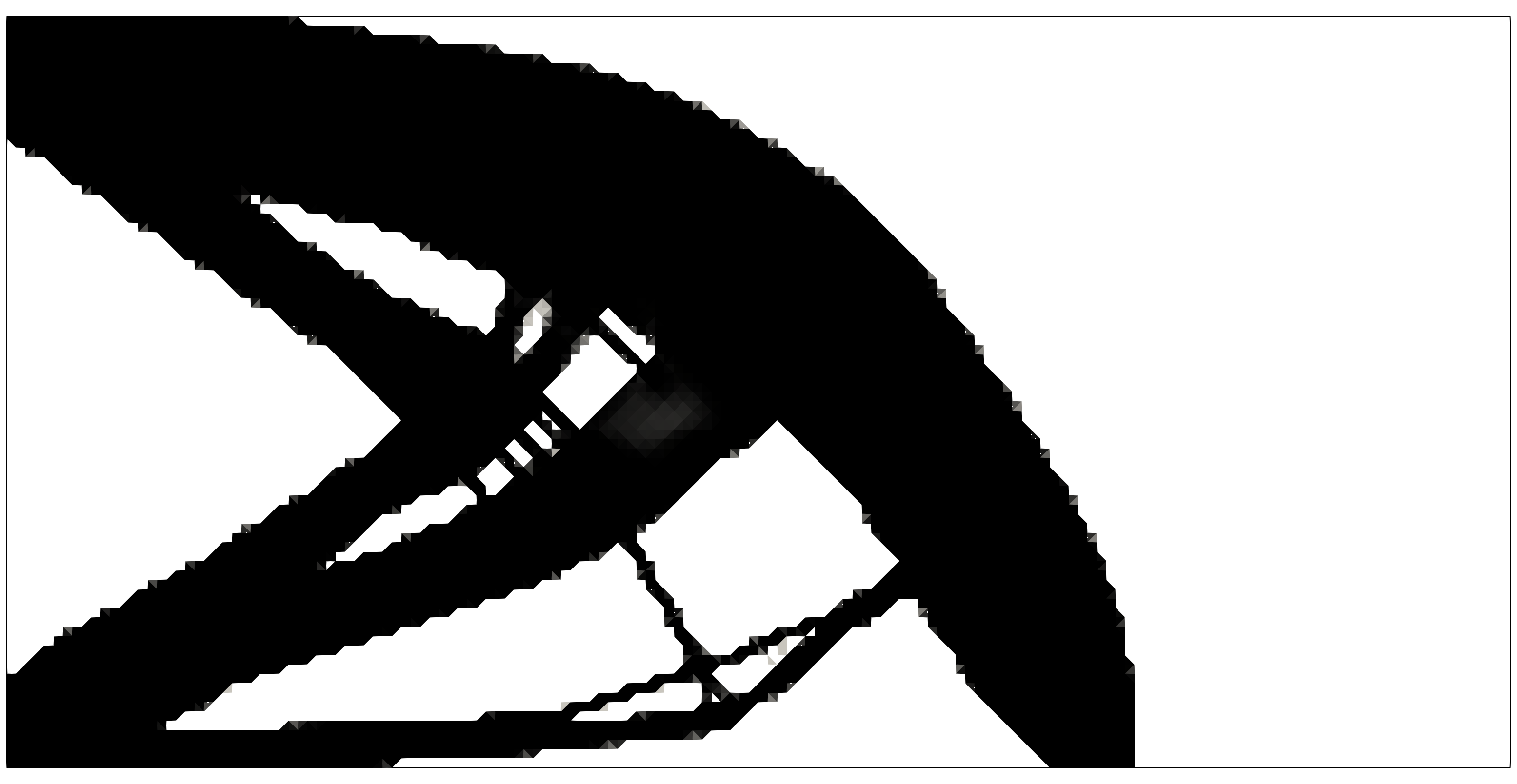}}
	
	\subfigure[Case 8: $\bthetaRef_{\text{opt}}$]{\includegraphics[width=0.3\textwidth]{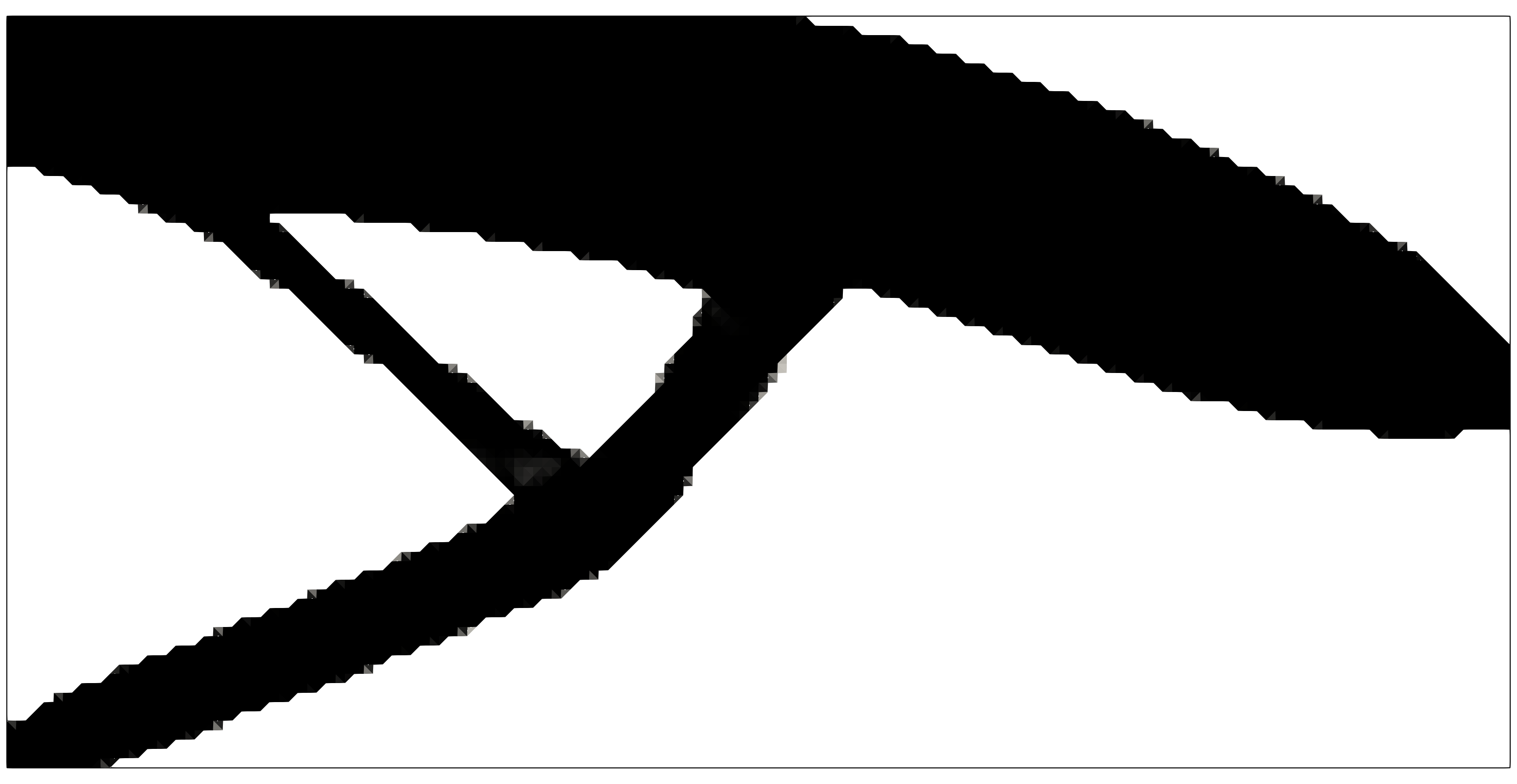}}
	\hspace{5pt}
	\subfigure[Case 8: $\bthetaEta$]{\includegraphics[width=0.3\textwidth]{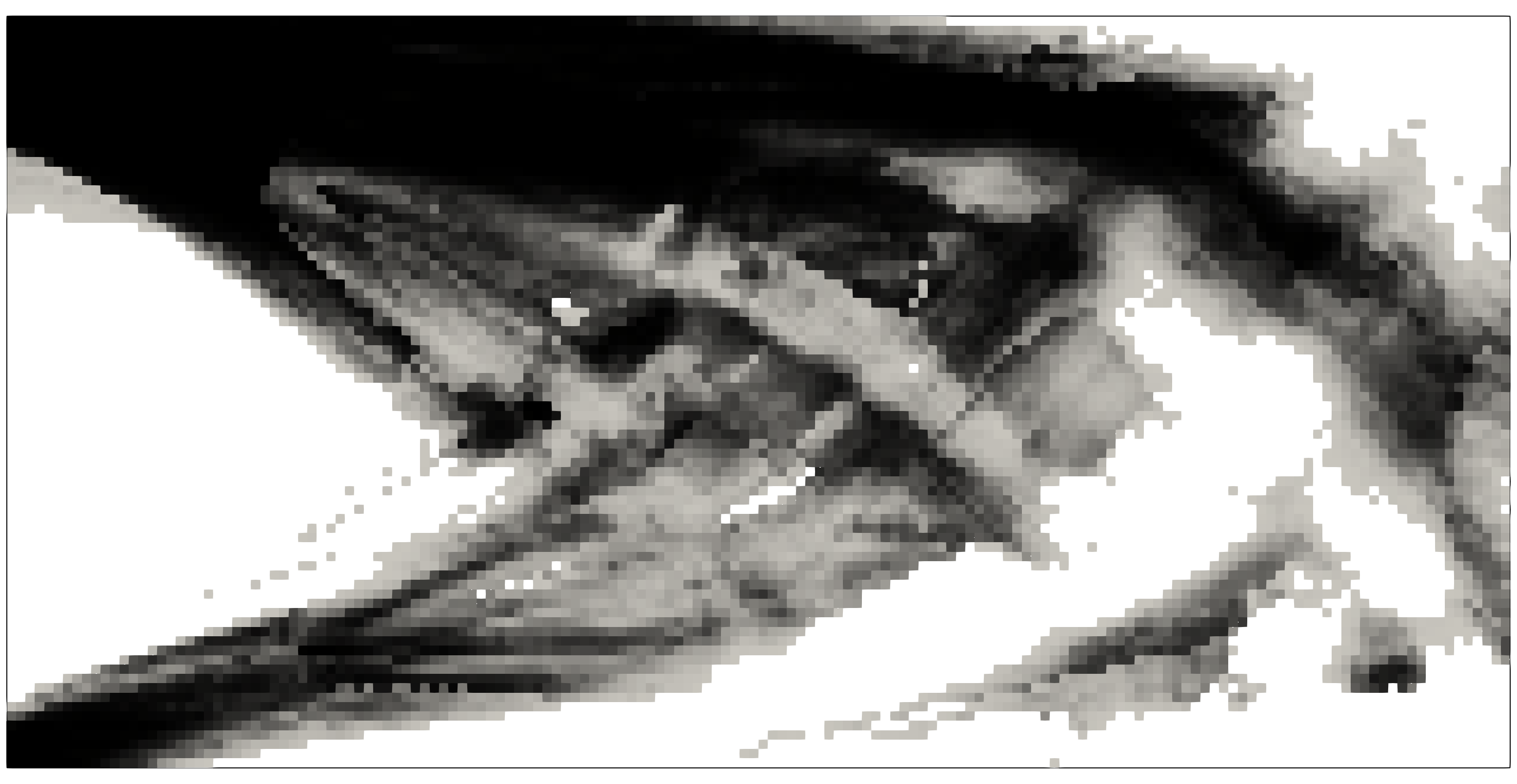}}
	\hspace{5pt}
	\subfigure[Case 8: $\bthetaEta_{\text{opt}}$]{\includegraphics[width=0.3\textwidth]{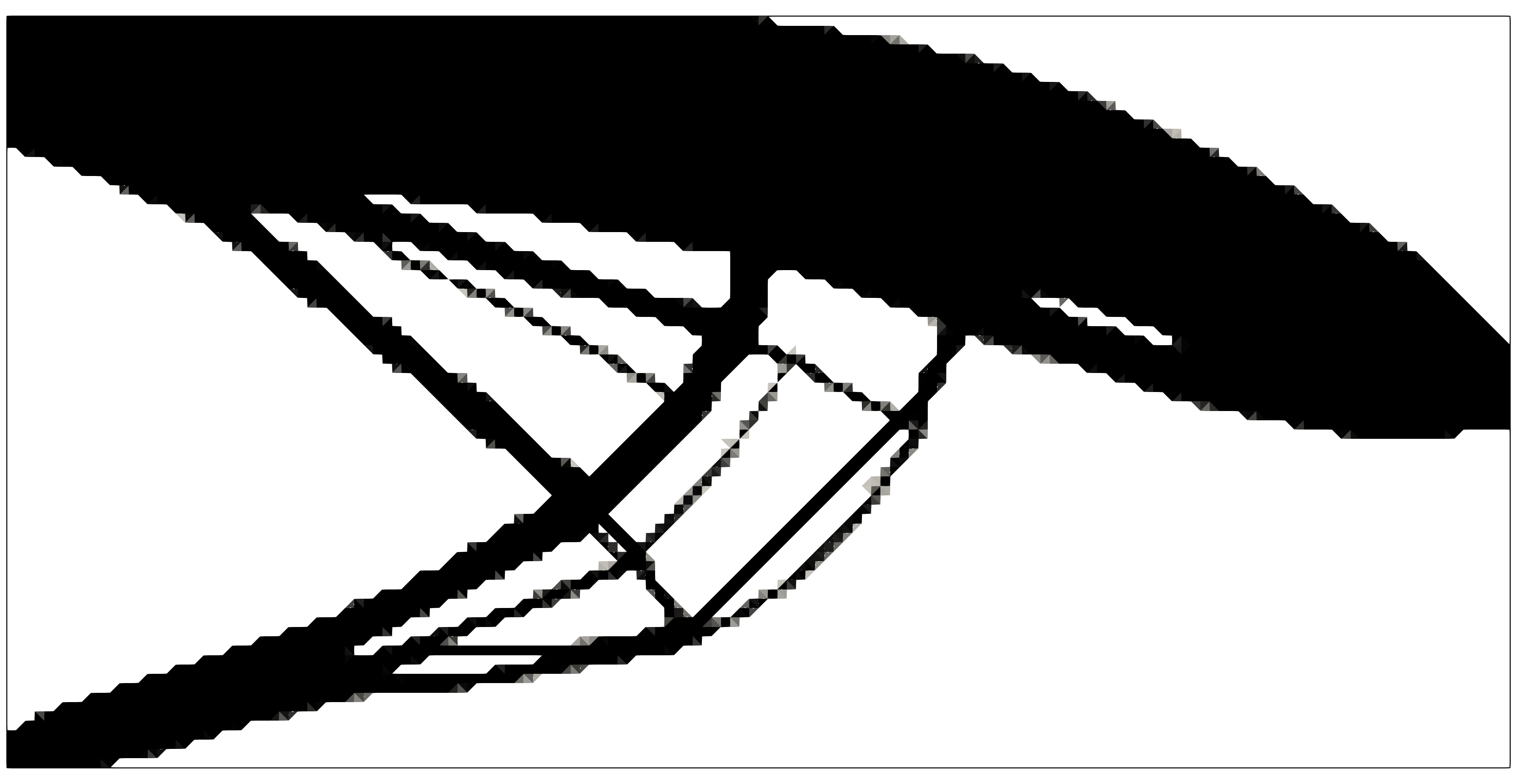}}
	
	\subfigure[Case 9: $\bthetaRef_{\text{opt}}$]{\includegraphics[width=0.3\textwidth]{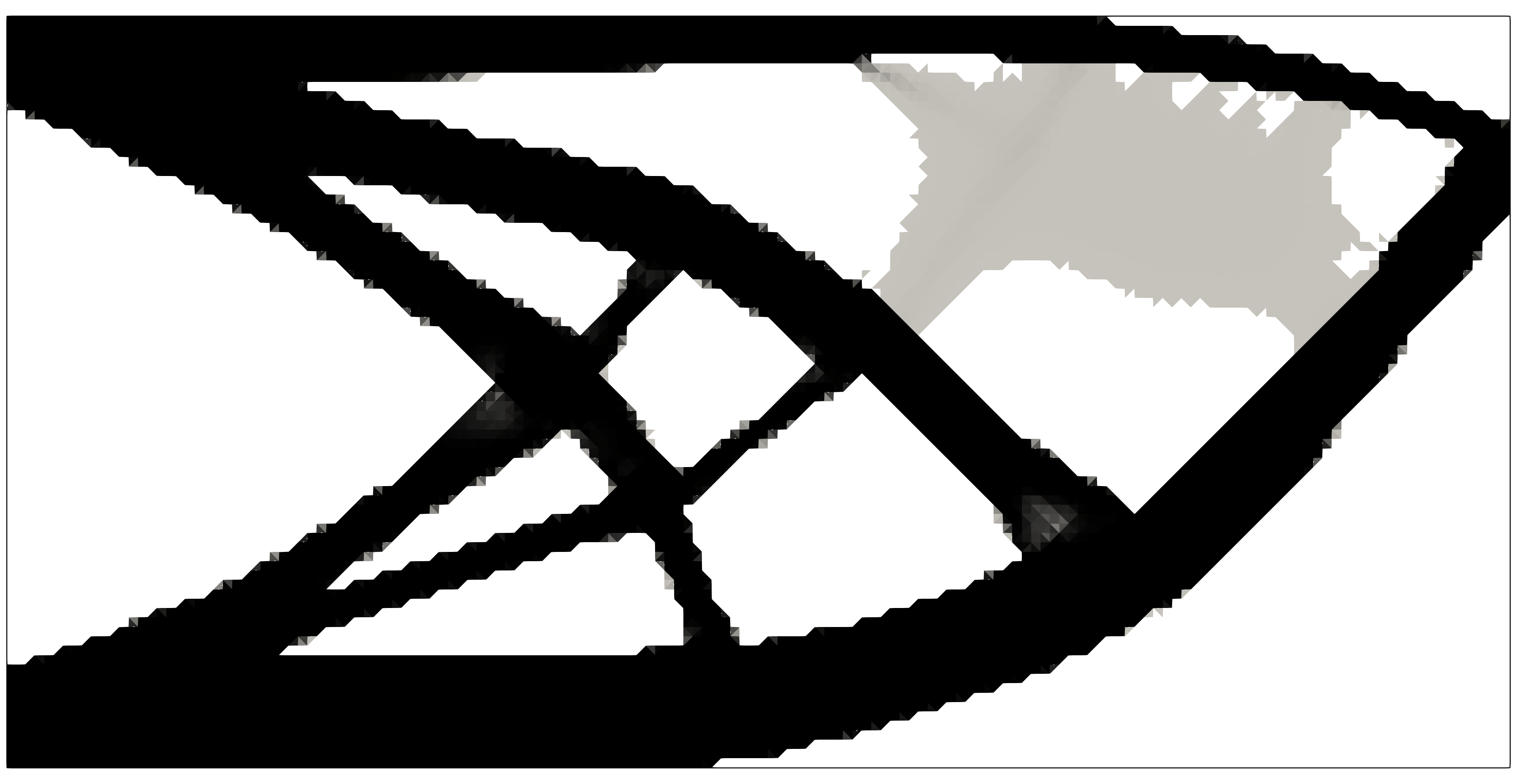}}
	\hspace{5pt}
	\subfigure[Case 9: $\bthetaEta$]{\includegraphics[width=0.3\textwidth]{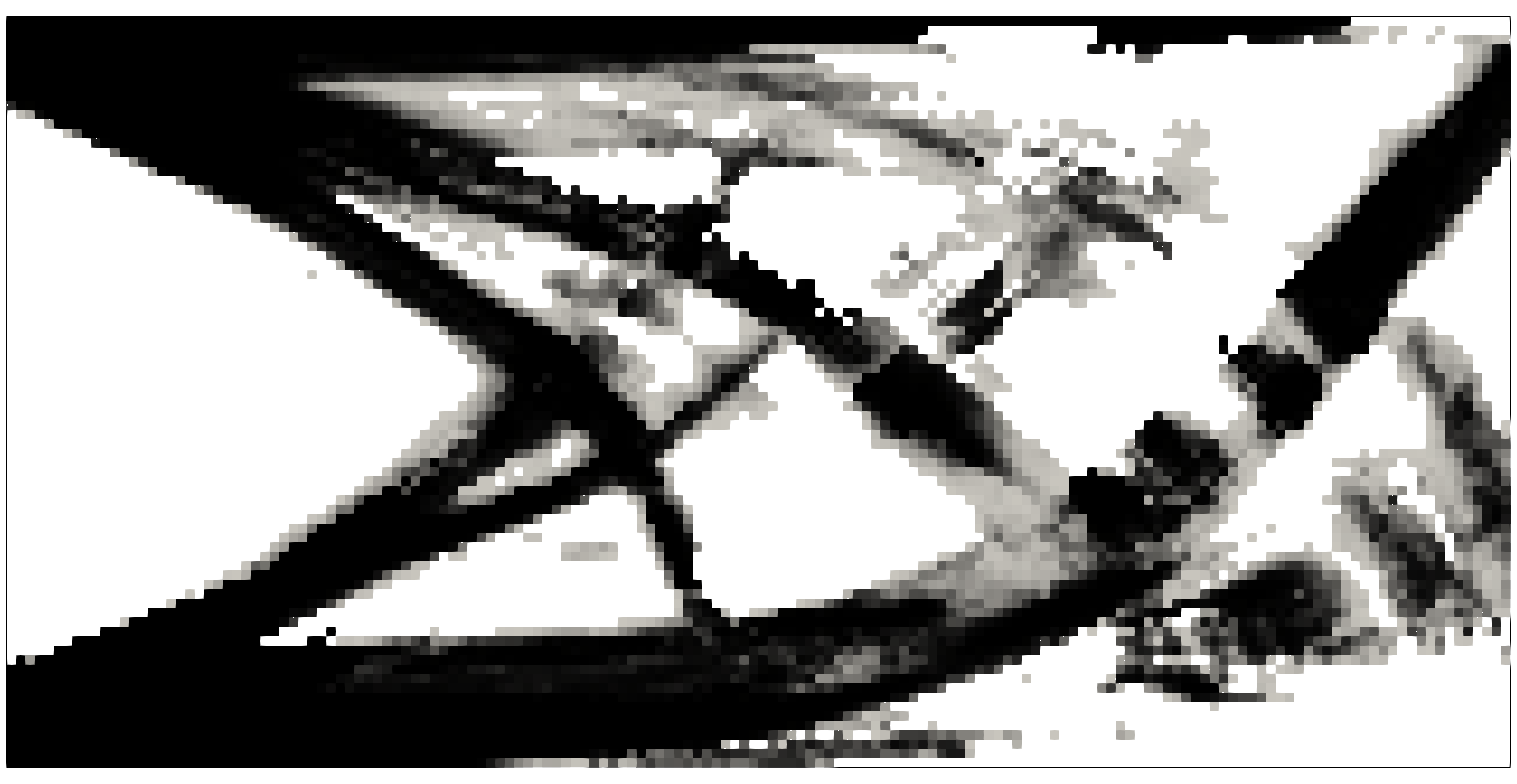}}
	\hspace{5pt}
	\subfigure[Case 9: $\bthetaEta_{\text{opt}}$]{\includegraphics[width=0.3\textwidth]{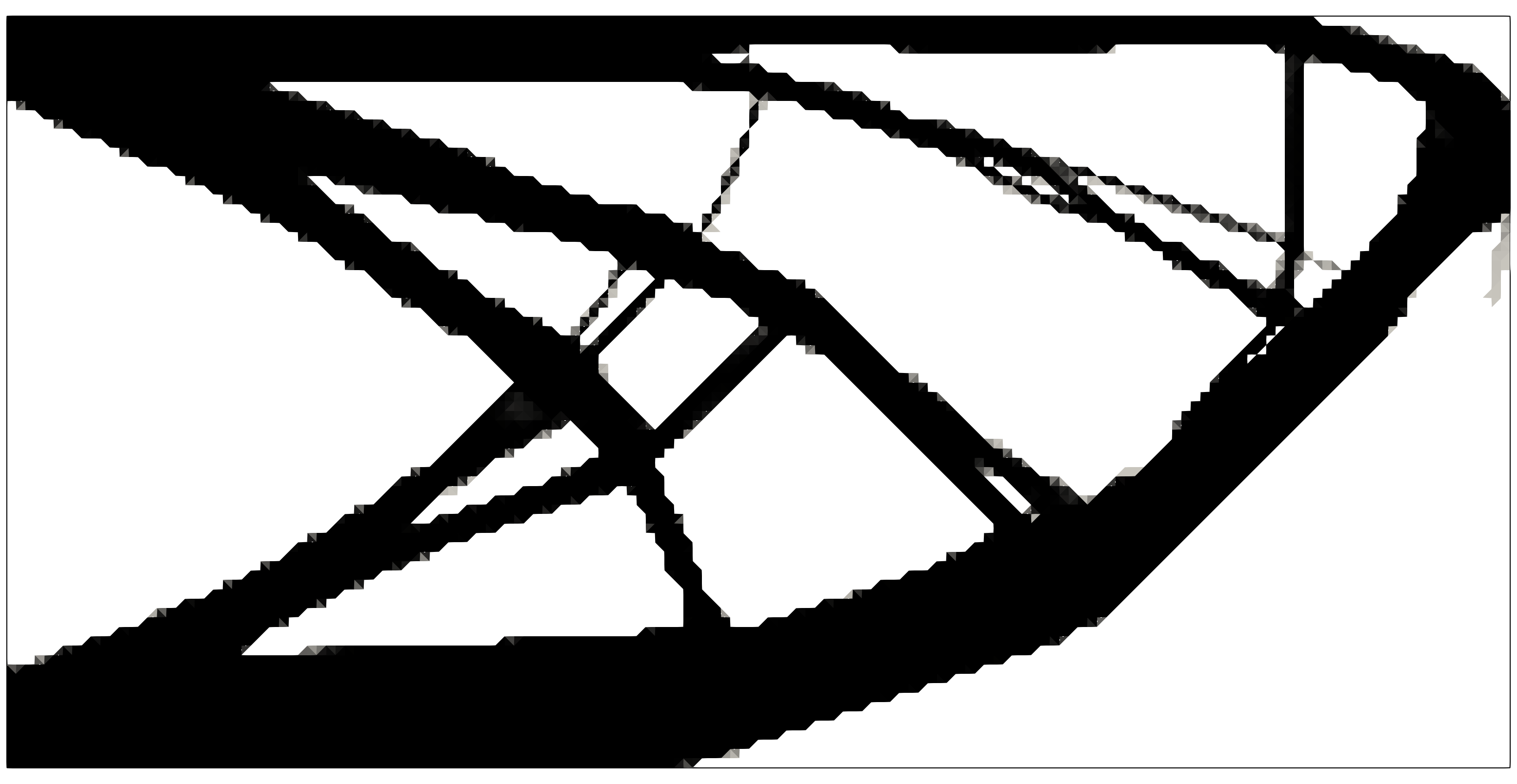}}
	
	\subfigure[Case 10: $\bthetaRef_{\text{opt}}$]{\includegraphics[width=0.3\textwidth]{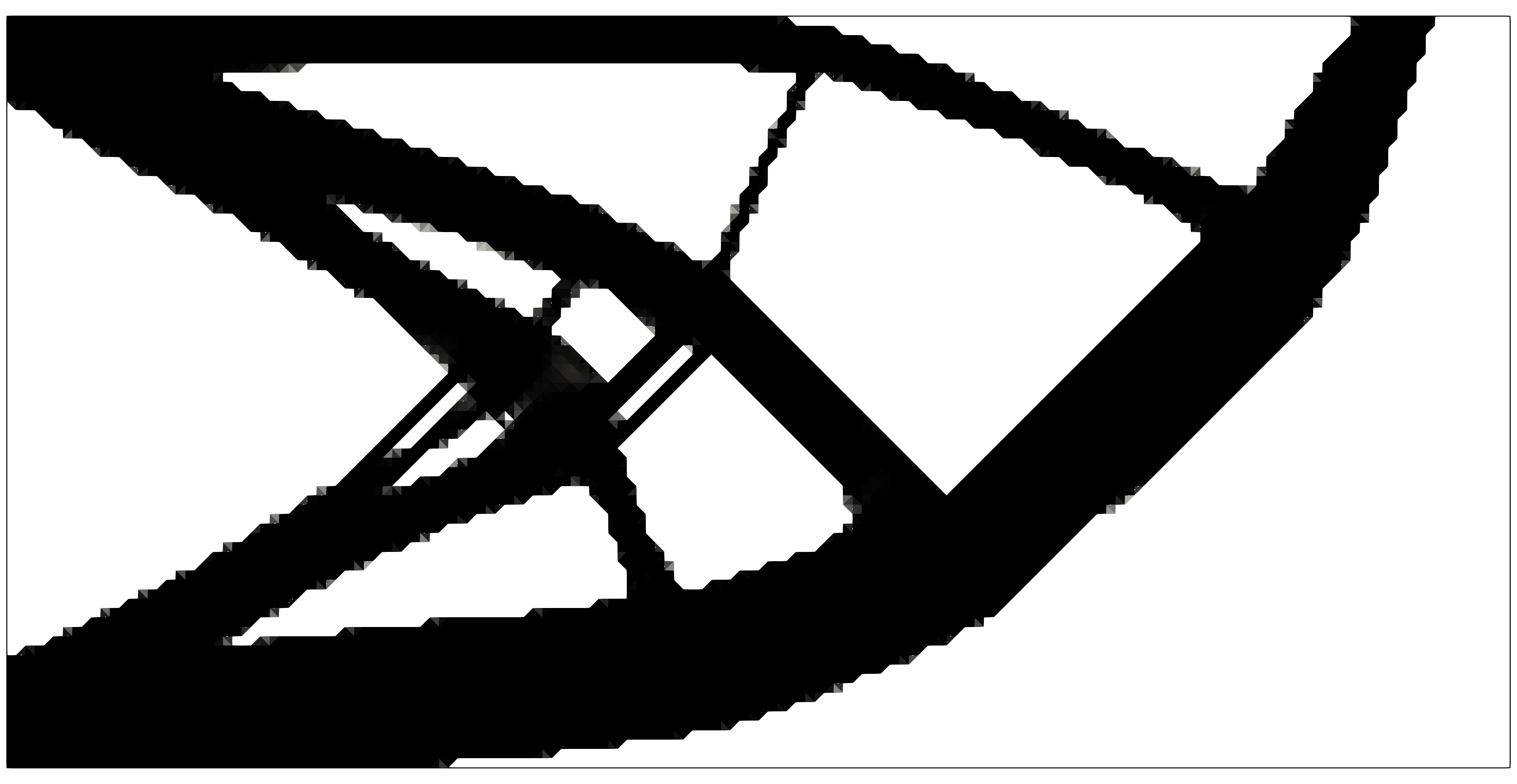}}
	\hspace{5pt}
	\subfigure[Case 10: $\bthetaEta$]{\includegraphics[width=0.3\textwidth]{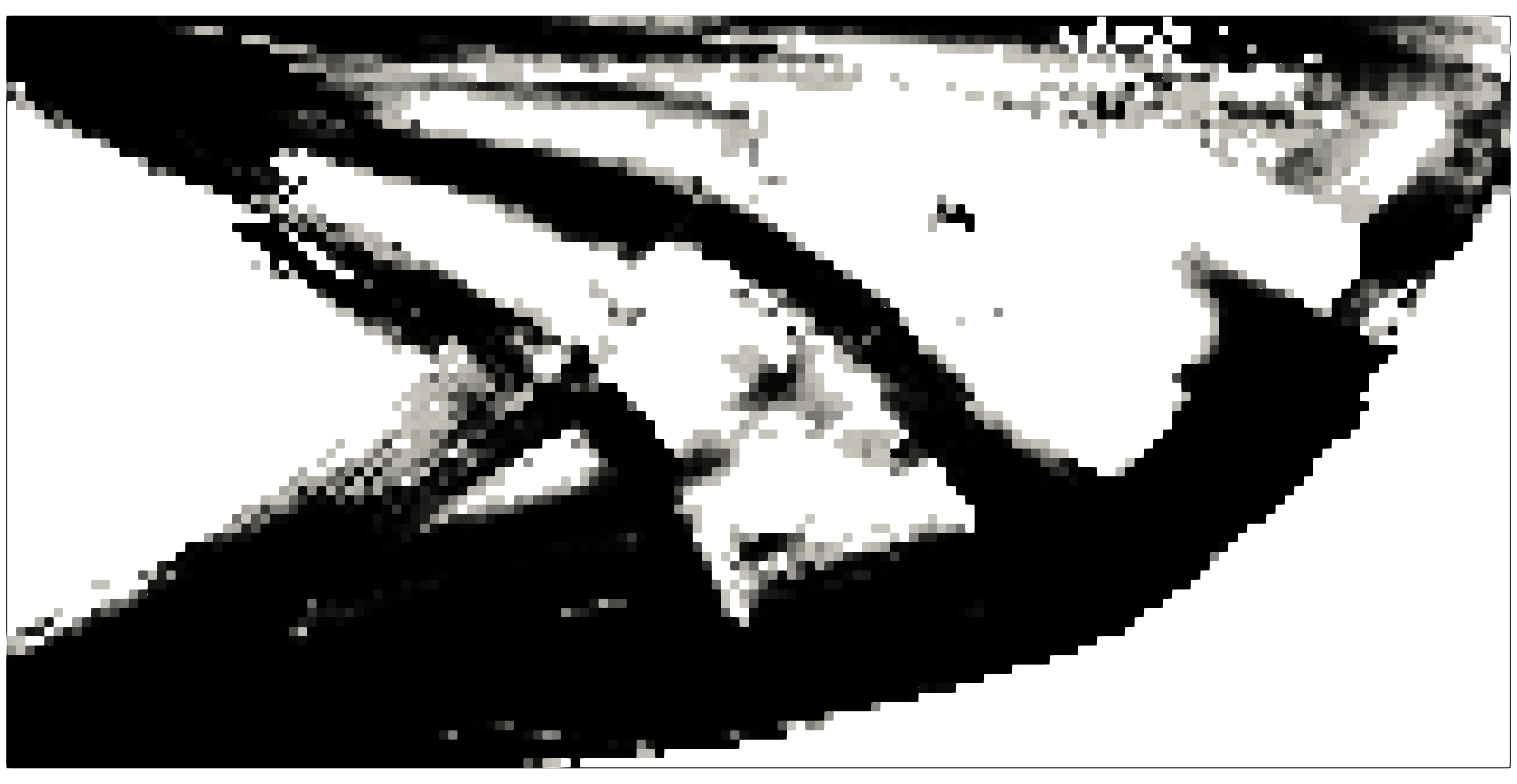}}
	\hspace{5pt}
	\subfigure[Case 10: $\bthetaEta_{\text{opt}}$]{\includegraphics[width=0.3\textwidth]{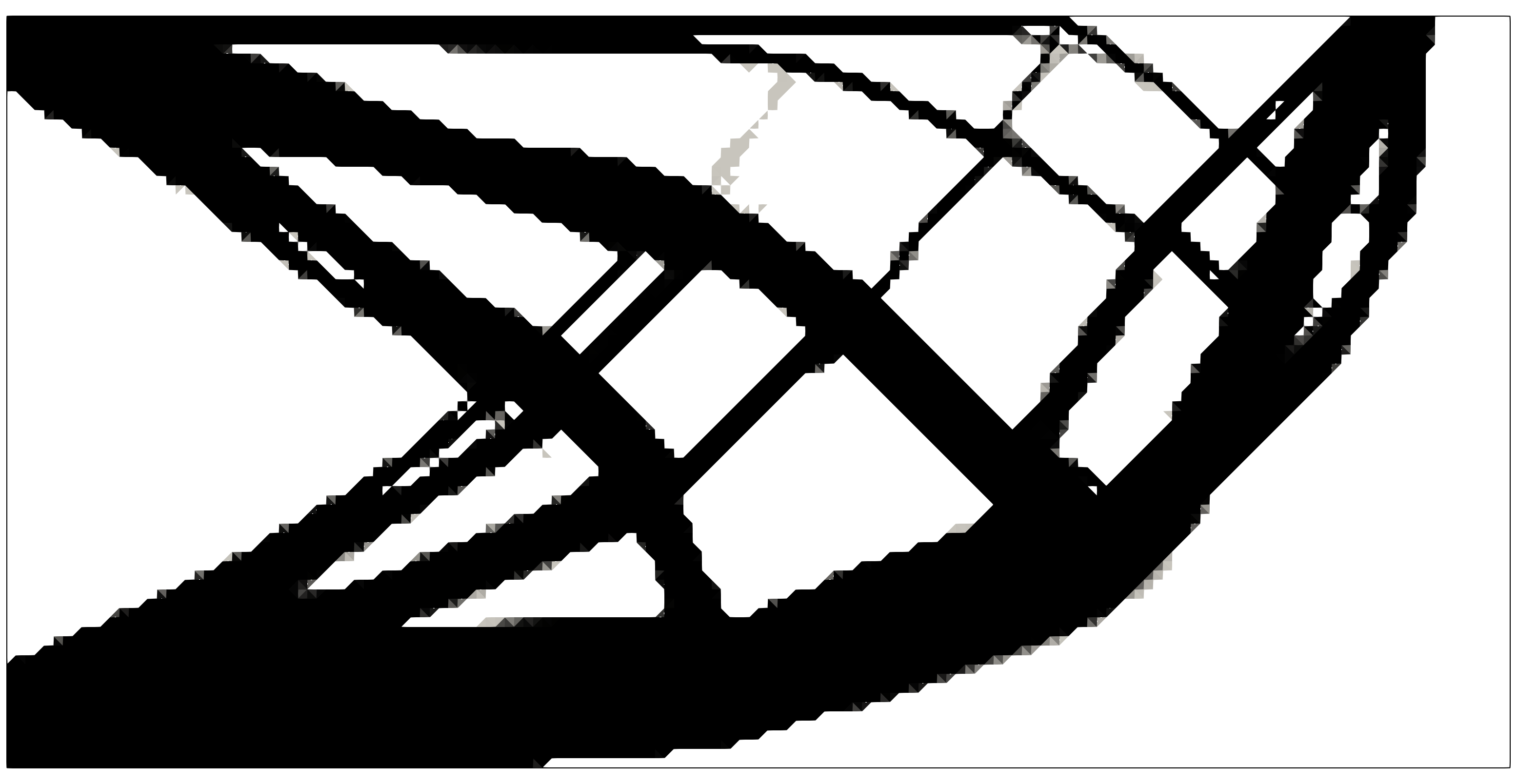}}
	
	\caption{Surrogate-based optimisation for extrapolated unseen cases. Left: \emph{ground truth} optimised topology. Centre: \emph{quasi-optimal} topology provided by the surrogate model. Right: topology optimised starting from the \emph{educated} initial guess in the centre.}
	\label{fig:ExtrapolSol}
\end{figure}

In Figure~\ref{fig:ExtrapolSol},  the \emph{quasi-optimal} topologies $\bthetaEta$ predicted by the $\texttt{FF}_{\!\eta}\texttt{-D}$ surrogate model (central column) and the outcome $\bthetaEta_{\text{opt}}$ of the surrogate-based optimisation (right column) are compared to the baseline material distributions $\bthetaRef_{\text{opt}}$ generated using the high-fidelity algorithm~\ref{alg:topOptHiFi} (left column).
The results confirm the ability of the surrogate model to extract the most relevant features of the structure including position of Dirichlet boundaries, force location, and some topological characteristics. Of course, the resolution achieved by $\texttt{FF}_{\!\eta}\texttt{-D}$ in the extrapolation setting is inferior with respect to the one displayed in Figure~\ref{fig:SeedSol} for the interpolation case. Nonetheless, when this is used as an \emph{educated} initial guess for the optimisation algorithm, the resulting topology converges to a mechanically-admissible structure mimicking the behaviour of the high-fidelity solution. 
It is worth noticing that the final topologies $\bthetaEta_{\text{opt}}$ tend to develop truss-like structures with finer details generated from the penalisation of the intermediate values of the material distribution. 

\begin{figure}[!htb]
	\centering
	\subfigure[High-fidelity compliance]{\includegraphics[width=0.48\textwidth]{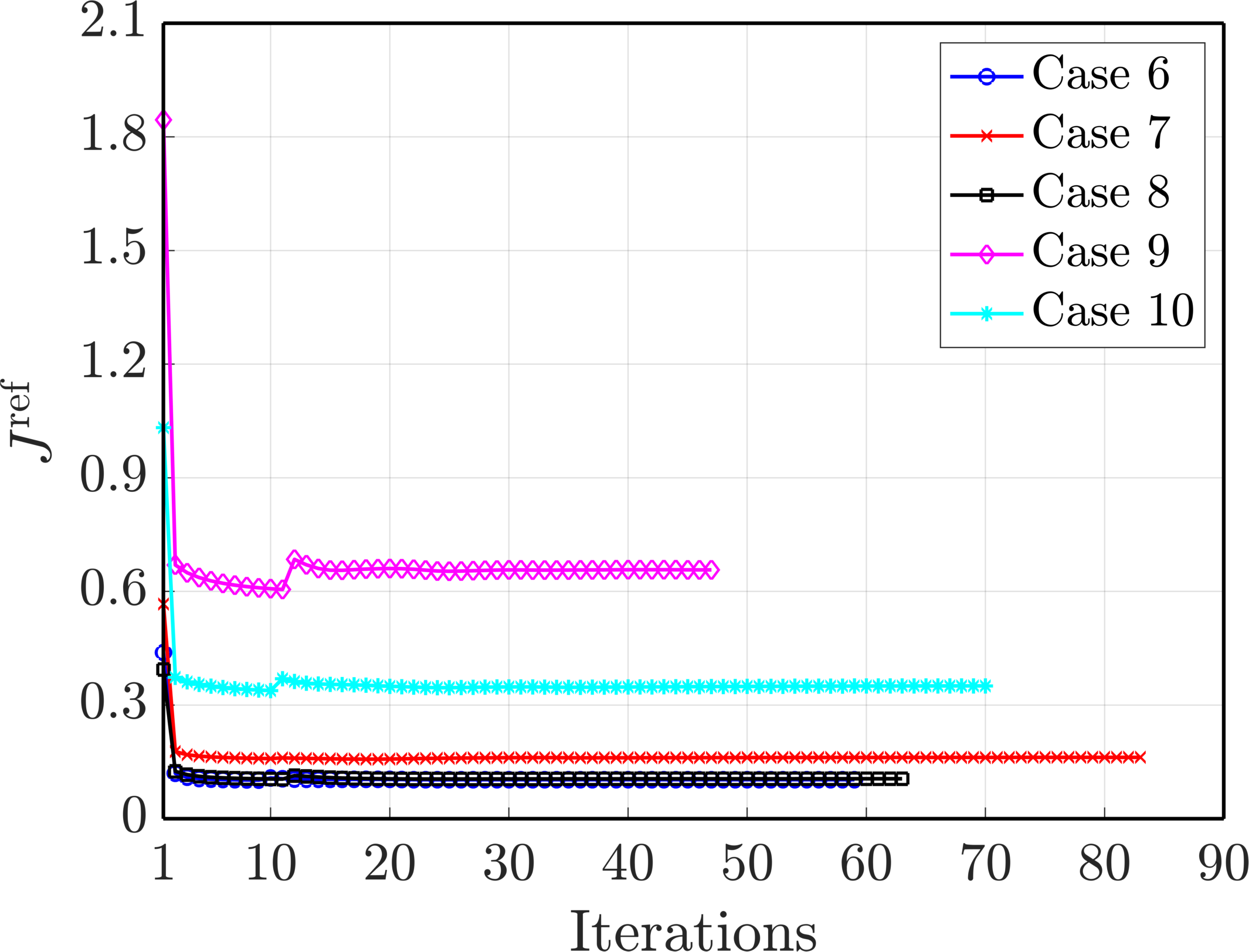}\label{fig:ExtrapolEvolutionA}}
	\hspace{5pt}
	\subfigure[Surrogate compliance]{\includegraphics[width=0.48\textwidth]{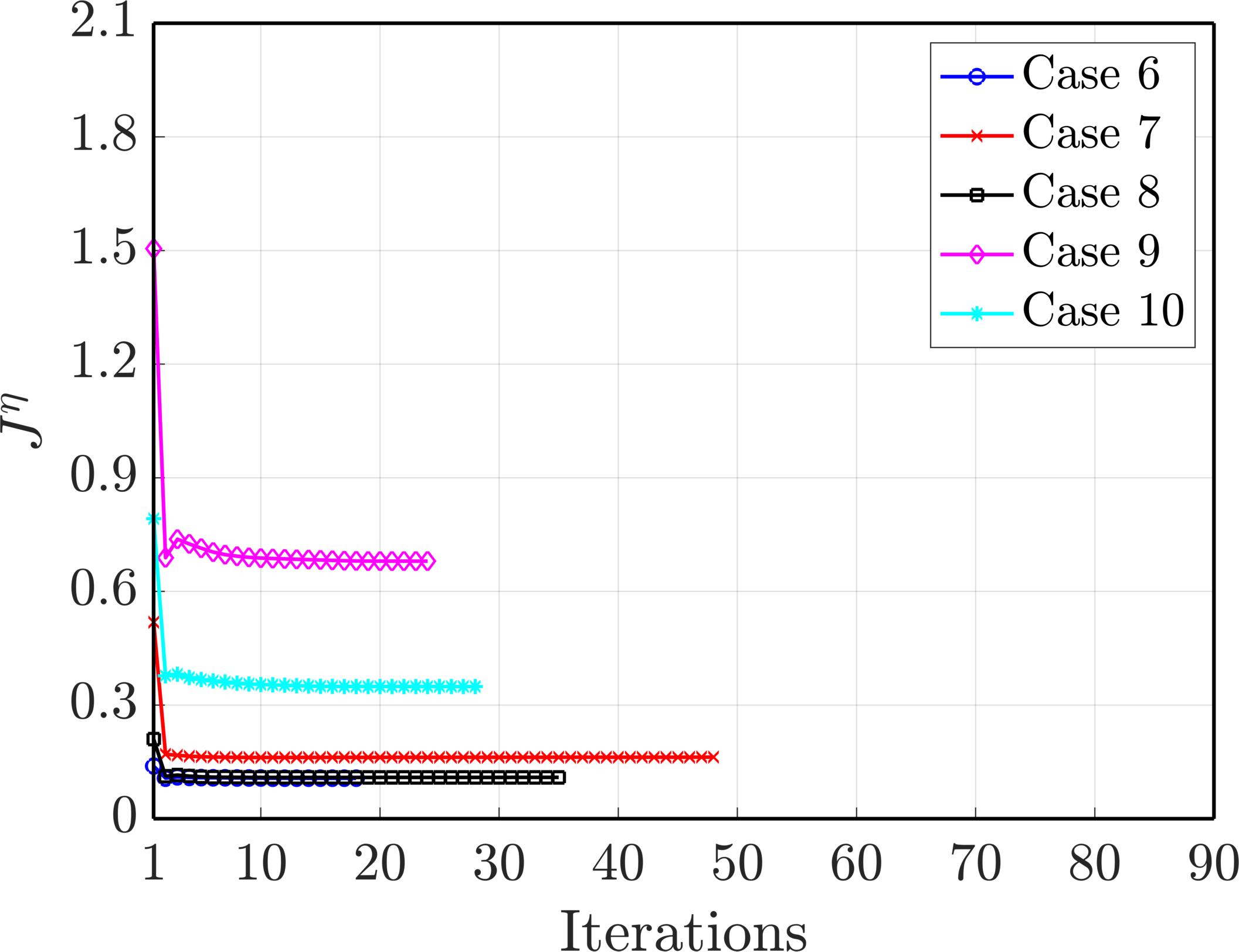}\label{fig:ExtrapolEvolutionB}}
	
	\subfigure[High-fidelity volume fraction]{\includegraphics[width=0.48\textwidth]{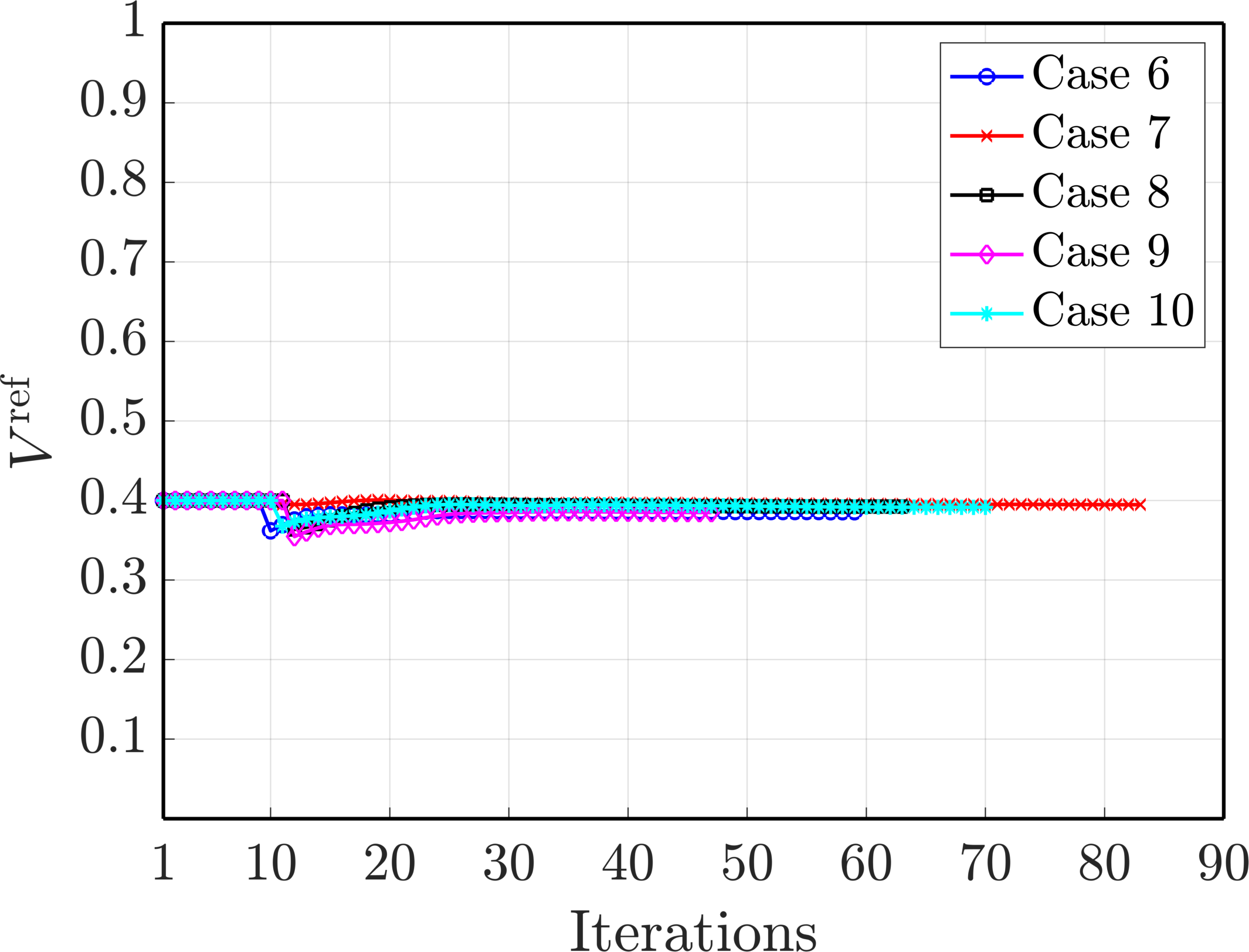}\label{fig:ExtrapolEvolutionC}}
	\hspace{5pt}
	\subfigure[Surrogate volume fraction]{\includegraphics[width=0.48\textwidth]{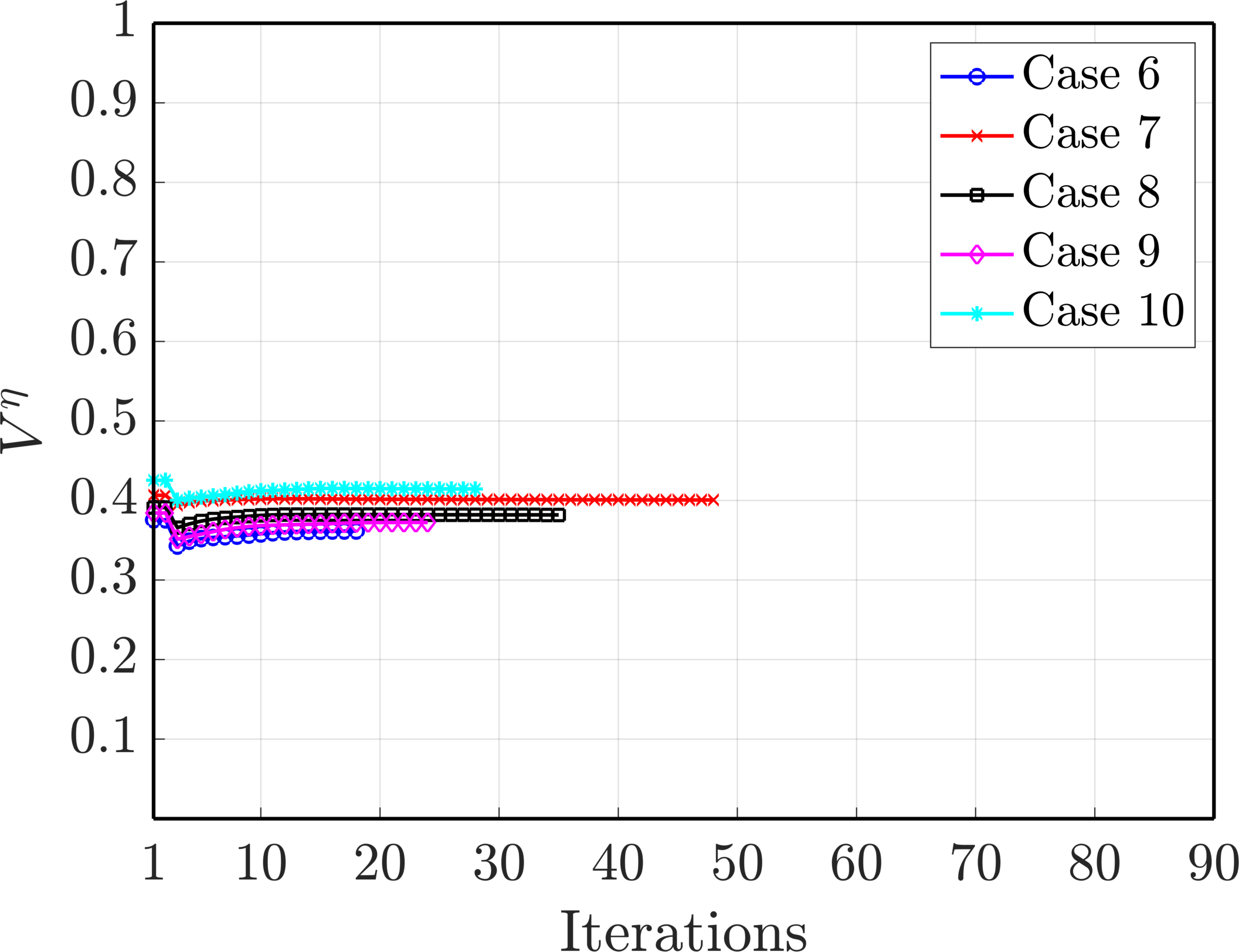}\label{fig:ExtrapolEvolutionD}}
	
	\caption{Evolution of the compliance (top) and the volume fraction (bottom) computed using the high-fidelity (left) and surrogate-based (right) algorithms for the extrapolated unseen cases. }
	\label{fig:ExtrapolEvolution}
\end{figure}

The computational gains observed in Section~\ref{sc:InterpolationOpti} are also confirmed in the case of extrapolation, when the surrogate model is trained without accessing data of an entire region of the parametric domain.
Despite the unseen test case being surrounded by many other unseen configurations, the surrogate model is still capable of learning relevant features of the structure and key information mapping the parameters $\bEta$ to the space of the optimised topologies.
Indeed, for all cases, the initial configuration provided by  $\texttt{FF}_{\!\eta}\texttt{-D}$ outperforms the initial material distribution of the high-fidelity optimiser in terms of the value of the compliance, see Figure~\ref{fig:ExtrapolEvolutionA} and~\ref{fig:ExtrapolEvolutionB}.
Moreover,  Table~\ref{tab:extrapolJV} reports the value $\Jeta[\text{opt}]$ of the objective functional upon convergence of the surrogate-based optimisation strategy, showcasing excellent agreement with the high-fidelity reference value $\Jref[\text{opt}]$, with mean error below $4\%$.
Finally, also in the extrapolation case, the use of the surrogate model $\texttt{FF}_{\!\eta}\texttt{-D}$ allows to reduce the computational cost of the optimisation pipeline, with an average reduction of the number of iterations by $53\%$, up to $69\%$ for case 6.
\begin{table}[!htb]
	\centering
	\begin{tabular}{| c || c | c | c || c | c | c || c | c |}
	\hline
	Case & $\Jeta[0]$ & $\Jeta[\text{opt}]$ & $\Jref[\text{opt}]$ & $\Veta[0]$ & $\Veta[\text{opt}]$ & $\Vref[\text{opt}]$ & $\niter[\eta]$ & $\niter[\text{ref}]$ \\
	\hline
	 6 & $0.139$ & $0.107$ & $0.103$ & $0.376$ & $0.361$ & $0.386$& $18$ & $59$\\
	\hline
	 7 & $0.519$ & $0.162$ & $0.162$ & $0.406$ & $0.400$ & $0.395$ & $48$ & $83$ \\
	\hline
	 8 & $0.210$ & $0.109$ & $0.105$ & $0.391$ & $0.382$ & $0.392$ & $35$ & $63$ \\
	\hline
	 9 & $1.505$ & $0.680$ & $0.657$ & $0.384$ & $0.372$ & $0.385$ & $24$ & $47$ \\
	\hline
	 10 & $0.792$ & $0.349$ & $0.351$ & $0.425$ & $0.414$ & $0.391$ & $28$ & $70$ \\
	\hline
	\end{tabular}

	\caption{Initial, optimal, and reference values of the compliance and the volume fraction computed using the high-fidelity and surrogate-based optimisation algorithms and corresponding number of iterations for the extrapolated unseen cases.}
	\label{tab:extrapolJV}
\end{table}

\section{Concluding remarks}
\label{sc:Conclusion}

In this work, a surrogate-based topology optimisation strategy was presented for linear elastic structures under parametric loads and boundary conditions.
The method goes beyond the state-of-the-art by proposing a novel paradigm for surrogate-based optimisation: instead of devising a surrogate model for the state (and adjoint, if needed) equations arising in the topology optimisation algorithm, this work proposes to \emph{reduce} the entire optimisation pipeline. 
This is achieved by training a surrogate model to learn the optimal topology for different values of the user-defined parameters.
The resulting \emph{quasi-optimal} topology is thus employed as an \emph{educated} initial guess for a novel optimisation algorithm, allowing to significantly reduce the computational cost of the overall optimisation procedure.

The method relies on encoder/decoder blocks to reduce the dimensionality of the parametric topology optimisation problem, while learning the mapping between the input parameters characterising the system and the low-dimensional latent space of the optimised topologies.
Different architectures of neural networks are proposed and numerically compared in terms of the accuracy of the prediction of unseen cases and the computational gains provided by the resulting surrogate-based optimisation algorithm.
Among the different proposed surrogate models, the architecture $\texttt{FF}_\eta\texttt{-D}$ provides the best trade-off between accuracy and efficiency. 
$\texttt{FF}_\eta\texttt{-D}$ combines a fully-connected feed-forward net $\texttt{FF}_\eta$ mapping the input parameters $\bEta$ to a latent space of dimension $\nA$ with a block $\texttt{D}$, mimicking a decoder, to upscale the information and reconstruct the topology in a space of dimension $\nT \gg \nA$.
The surrogate model concurrently learns the parametric dependence and achieves a compression factor of more than $700$, with approximately half of the unknowns required by a classic autoencoder using a decoder of the same size.

The resulting surrogate model is tested by means of cross-validation, showing an average relative mean squared error in the prediction of the \emph{quasi-optimal} topologies of $11\%$.  When used to initialise the surrogate-based topology optimisation pipeline, $\texttt{FF}_\eta\texttt{-D}$ attains accuracies in the estimated compliance between $1\%$ and $7\%$, with a reduction of the required number of iterations up to $64\%$.
The approximation and generalisation capabilities of $\texttt{FF}_\eta\texttt{-D}$ are further tested by considering the challenging scenario of extrapolation, that is,  the surrogate model is trained and validated without having access to the information of an entire portion of the parametric space which is then used for testing.
In this case, the accuracy of the \emph{quasi-optimal} predictions clearly worsens, with the average rMSE growing from $11\%$ to $29\%$, suffering similar difficulties as classic autoencoders (average error grows from $6\%$ to $21\%$). Nonetheless, the \emph{quasi-optimal} topologies still allow to improve the performance of high-fidelity optimisers by reducing the average  number of optimisation iterations by $53\%$ while achieving discrepancies below $4\%$ in the optimal value of the objective functional.

\hl{
To summarise, this study provides a \emph{demonstrator} of the feasibility for non-intrusive surrogate models to leverage existing data of optimised topologies and learn relevant features of previously optimised designs.
By training the surrogate model on a dataset of optimised topologies, only the final design configurations need to be saved, significantly reducing storage requirements, and eliminating the need to access the physical solutions of the state and adjoint equations at each optimisation iteration.
Nonetheless, consistency with the underlying physics is guaranteed by using this purely data-driven \emph{quasi-optimal} topology as initial condition of a surrogate-based optimisation procedure: in this algorithm, governing equations are enforced while the overall number of optimisation iterations is significantly reduced.
The work showcases a proof-of-concept of the surrogate model using a 4-dimensional benchmark problem in linear elasticity, namely, the design of a 2D cantilever beam with minimum compliance and maximum volume constraint, with a parametrised external load.
In order to fully evaluate the computational advantage of the proposed strategy to couple surrogate models with high-fidelity optimisers, future works need to take into account more realistic scenarios, such as three-dimensional systems, nonlinear models, and larger numbers of parameters.
In this context, the proposed surrogate model strategy is expected to benefit from contributions on ML architectures capable of handling scarce data~\cite{Choi-CCSK-24}, adaptive sampling algorithms specifically balancing exploitation and exploration of high-dimensional parametric spaces~\cite{Liu-LPC-18}, data augmentation procedures to engineer parsimonious artificial dataset entries~\cite{Muixi-MZGD-25}, while possibly leveraging data of different nature according to multi-fidelity paradigms~\cite{Park-PHK-17-reviewMF,Gunzburger-PWG-18}.
}

\section*{Acknowledgements}
The authors would like to express their sincere gratitude to Dr. Guillem Barroso, whose contributions were instrumental in initiating and shaping this work. His efforts during the early stages of the project were essential to its development.

This work was supported by MCIU/AEI/10.13039/501100011033,  Spanish Ministry of Science, Innovation and Universities and Spanish State Research Agency  (Grants No. TED2021-132021B-I00, PID2023-149979OB-I00) and by the Generalitat de Catalunya (Grant No. 2021-SGR-01049).
MG is Fellow of the Serra H\'unter Programme of the Generalitat de Catalunya.

\bibliographystyle{ieeetr}
\bibliography{Ref-surrogateTopOpt}

\appendix

\section{Computation of the homogenised elasticity tensor}
\label{sc:appHooke}

Using Voigt notation~\cite{FishBelytschko-book-07}, $\nsd \times \nsd$ symmetric second-order tensors are stored in vector form accounting only for the $\msd := \nsd(\nsd+1)/2$ non-redundant components.
Similarly, symmetric fourth-order tensors can be written in the form of $\msd \times \msd$ matrices, instead of storing $\nsd^4$ entries.

For the two-dimensional case under analysis, the strain tensor in Voigt notation is thus defined as 
\begin{equation} \label{eq:strainVoigt}
\strainV := \bigl[\varepsilon_{11} ,\; \varepsilon_{22} ,\; \varepsilon_{12} \bigr]^\top ,
\end{equation}
with components
\begin{equation} \label{eq:NormalShear}
\varepsilon_{ij} := \frac{\partial u_i}{\partial x_j} + (1-\delta_{ij}) \frac{\partial u_j}{\partial x_i}, \quad \text{for } i,j = 1,2 ,
\end{equation}
$\delta_{ij}$ being the Kronecker delta.
Moreover, given the Voigt formulation of the fourth-order elasticity tensor
\begin{equation}\label{eq:compA}
\matA :=
\begin{bmatrix}
2\mu+\lambda & \lambda & 0 \\
\lambda & 2\mu+\lambda & 0 \\
0 & 0 & \mu
\end{bmatrix} ,
\end{equation}
the resulting expression for the stress tensor expressed is given by $\stressV = \matA \strainV$.

Similarly, the elasticity tensor $\matAhs$ fulfilling the Hashin-Shtrikman bounds can be written as
\begin{equation}\label{eq:compAhs}
\matAhs :=
\begin{bmatrix}
2\muE+\lambdaE & \lambdaE & 0 \\
\lambdaE & 2\muE+\lambdaE & 0 \\
0 & 0 & \muE
\end{bmatrix} ,
\end{equation}
with the effective moduli
\begin{subequations}
\begin{align}
\lambdaE &:= \frac{\thetaBar \, \mu(\mu + \lambda)\Big(\lambda + 2(1-\thetaBar)\mu\Big)}{\Big(\mu + \lambda + (1-\thetaBar)(3\mu+\lambda)\Big)\Big(\mu + (1-\thetaBar)(\mu+\lambda)\Big)} ,  \\
\muE &:= \frac{\thetaBar \, \mu(\mu + \lambda)}{\mu + \lambda + (1-\thetaBar)(3\mu+\lambda)} ,
\end{align}
\end{subequations}
and $\thetaBar$ and $1-\thetaBar$ being the volume fractions of the two phases.

Following the same rationale, the components of the fourth-order elasticity tensor $\matAc{i}$ defined in~\eqref{eq:Ac} for the $i$-th phase of the sequential laminate are given by
\begin{equation}\label{eq:compAc}
\begin{aligned}
[\matAc{i}]_{11} :=& 2\mu+\lambda - \frac{1}{\mu}\Big((2\mu+\lambda)^2[\vSig{i}]_1^2+\lambda^2[\vSig{i}]_2^2\Big) \\
&+ \frac{\mu+\lambda}{\mu(2\mu+\lambda)}\Big((2\mu+\lambda)[\vSig{i}]_1^2+\lambda[\vSig{i}]_2^2\Big)^2 , \\
[\matAc{i}]_{22} :=& 2\mu+\lambda - \frac{1}{\mu}\Big(\lambda^2[\vSig{i}]_1^2+(2\mu+\lambda)^2[\vSig{i}]_2^2\Big) \\
&+ \frac{\mu+\lambda}{\mu(2\mu+\lambda)}\Big(\lambda[\vSig{i}]_1^2+(2\mu+\lambda)[\vSig{i}]_2^2\Big)^2 , \\
[\matAc{i}]_{33} :=& \mu - \frac{1}{\mu}\mu^2 + \frac{\mu+\lambda}{\mu(2\mu+\lambda)}\Big(2\mu[\vSig{i}]_1[\vSig{i}]_2\Big)^2 , \\
[\matAc{i}]_{12} :=& \lambda - \frac{1}{\mu}\lambda(2\mu+\lambda) \\
&+ \frac{\mu+\lambda}{\mu(2\mu+\lambda)}\Big(\lambda[\vSig{i}]_1^2+(2\mu+\lambda)[\vSig{i}]_2^2\Big)\Big((2\mu+\lambda)[\vSig{i}]_1^2+\lambda[\vSig{i}]_2^2\Big) , \\
[\matAc{i}]_{23} :=& -\frac{1}{\mu}(\mu+\lambda)\Big(2\mu[\vSig{i}]_1[\vSig{i}]_2\Big) \\
&+ \frac{\mu+\lambda}{\mu(2\mu+\lambda)}\Big(\lambda[\vSig{i}]_1^2+(2\mu+\lambda)[\vSig{i}]_2^2\Big)\Big(2\mu[\vSig{i}]_1[\vSig{i}]_2\Big) , \\
[\matAc{i}]_{13} :=& -\frac{1}{\mu}(\mu+\lambda)\Big(2\mu[\vSig{i}]_1[\vSig{i}]_2\Big) \\
&+ \frac{\mu+\lambda}{\mu(2\mu+\lambda)}\Big((2\mu+\lambda)[\vSig{i}]_1^2+\lambda[\vSig{i}]_2^2\Big)\Big(2\mu[\vSig{i}]_1[\vSig{i}]_2\Big) .
\end{aligned}
\end{equation}

The homogenised compliance tensor is obtained from equation~\eqref{eq:optA} by the linear combination of the inverse of $\matA$ and the inverse of $[m_1 \matAc{1} + m_2 \matAc{2}]$, respectively with weights $1$ and $(1-\theta)/\theta$.
In order to solve the elastic problem~\eqref{eq:elastWeak}, the corresponding homogenised elasticity tensor is computed by inverting $[\matAs]^{-1}$.
Note that these operations are straightforward by recalling that, for an invertible matrix $\matC$, it holds $\matC^{-1} = \adj(\matC)/\Det(\matC)$, with $\adj$ and $\Det$ denoting the adjoint and the determinant of the matrix, respectively.

\section{Architecture of the networks}
\label{sc:appNets}

This appendix reports some technical details on the architectures and settings employed for the surrogate models based on parametric autoencoders presented in this work.

The encoder block $\texttt{E}$ has an input layer of dimension $12,800$ and three inner layers featuring $200$, $100$,  and $25$ neurons, respectively.
The decoder block $\texttt{D}$ consists of three inner layers of dimension $25$, $100$, and $200$ and an output layer with $12,800$ neurons.
The maximum dimension of the latent space is set to $\nA^{\text{max}} = 25$.
The feed-forward net $\texttt{FF}_{\!\eta}$ mapping the parameters $\bEta$ onto the latent space features an input layer of dimension $2$ and two layers with $50$ and $25$ neurons.

ReLU activation functions are employed in all hidden layers, whereas a sigmoid function is used in the output layer.
The parameters of the models are initialised by means of the Kaiming uniform method.

The regularisation coefficients are set according to the values in Table~\ref{tab:coeff}.
\begin{table}[!htb]
	\centering
	\begin{tabular}{| l || c | c | c | c | c |}
	\hline
 	& $\texttt{FF}_{\!\eta}\texttt{-D}$ & $\texttt{E-FF}_{\!\eta}\texttt{-D}$ & $\texttt{S-AE-FF}_{\!\eta}\texttt{-D}$ & $\texttt{S-AE-FF}_{\!\eta}$ & $\texttt{AE}$ \\
	\hline
	Loss & \eqref{eq:lossParametric} & \eqref{eq:lossCombined} & \eqref{eq:lossStagWithD} & \eqref{eq:lossStagComb} & \eqref{eq:LAE} \\
	\hline
	$\NNpar{\alpha}$ & - & $10^{-3}$ & - & - & - \\
	\hline
	$\NNpar{r}$ & $10^{-4}$ & $10^{-4}$ & - & - & $10^{-4}$ \\
	\hline
	$\NNpar{r,\mathrm{I}}$ & - & - & $10^{-4}$ & $10^{-4}$ & - \\
	\hline
	$\NNpar{r,\mathrm{II}}$ & - & - & $10^{-4}$ & $10$ & - \\
	\hline
	\end{tabular}
	
	\caption{Loss functions and regularisation coefficients for each model.}
	\label{tab:coeff}
\end{table}

Training of the different architectures presented in Section~\ref{sc:Surrogate} is performed during $5,000$ epochs using the Adam optimiser with batches of size $600$ and learning rate $10^{-3}$.  
Early stopping tolerance is set to $10^{-3}$, early stopping patience to $1,000$ epochs and no dropout is applied.

\end{document}